\documentclass[a4paper,12pt]{amsart}

\usepackage{iftex}
\ifPDFTeX
  \usepackage[T1]{fontenc}
  \usepackage{lmodern}
\else
  \usepackage{fontspec}
\fi
\usepackage[english]{babel}
\usepackage{microtype}

\usepackage[margin=25mm]{geometry}
\usepackage{setspace}

\usepackage{tikz}
\usetikzlibrary{arrows.meta}

\makeatletter
\renewcommand{\section}{%
  \@startsection{section}{1}%
    \z@{.7\linespacing\@plus\linespacing}{.5\linespacing}%
    {\normalfont\large\bfseries\centering}}
\makeatother

\usepackage{amsmath,amssymb,amsthm}
\numberwithin{equation}{section}
\newtheorem{theorem}{Theorem}[section]
\newtheorem{lemma}[theorem]{Lemma}
\newtheorem{proposition}[theorem]{Proposition}
\newtheorem{corollary}[theorem]{Corollary}

\newtheorem*{restatedwaist}{Theorem \ref{theorem: new waist inequality}}
\newtheorem*{restatedavoidance}
  {Proposition \ref{proposition: quantitative avoidance}}

\newtheorem*{restatedapproximation}
  {Proposition \ref{proposition: low dimensional approximation}}

\theoremstyle{definition}

\theoremstyle{remark}
\newtheorem{remark}[theorem]{Remark}
\newtheorem{problem}[theorem]{Problem}

\usepackage{graphicx}
\usepackage{booktabs}
\usepackage{array}
\usepackage{tabularx}
\usepackage{caption}
\usepackage{subcaption}

\usepackage{tikz-cd}

\usepackage[numbers,sort&compress]{natbib}

\usepackage[hidelinks]{hyperref}

\newcommand{\mdim}{\mathrm{mdim}}
\newcommand{\diam}{\mathrm{diam}}
\newcommand{\widim}{\mathrm{Widim}}

\newcommand{\supp}{\mathrm{supp}}
\newcommand{\norm}[1]{\left|\!\left|#1\right|\!\right|}
\newcommand{\lmdimM}{\underline{\mathrm{mdim}}_\mathrm{M}}

\title[Waist inequality and mean dimension]
{Waist inequality and a 
Hurewicz-type theorem for mean dimension}

\author{Qiang Huo, Kangbo Ouyang,
Ruxi Shi, Masaki Tsukamoto} 

\date{} 

\address[Qiang Huo]
{School of Mathematical Sciences and
Key Laboratory of Mathematics and Information Networks 
(Ministry of Education), 
Beijing University of Posts and Telecommunications, Beijing 100876, China}

\email{qianghuo@ustc.edu.cn}

\address[Kangbo Ouyang]
{School of Mathematical Sciences, 
University of Science and Technology of China, 
Hefei 230026, China}

\email{oy19981231@mail.ustc.edu.cn}

\address[Ruxi Shi]
{Shanghai Center for Mathematical Sciences, 
Fudan University, Shanghai 200438, China}

\email{ruxishi@fudan.edu.cn}

\address[Masaki Tsukamoto]
{Department of Mathematics, 
Kyoto University, Kyoto 606-8502, Japan}

\email{tsukamoto@math.kyoto-u.ac.jp}



\thanks{
Q.H. was partially supported by National Key Research 
and Development Program of China 2024YFA1013600.
K.O. was partially supported 
by the National Natural Science Foundation of China 
(Grant Nos. 12371196 and 12426201).
R.S. was supported by the NSFC Nos.~12571198 
and 12231013.
M.T. was supported by JSPS KAKENHI JP25K06974.}

\begin{document}

\begin{abstract}
We establish a Hurewicz-type inequality 
relating mean dimension
and conditional 
metric mean dimension.
For every equivariant continuous map 
between dynamical systems, 
the mean dimension of the domain is bounded
above by the sum of the mean dimension 
of the target and the
conditional metric mean dimension of the map.
In particular, 
any equivariant continuous map 
that decreases mean dimension 
must have infinite topological conditional entropy.
The key ingredient is a new waist inequality, 
which we believe is of independent geometric interest 
beyond its application to dynamical systems.
Compared with Gromov's original waist inequality,
a distinctive feature of this inequality
is that it applies to finite simplicial complexes 
of arbitrary topology, 
subject only to dimensional restrictions.
\end{abstract}

\maketitle

\section{Introduction} \label{section: introduction}

\subsection{Entropy and mean dimension} 
\label{subsection: entropy and mean dimension}

The main purpose of this paper 
is to establish a Hurewicz-type inequality 
relating mean dimension and 
conditional metric mean dimension 
for dynamical systems. 
In the course of this study, 
we also establish a new waist inequality, 
which may be of independent interest.

We begin by explaining a main consequence of 
our Hurewicz-type inequality concerning entropy 
and mean dimension of dynamical systems. 
We present this consequence first 
because it requires little preparation and 
gives a particularly clean formulation 
of the main phenomenon studied in this paper.
The Hurewicz-type inequality itself will be stated in 
\S\ref{subsection: mean dimension and 
conditional metric mean dimension}, 
while the new waist inequality will be stated in 
\S\ref{subsection: waist inequality}.

A pair $(X,T)$ is called a \textbf{dynamical system} 
if $X$ is a compact metrizable space and 
$T\colon X\to X$ is a homeomorphism. 
We denote its topological entropy and mean dimension 
by
\[
h_{\mathrm{top}}(X,T) \quad\text{and}\quad \mdim(X,T),
\]
respectively. 
Topological entropy was introduced by 
Adler, Konheim, and McAndrew \cite{AKM65}, and 
mean dimension was introduced by Gromov \cite{Gromov}.
Roughly speaking, the topological entropy
$h_{\mathrm{top}}(X,T)$ measures the number of nats 
per unit time needed to describe the orbits of the system $(X,T)$, 
whereas the mean dimension $\mdim(X,T)$ measures 
the number of parameters per unit time needed to 
describe them.
We give precise definitions of these quantities in
\S\ref{subsection: mean dimension and conditional metric mean dimension}.

A fundamental result of Lindenstrauss and Weiss
\cite[Theorem 4.2]{Lindenstrauss--Weiss} states that
finite topological entropy implies zero mean dimension.

\begin{theorem}[Lindenstrauss--Weiss, 2000]
  \label{theorem: Lindenstrauss--Weiss}
 If a dynamical system $(X,T)$ has finite topological entropy
 \[ h_{\mathrm{top}}(X,T) < \infty, \]
 then it has zero mean dimension $\mdim(X, T) = 0$.
\end{theorem}  

This may be viewed as a dynamical analogue of the elementary fact that
a finite set has zero topological dimension.
A further study of the relationship between topological entropy 
and mean dimension was carried out by 
Lindenstrauss \cite{Lindenstrauss}.

One of the initial motivations for this paper was a question posed by
the third and fourth authors in \cite[Problem 1.7]{Shi--Tsukamoto}
concerning a \textit{relative version} 
of Theorem \ref{theorem: Lindenstrauss--Weiss}.
Let $(X,T)$ and $(Y,S)$ be dynamical systems, and let 
$\pi\colon X\to Y$ be an equivariant continuous map.
Here \lq\lq{}equivariant\rq\rq{} means that $\pi\circ T = S\circ \pi$.
We often write $\pi\colon (X,T)\to (Y,S)$ to emphasize the underlying 
dynamics $T$ and $S$.
We denote the topological conditional entropy of $\pi$ by 
(see \S \ref{subsection: mean dimension and conditional 
metric mean dimension} 
for its precise definition)
\[ h_{\mathrm{top}}(\pi, T) := 
   \sup_{y\in Y} h_{\mathrm{top}}(\pi^{-1}(y), T). \]

The following question was raised in \cite[Problem 1.7]{Shi--Tsukamoto}.

\begin{problem}[Shi--Tsukamoto, 2023]
  \label{problem: Shi--Tsukamoto}
  Suppose $\mdim(X, T) > \mdim(Y, S)$.
  Does every equivariant continuous map $\pi\colon (X,T)\to (Y,S)$ have
  infinite topological conditional entropy?
\end{problem}

The paper \cite[Theorem 1.1]{Shi--Tsukamoto} gave a positive answer when 
\[ (X, T) = \left(\left([0,1]^a\right)^{\mathbb{Z}}, \mathrm{shift}\right), 
\quad 
   (Y, S) = \left(\left([0,1]^b\right)^{\mathbb{Z}}, \mathrm{shift} \right), 
   \quad 
    (a > b). \]
In this case, mean dimension is given by $\mdim(X, T) = a$ 
and $\mdim(Y, S) = b$.
Another paper \cite[Theorem 1.2]{Tsukamoto_waist2} 
further gave a positive answer when
$(X, T) = \left(\left([0,1]^a\right)^{\mathbb{Z}}, \mathrm{shift}\right)$ 
and $(Y, S)$ is an arbitrary dynamical system with $\mdim(Y, S) < a$.

Our first main result
affirmatively solves Problem \ref{problem: Shi--Tsukamoto}
in full generality:

\begin{theorem}
  \label{theorem: mean dimension and conditional entropy}
Let $(X, T)$ and $(Y, S)$ be dynamical systems.
If there exists an equivariant continuous map $\pi\colon (X,T)\to (Y,S)$
with finite topological conditional entropy
\[ h_{\mathrm{top}}(\pi, T) < \infty, \]
then we have
\[ \mdim(X, T) \leq \mdim(Y, S). \]
\end{theorem}

To the best of our knowledge, 
this is the first general result 
relating topological conditional entropy to mean dimension.
When $(Y, S)$ is a single point, the theorem reduces to
Theorem \ref{theorem: Lindenstrauss--Weiss} 
of Lindenstrauss and Weiss.
Theorem 
\ref{theorem: mean dimension and conditional entropy} 
may also be viewed as a dynamical analogue of
the fact that
a continuous map between compact metrizable spaces 
with uniformly bounded fiber cardinalities 
does not decrease topological dimension.

\subsection{Mean dimension and conditional metric mean dimension}
\label{subsection: mean dimension and conditional metric mean dimension}

Let $f\colon X\to Y$ be a continuous map 
between compact metrizable spaces $X$ and $Y$.
Then the Hurewicz theorem 
\cite[p.91, Theorem VI. 7]{Hurewicz--Wallman} 
in classical dimension theory
states that
\begin{equation} 
  \label{eq: classical Hurewicz}
  \dim X \leq \dim Y + \sup_{y\in Y} \dim f^{-1}(y). 
\end{equation}
Here $\dim(\cdot)$ denotes the topological dimension.

It is natural to expect an analogous inequality 
to hold for mean dimension.
However, it turns out that a naive analogue fails; for 
an equivariant continuous map $\pi\colon (X,T)\to (Y,S)$ 
between dynamical systems, the inequality
\begin{equation}
  \label{equation: naive Hurewicz inequality}
 \mdim(X, T) \leq \mdim(Y, S) 
 + \sup_{y\in Y} \mdim(\pi^{-1}(y), T) 
\end{equation} 
does not hold in general 
\cite[Theorem 1.3]{Tsukamoto_Hurewicz}.
Our main result below shows that 
a valid Hurewicz-type inequality
can be obtained by replacing
$\sup_{y\in Y}\mdim(\pi^{-1}(y),T)$
with the \textit{conditional metric mean dimension}.
Before stating our main result, 
we introduce the precise definitions of 
several quantities related to entropy 
and mean dimension.

Let $(X,\mathbf{d})$ be a compact metric space 
and let $\varepsilon>0$.
We define the \textbf{$\varepsilon$-covering number}
$\#(X,\mathbf{d},\varepsilon)$
to be the minimum number of open sets 
of diameter less than $\varepsilon$ needed to cover $X$.
Let $Y$ be a topological space and 
let $f\colon X\to Y$ be a continuous map.
We say that $f$ is an \textbf{$\varepsilon$-embedding}
if every fiber $f^{-1}(y)$ 
has diameter smaller than $\varepsilon$.
We define the \textbf{$\varepsilon$-width dimension}
$\widim_\varepsilon(X,\mathbf{d})$ 
to be the minimum integer $n$ such that 
there exist an $n$-dimensional finite simplicial complex 
$P$ and an $\varepsilon$-embedding $f\colon X\to P$.
(A simplicial complex $P$ is said to be 
\lq\lq{}finite\rq\rq{} 
if it consists of finitely many simplices.)

Let $(X,T)$ be a dynamical system 
and let $\mathbf{d}$ be a compatible metric on $X$.
For a positive integer $N$, 
we define the Bowen metric $\mathbf{d}_N$ on $X$ by
\[ \mathbf{d}_N(x, y) 
= \max_{0\leq n < N} \mathbf{d}(T^n x, T^n y), \quad 
 (x, y\in X). \]
We define the \textbf{topological entropy} of $(X,T)$ by
\[ h_{\mathrm{top}}(X,T)
= \lim_{\varepsilon\to 0} \left(\lim_{N\to \infty}
\frac{\log \#(X, \mathbf{d}_N, \varepsilon)}{N}\right). \]
We define the \textbf{mean dimension} of $(X,T)$ by 
\[ \mdim(X, T) = 
\lim_{\varepsilon\to 0} \left(\lim_{N\to \infty} 
\frac{\widim_\varepsilon(X, \mathbf{d}_N)}{N}\right). \]
We can also define the topological entropy and mean dimension for 
(not necessarily invariant) closed subsets $A\subset X$ by 
\begin{align*}
 h_{\mathrm{top}}(A,T) &= 
 \lim_{\varepsilon\to 0} \left(\limsup_{N\to \infty}
\frac{\log \#(A, \mathbf{d}_N, \varepsilon)}{N}\right), \\
\mdim(A, T) &=
\lim_{\varepsilon\to 0} \left(\limsup_{N\to \infty}
\frac{\widim_\varepsilon(A, \mathbf{d}_N)}{N}\right). 
\end{align*}
When $A=\emptyset$, 
we define both the quantities to be $-\infty$.
Topological entropy and mean dimension are 
topological invariants of dynamical systems; they are 
independent of the choice of a
compatible metric $\mathbf{d}$ on $X$.

We define the \textbf{lower metric mean dimension} 
of $(X,T)$ by
\[ \lmdimM(X, T, \mathbf{d}) =
\liminf_{\varepsilon\to 0} \left(\lim_{N\to \infty}
\frac{\log \#(X, \mathbf{d}_N, \varepsilon)}
{N\log(1/\varepsilon)}\right). \]
This is a dynamical analogue 
of the lower Minkowski dimension 
in classical dimension theory.
One can also consider the 
\textbf{upper metric mean dimension}, 
but we do not need it here.
If $h_{\mathrm{top}}(X,T) < \infty$, then
$\lmdimM(X, T, \mathbf{d}) = 0$.
Unlike topological entropy and mean dimension,
metric mean dimension generally depends on the choice of
the compatible metric $\mathbf{d}$.

Metric mean dimension was introduced 
by Lindenstrauss and Weiss \cite{Lindenstrauss--Weiss}.
They proved the following fundamental inequality
\cite[Theorem 4.2]{Lindenstrauss--Weiss}.
\begin{theorem}[Lindenstrauss--Weiss, 2000]
  \label{theorem: Lindenstrauss--Weiss inequality}
  Let $(X,T)$ be a dynamical system 
  with a compatible metric $\mathbf{d}$ 
  on $X$. Then we have
  \[ \mdim(X, T) \leq \lmdimM(X, T, \mathbf{d}). \]
\end{theorem}
Theorem \ref{theorem: Lindenstrauss--Weiss} 
in \S \ref{subsection: entropy and mean dimension} 
immediately follows from this result.

We next turn to the relative setting.
Let $(X,T)$ and $(Y,S)$ be dynamical systems,
and let $\pi\colon (X,T)\to (Y,S)$ 
be an equivariant continuous map.
Let $\mathbf{d}$ be a compatible metric on $X$.
We define the \textbf{topological conditional entropy} 
of $\pi$ by
\[ h_{\mathrm{top}}(\pi, T) 
= \lim_{\varepsilon\to 0} \left(\lim_{N\to \infty}
   \frac{\log 
   \left(\sup_{y\in Y}
   \#(\pi^{-1}(y), 
   \mathbf{d}_N, \varepsilon)\right)}{N}\right). \]   
It is known that \cite[Theorem 6.8.3]{Downarowicz}
\[ h_{\mathrm{top}}(\pi, T) 
= \sup_{y\in Y} h_{\mathrm{top}}(\pi^{-1}(y), T). \]
We define the \textbf{lower conditional metric mean dimension} of $\pi$ by
\[ \lmdimM(\pi, T, \mathbf{d}) 
= \liminf_{\varepsilon\to 0} \left(\lim_{N\to \infty}
   \frac{\log \left(\sup_{y\in Y}
   \#(\pi^{-1}(y), \mathbf{d}_N, \varepsilon)\right)}
   {N\log(1/\varepsilon)}\right). \]
This was first introduced by 
Liang \cite{Liang}.
If $h_{\mathrm{top}}(\pi, T) < \infty$, then
$\lmdimM(\pi, T, \mathbf{d}) = 0$.

The following is the main theorem of this paper.

\begin{theorem}[Main theorem]
  \label{theorem: main theorem}
 Let $\pi\colon (X, T)\to (Y,S)$ be an equivariant continuous map
 between dynamical systems, 
 and let $\mathbf{d}$ be a compatible metric on $X$.
 Then we have 
 \begin{equation} \label{eq: main theorem}
  \mdim(X, T) \leq \mdim(Y, S) + \lmdimM(\pi, T, \mathbf{d}). 
 \end{equation} 
\end{theorem}

Theorem \ref{theorem: mean dimension and conditional entropy} in 
\S \ref{subsection: entropy and mean dimension}
is an immediate corollary of this theorem
since $\lmdimM(\pi, T, \mathbf{d}) = 0$ when
$h_{\mathrm{top}}(\pi, T) < \infty$.

The third and fourth authors \cite[Theorem 1.3]{Shi--Tsukamoto} proved 
\eqref{eq: main theorem} in the case that 
\begin{equation}
\label{eq: full shift case}
(X, T) = 
\left(\left([0,1]^a\right)^{\mathbb{Z}}, \mathrm{shift}\right),  
\quad 
(Y, S) = 
\left(\left([0,1]^b\right)^{\mathbb{Z}}, 
\mathrm{shift} \right)
\end{equation}
with the metric $\mathbf{d}$ on $X$ defined by 
\begin{equation}
  \label{eq: metric on full shift in introduction}
  \mathbf{d}(x, x') 
  = \sum_{n\in \mathbb{Z}} 2^{-|n|} \norm{x_n - x'_n}.
\end{equation}
Here $\norm{\cdot}$ 
is the standard Euclidean norm on $\mathbb{R}^a$.
Another paper \cite{Tsukamoto_waist2} proved it
when $(X, T) = 
\left(\left([0,1]^a\right)^{\mathbb{Z}}, 
\mathrm{shift}\right)$
with the above metric $\mathbf{d}$ 
and $(Y, S)$ is an arbitrary dynamical system.
Now Theorem \ref{theorem: main theorem} 
establishes it in full generality.

It is easy to see that the inequality 
\eqref{eq: main theorem} 
is sharp:
Let $(X, T)$ and $(Y, S)$ be given by 
\eqref{eq: full shift case} with $a\geq b$.
We equip $X$ with the metric $\mathbf{d}$ defined by  
\eqref{eq: metric on full shift in introduction}.
Let $\varphi\colon [0,1]^a\to [0,1]^b$ be the projection 
onto the first $b$ coordinates, and define 
$\pi\colon (X, T)\to (Y,S)$ by
\[ \pi\left((x_n)_{n\in \mathbb{Z}}\right)
= \left(\varphi(x_n)\right)_{n\in \mathbb{Z}}. \]
Then we have
\[ \mdim(X, T) = a, \quad 
   \mdim(Y, S) = b, \quad 
   \lmdimM(\pi, T, \mathbf{d}) = a-b. \]

\subsection{Waist inequality}
\label{subsection: waist inequality}

The earlier works \cite{Shi--Tsukamoto,Tsukamoto_waist2}
established special cases of Theorem \ref{theorem: main theorem}
using waist inequalities.
The original waist inequality was discovered by 
Gromov \cite{Gromov_waist}.
Let $S^n$ be the standard $n$-dimensional sphere.
Gromov proved that, for every continuous map 
$f\colon S^n  \to \mathbb{R}^m$ 
$(m\leq n)$, there exists a point 
$y\in \mathbb{R}^m$ such that 
\begin{equation}
 \label{eq: Gromov waist inequality} 
 \mathrm{Vol}
 \left(\text{the $r$-neighborhood of }f^{-1}(y)\right) 
 \geq \mathrm{Vol}
 \left(\text{the $r$-neighborhood of }S^{n-m}\right)
\end{equation}
for all $r>0$. 
Here $\mathrm{Vol}(\cdot)$ 
denotes the $n$-dimensional volume,
and $S^{n-m} \subset S^n$ is 
an equatorial $(n-m)$-dimensional sphere.
Roughly speaking, the waist inequality guarantees that
every continuous map $f\colon S^n\to\mathbb{R}^m$
has a \lq\lq{}large\rq\rq{} fiber $f^{-1}(y)$.

Klartag \cite{Klartag} developed
waist inequalities for convex bodies.
The proofs in \cite{Shi--Tsukamoto,Tsukamoto_waist2}
are based on one of Klartag\rq{}s inequalities.

The proof of Theorem \ref{theorem: main theorem} 
also relies on a waist inequality. 
While existing waist inequalities were sufficient 
for the earlier works \cite{Shi--Tsukamoto,Tsukamoto_waist2}, 
they are not sufficient for our present purpose. 
We therefore establish a new waist inequality.

Let $P$ and $Q$ be finite simplicial complexes, 
and let $f\colon P\to Q$ be a continuous map. 
We would like to show that some fiber 
$f^{-1}(y)$ $(y\in Q)$ is 
\lq\lq{}large\rq\rq{} 
in an appropriate sense. 
The first question is how to make 
this notion of largeness precise. 
Perhaps the simplest approach is 
to equip $P$ with a metric and 
measure the size of the fibers with respect to this metric.
However, 
it turns out that a different approach is more fruitful 
for our purposes. 
Let $\Phi\colon P\to\mathbb{R}^n$ 
be an auxiliary continuous map. 
We regard $\Phi$ as an \lq\lq{}observable\rq\rq{} 
and measure the size of each fiber $f^{-1}(y)$ 
through the size of 
its image $\Phi(f^{-1}(y))$ in $\mathbb{R}^n$.\footnote{
This viewpoint is closely related to equipping
$P$ with the pseudometric $\mathbf{d}$ defined by
$\mathbf{d}(x,x'):=\norm{\Phi(x)-\Phi(x')}$.
}
We therefore seek to formulate a theorem asserting that, 
for every continuous map $f\colon P\to Q$, 
there exists a point $y\in Q$ 
such that $\Phi(f^{-1}(y))$ is large.

Let $\norm{\cdot}$ be the standard 
Euclidean norm on $\mathbb{R}^n$.
It is convenient to use 
\[ \rho_n(x) := \frac{\norm{x}}{\sqrt{n}}\quad 
    (x\in \mathbb{R}^n) \]
as a norm on $\mathbb{R}^n$.
For $r>0$ and a subset $A\subset \mathbb{R}^n$, let
$\#\left(A, \rho_n, r\right)$ denote the minimum number of 
open sets of $\rho_n$-diameter 
less than $r$ needed to cover 
$A$.

Let $P$ be a finite simplicial complex and let 
$\Phi\colon P\to \mathbb{R}^n$ be a continuous map.
For an integer $d\geq 0$, 
we define the \textbf{$d$-radius} of $\Phi$,
following Gromov \cite[p.~127]{Gromov_filling}, by
\begin{equation}
  \label{eq: d-radius}
\begin{split}
&\mathcal{R}_d(\Phi) \\
&:=
\inf\left\{
\sup_{x\in P}
\rho_n\left(\Phi(x)-(\psi\circ\varphi)(x)\right)
\;\middle|\;
\begin{gathered}
K \text{ is a finite simplicial complex},
\, \dim K \leq d\\
\varphi\colon P\to K,\quad
\psi\colon K\to\mathbb{R}^n
\text{ are continuous}
\end{gathered}
\right\}.
\end{split}
\end{equation}
The quantity $\mathcal{R}_d(\Phi)$ measures 
how well $\Phi$ can be uniformly approximated 
by maps that factor 
through finite simplicial complexes 
of dimension at most $d$.

We have $\mathcal{R}_d(\Phi) = 0$
for $d\geq \dim P$
since we can take $K = P$, $\varphi = \mathrm{id}$ and 
$\psi = \Phi$ in the above definition.
Similarly we also have $\mathcal{R}_d(\Phi)=0$ for $d\geq n$.

The following is our new waist inequality.

\begin{theorem}[New waist inequality]
  \label{theorem: new waist inequality}
 Let $3\leq d < n$ be integers.
 Let $P$ and $Q$ be finite simplicial complexes.
 Let $f\colon P\to Q$ and $\Phi\colon P\to \mathbb{R}^n$ 
 be continuous maps. Suppose 
 \begin{equation} \label{eq: dimensional assumption}
   \dim P \leq 2d-1, \quad \dim Q \leq d-3. 
 \end{equation}  
 Let $s$ be a real number with $0< s < d-\dim Q -2$ and set 
 \[ \kappa := 1 - \frac{s+2}{d-\dim Q} >0. \]
 Then for every $r>0$ we have
 \begin{equation}
   \label{eq: new waist inequality}
   \sup_{y\in Q} \#\left(\Phi(f^{-1}(y)), \rho_n, r\right)
  \geq \frac{1}{\dim Q + 1} 
  \left(\frac{\kappa \mathcal{R}_d(\Phi)}{r}\right)^s.
 \end{equation}
\end{theorem}

The statement of this theorem may look 
somewhat complicated at first sight. 
Let us explain its meaning 
and the roles of the various parameters.

In the application to 
Theorem \ref{theorem: main theorem}, 
the ambient dimension $n$ is larger than $\dim P$, and
$\dim P$ is much larger than $\dim Q$.
We choose $d$ close to $\dim P$. 
Thus the dimensional assumptions 
\eqref{eq: dimensional assumption} are easily satisfied.

For example, suppose that
\[
d=\dim P-1,\qquad
s=d-\dim Q-3=\dim P-\dim Q-4>0.
\]
Then, 
\[ \kappa = \frac{1}{\dim P - \dim Q -1}. \]
Provided that $\mathcal{R}_d(\Phi)>0$, 
the inequality \eqref{eq: new waist inequality} gives
\[
\sup_{y\in Q}
\#\left(\Phi(f^{-1}(y)),\rho_n,r\right)
\geq
C\left(\frac{1}{r}\right)^{\dim P-\dim Q-4},
\]
where $C>0$ is independent of $r$.
Heuristically,
one might expect the exponent on the right-hand side 
to be $\dim P-\dim Q$, 
the expected dimension of a fiber of $f$. 
Our theorem gives an exponent close to this heuristic value.
The key point is that the constant $C$ depends 
only on $\dim P$, $\dim Q$, and $\mathcal{R}_d(\Phi)$. 
Thus, apart from its dimension, 
the potentially complicated topology of $P$ 
enters the estimate only through the single quantity 
$\mathcal{R}_d(\Phi)$.

Compared with Gromov's original waist inequality
\eqref{eq: Gromov waist inequality},
our inequality is weaker in two respects.
First, neither the exponent nor the constant in
\eqref{eq: new waist inequality} is optimal,
at least in the special case where $P=S^{n-1}$,
$Q=[0,1]^m$, and
$\Phi\colon S^{n-1}\hookrightarrow\mathbb{R}^n$
is the standard inclusion.
Second, the point $y\in Q$ giving a large fiber
may depend on $r>0$,
whereas Gromov's inequality guarantees a single point $y$
that works for all $r>0$.

On the other hand,
the main strength of our inequality is
its topological generality:
$P$ and $Q$ can be arbitrary finite simplicial complexes
satisfying the dimensional assumptions
\eqref{eq: dimensional assumption}.
Although the exponent and the constant 
appearing in our estimate are unlikely to be optimal, 
the estimate is strong enough 
for our application to Theorem \ref{theorem: main theorem}, 
for which this topological generality is essential.

\begin{remark}
  \begin{enumerate}
  \item We do not know to what extent 
  the dimensional assumptions
\eqref{eq: dimensional assumption} are necessary.
However, it is easy to see that the condition
$\dim Q\leq d-3$ cannot be relaxed beyond $\dim Q\leq d$.
Indeed, let 
$P=Q=[0,1]^{d+1}$.
Let $f\colon P\to Q$ be the identity map, and let
$\Phi\colon P\hookrightarrow\mathbb{R}^{d+1}$
be the standard inclusion.
Then, for every $r>0$,
\[
\sup_{y\in Q}
\#\left(\Phi(f^{-1}(y)),\rho_{d+1},r\right)=1,
\]
whereas $\mathcal{R}_d(\Phi)>0$.
Hence no inequality of the form
\eqref{eq: new waist inequality}
with a positive exponent $s$ and a positive constant
$\kappa$
can hold in general when $\dim Q=d+1$.
Notice that $\dim P=d+1\leq 2d-1$
for $d\geq 2$ in this example.

On the other hand, we do not know whether the condition
$\dim P\leq 2d-1$ is necessary for
\eqref{eq: new waist inequality}.
Our proof uses this assumption
through the application of
the Freudenthal suspension theorem
in \S
\ref{section: proof of proposition low dimensional approximation}.
In \S \ref{subsection: weaker waist inequality},
we prove a weaker waist inequality
(Theorem \ref{theorem: weaker waist inequality})
without any restriction on $\dim P$.
Although its estimate has less satisfactory constants,
it is still sufficient to prove
Theorem \ref{theorem: main theorem}.
\item 
Even when $Q$ is a single point,
Theorem~\ref{theorem: new waist inequality}
gives a nontrivial covering estimate.
Indeed, in this case the inequality
\eqref{eq: new waist inequality}
becomes
\[
\#(\Phi(P),\rho_n,r)
\ge
\left(
\frac{\kappa\mathcal R_d(\Phi)}{r}
\right)^s,
\qquad
\kappa=1-\frac{s+2}{d},
\quad 0<s<d-2.
\]
For comparison, Nabutovsky's
inequality on the Alexandrov width
\cite[Theorem~1.4(2)]{Nabutovsky}\footnote{
This theorem is a more general result 
and originally formulated in terms of
Alexandrov width and Hausdorff content;
\eqref{eq: Nabutovsky inequality}
is a consequence of that theorem.
}
implies
\begin{equation}
  \label{eq: Nabutovsky inequality}
\#(\Phi(P),\rho_n,r)
\ge
\left(
\frac{\mathcal R_d(\Phi)}
{4(d+1)r}
\right)^{d+1}.
\end{equation}
Although the latter has the larger exponent $d+1$,
its scale factor $1/[4(d+1)]$ tends to zero
as $d\to\infty$.
In contrast, our scale factor $\kappa$
satisfies $\kappa\ge\delta-2/d$
whenever $s\le(1-\delta)d$
for a fixed $\delta>0$.
This feature makes our inequality
particularly suitable for applications
to dynamical systems,
where one naturally fixes the scale $r$
and lets $d\to\infty$.
  \end{enumerate}
\end{remark}

For the proof of Theorem \ref{theorem: main theorem}, 
it is useful to formulate an $\ell^\infty$-version of 
Theorem \ref{theorem: new waist inequality}.
Let $\norm{x}_\infty$ denote the max-norm:
\[ \norm{x}_\infty
  := \max_{1\leq i \leq n} |x_i| \quad 
  \text{for } x= (x_1, \dots, x_n) \in \mathbb{R}^n. \]
Let $P$ be a finite simplicial complex and 
$\Phi\colon P\to \mathbb{R}^n$ be a continuous map.
Let $\mathcal{R}^\infty_d(\Phi)$
be defined as in \eqref{eq: d-radius},
with $\rho_n$ replaced by $\norm{\cdot}_\infty$.
Thus it is the infimum of the uniform max-norm errors
in approximating $\Phi$ by maps that factor through
finite simplicial complexes of dimension at most $d$.

The following lemma provides a bridge 
between the $\rho_n$-version 
and the $\ell^\infty$-version of the waist inequality.

\begin{lemma}
 \label{lemma: rho_n and ell^infty}
 For all integers $d, k\geq 0$, 
 \[ \mathcal{R}^\infty_{d+k}(\Phi) \leq
     \sqrt{\frac{n}{k+1}} \, \mathcal{R}_d(\Phi). \]
\end{lemma}

\begin{proof}
  For a real number $t$, we define 
  \[ t_+ := \max(t, 0), \quad 
    \mathrm{sign}(t) = \begin{cases}
    1 & (t>0) \\ 0 & (t=0) \\ -1 & (t<0)
    \end{cases}. \]
 Let $r>\mathcal{R}_d(\Phi)$ be arbitrary.
There are a finite simplicial 
complex $K$ of $\dim K \leq d$ and 
continuous maps $\varphi\colon P\to K$ and 
$\psi\colon K\to \mathbb{R}^n$ such that 
\[ \rho_n\left(\Phi(x) - 
\psi\left(\varphi(x)\right)\right) < r 
  \qquad 
   (\forall x\in P). \]
Let $x\in P$ and set 
\[ u = (u_1, \dots, u_n) 
:= \Phi(x) - \psi\left(\varphi(x)\right). \]
We have $\rho_n(u) < r$. Put 
\[ R = r \sqrt{\frac{n}{k+1}}, \qquad 
   v_i = \mathrm{sign}(u_i) \left(|u_i|-R\right)_+ \quad 
  (1\leq i \leq n). \]
Set $v = (v_1, \dots, v_n)$.
We have $\norm{u-v}_\infty \leq R$.
If $k+1$ coordinates of $v$ were nonzero, then 
$\norm{u}^2 \geq (k+1) R^2  = n r^2$; a contradiction.
Hence $v$ has at most $k$ nonzero coordinates and belongs to 
a finite union $\Sigma_k$ of coordinate cubes of dimension at 
most $k$.\footnote{More rigorously, 
$\Sigma_k$ is defined as follows:
We take $M>0$ so that $v(P) \subset [-M, M]^n$ 
and define 
\[ \Sigma_k 
= \{w\in [-M, M]^n\mid |\{i\mid w_i\neq 0\}|\leq k\}.\]
This becomes a finite simplicial complex of dimension 
$\min(k,n)$.}

Define $\varphi^\prime \colon P\to K\times \Sigma_k$ and 
$\psi^\prime \colon K\times \Sigma_k \to \mathbb{R}^n$ by 
\[ \varphi^\prime(x) = (\varphi(x), v), \qquad 
    \psi^\prime(z, w) = \psi(z) + w. \]
Then 
\[ \Phi(x)-\psi^\prime\left(\varphi^\prime(x)\right)
   = \Phi(x) - \psi(\varphi(x)) -v = u-v. \]   
Since $\norm{u-v}_\infty \leq R$,
\[ \sup_{x\in P} 
\norm{\Phi(x)-\psi^\prime
\left(\varphi^\prime(x)\right)}_\infty
\leq R =  r \sqrt{\frac{n}{k+1}}. \]
Letting $r\to \mathcal{R}_d(\Phi)$, 
this proves the claim.
\end{proof}

The following is the $\ell^\infty$-version 
of our waist inequality.
This will be used 
in the proof of Theorem \ref{theorem: main theorem}.

\begin{corollary}[$\ell^\infty$-waist inequality]
  \label{corollary: ell^infty waist inequality}
   Let $3\leq d < n$ be integers.
 Let $P$ and $Q$ be finite simplicial complexes.
 Let $f\colon P\to Q$ and $\Phi\colon P\to \mathbb{R}^n$ 
 be continuous maps. Suppose 
 \begin{equation*} 
   \dim P \leq 2d-1, \quad \dim Q \leq d-3. 
 \end{equation*}  
 Let $s$ be a real number with $0< s < d-\dim Q -2$ and set 
 \[ \kappa = 1 - \frac{s+2}{d-\dim Q} >0. \]
 Then for every $r>0$ and every integer $k\geq 0$ we have
 \begin{equation}
   \label{eq: ell^infty waist inequality}
   \sup_{y\in Q} \#\left(\Phi(f^{-1}(y)), 
   \norm{\cdot}_\infty, r\right)
  \geq \frac{1}{\dim Q + 1} 
  \left(\kappa \sqrt{\frac{k+1}{n}} \, 
  \frac{\mathcal{R}^\infty_{d+k}(\Phi)}{r}\right)^s.
 \end{equation}
 Here $\#\left(\Phi(f^{-1}(y)), 
   \norm{\cdot}_\infty, r\right)$ is the minimum number of 
   open sets of $\norm{\cdot}_\infty$-diameter less than 
   $r$ needed to cover $\Phi(f^{-1}(y))$.
\end{corollary}

\begin{proof}
Since $\rho_n \leq \norm{\cdot}_\infty$, we have 
\[ \#\left(\Phi(f^{-1}(y)), 
   \norm{\cdot}_\infty, r\right) \geq
   \#\left(\Phi(f^{-1}(y)), 
   \rho_n, r\right). \]
From Lemma \ref{lemma: rho_n and ell^infty}
\[ \mathcal{R}_{d}(\Phi) \geq 
    \sqrt{\frac{k+1}{n}} \, 
    \mathcal{R}^\infty_{d+k}(\Phi). \]
Then \eqref{eq: ell^infty waist inequality}
follows from Theorem 
\ref{theorem: new waist inequality}.
\end{proof}

\subsection{Organization of the paper}
\label{subsection: organization of the paper}

In \S\ref{section: proof of main theorem},
we prove the main theorem assuming our new waist inequality.
In \S\ref{section: two main proposition and the proof of waist inequality},
we derive the waist inequality from two propositions:
quantitative avoidance and low-dimensional approximation.
The quantitative avoidance proposition is proved in
\S\ref{section: proof of proposition avoidance}
using Morse theory on spheres.
The low-dimensional approximation proposition is proved in
\S\ref{section: proof of proposition low dimensional approximation},
where we first explain the motivation and review
the necessary background from algebraic topology.
\S\ref{section: proof of proposition low dimensional approximation}
is probably the most technical part of the paper.

For readers primarily interested in mean dimension theory,
\S\ref{subsection: weaker waist inequality}
provides an alternative proof of the main theorem
based on a weaker waist inequality.
This approach avoids the low-dimensional approximation
proposition and its algebraic-topological proof.
\S\ref{subsection: weaker waist inequality} 
is not needed for the rest of the paper
and may be skipped.

\subsection{Use of generative AI}
  \label{subsection: use of generative}

We briefly explain the role of generative AI 
in the development
of this research and the preparation of the manuscript.  

The authors identified the research problem 
and recognized the need for a new waist inequality. 
In exploring this problem, 
the authors initially consulted ChatGPT (GPT-5.4) 
and subsequently used ChatGPT (GPT-5.6 Sol). 
Through interactions with the authors, 
GPT-5.6 Sol generated key ideas for the 
proof of the new waist inequality, 
particularly those concerning 
quantitative avoidance and low-dimensional approximation. 
These ideas were further developed 
and refined through subsequent discussions. 

ChatGPT (GPT-6 Astra) was also extensively 
used in preparing the manuscript, including the figures. 
The authors take full responsibility 
for the content of the paper.

 \section{Proof of Theorem \ref{theorem: main theorem}}
\label{section: proof of main theorem}

In this section, we prove Theorem \ref{theorem: main theorem}
using the new waist inequality
(Theorem \ref{theorem: new waist inequality}),
which will be established in later sections.
Let $\pi\colon (X, T)\to (Y,S)$ be an equivariant 
continuous map between dynamical systems, and let 
$\mathbf{d}$ be a compatible metric on $X$.
We would like to prove 
\begin{equation}
  \label{eq: Hurewicz type inequality for mean dimension, revisited}
 \mdim(X, T) \leq \mdim(Y, S) + 
   \lmdimM(\pi, T, \mathbf{d}). 
\end{equation}   
If $\mdim(X, T)=0$ or $\mdim(Y, S)=\infty$ 
then 
\eqref{eq: Hurewicz type inequality for mean dimension, revisited}
is trivial.
So we assume $\mdim(X,T)>0$ and $\mdim(Y,S)<\infty$ 
throughout this section.

We first treat the case $\mdim(X,T)<\infty$
in \S\ref{subsection: the case of finite mean dimension}.
The case $\mdim(X,T)=\infty$ requires an additional argument
and will be treated in
\S\ref{subsection: the case of infinite mean dimension}.

\subsection{The case of $\mdim(X,T)<\infty$}
\label{subsection: the case of finite mean dimension}
In this subsection, we assume that 
$\pi\colon (X, T)\to (Y,S)$ is an equivariant continuous 
map between dynamical systems such that 
\[ 0< \mdim(X, T) < \infty, \quad \mdim(Y, S) < \infty. \]
We take a compatible metric $\mathbf{d}$ on $X$.

We take arbitrary positive numbers $a$ and $b$ with 
\[ a < \mdim(X,T), \quad \mdim(Y, S) < b. \]
We would like to prove 
\[ a \leq b + \lmdimM(\pi, T, \mathbf{d}). \]

For $r>0$ we set 
\[ S_r(\pi) := \lim_{N\to \infty} 
   \frac{\log \left(\sup_{y\in Y}
   \#\left(\pi^{-1}(y), \mathbf{d}_N, r\right)\right)}{N}. 
   \]
Recall that 
\[ \lmdimM(\pi, T, \mathbf{d}) 
    = \liminf_{r\to 0} \frac{S_r(\pi)}{\log (1/r)}.\]

    For a map $F\colon X\to [0,1]^{\ell}$ and $N>0$, we set 
\[ F_N(x) = \left(F(x), F(T x), \dots, F(T^{N-1}x)\right),
    \qquad (x\in X). \]    

The following proposition is the key step 
in reducing our dynamical problem 
to a finite-dimensional geometric problem, 
to which the waist inequality can be applied.
The basic idea is to model the map $\pi\colon X\to Y$, 
viewed over long orbit segments, 
by a map between finite-dimensional simplicial complexes. 
The point is to construct such a model 
without losing the quantitative information 
about the spaces and their fibers.

\begin{proposition}[Polyhedral modeling]
  \label{proposition: polyhedral modeling}
  Let $F\colon (X, \mathbf{d})\to 
  ([0,1]^{\ell}, \norm{\cdot}_\infty)$ 
  be a $1$-Lipschitz map.
  For any positive numbers $\delta$ and $\eta$, there are 
  $\varepsilon>0$ and $L\in \mathbb{N}$ such that 
  for every $N\in L\mathbb{N}$ there are finite simplicial 
  complexes $P_N, Q_N$ and continuous maps 
  \[ f_N\colon X\to P_N, \quad 
    g_N\colon P_N \to Q_N, \quad 
    \Phi_N \colon P_N\to [0,1]^{\ell N}, \]
  that satisfy the following conditions. 
   \begin{itemize}
     \item $f_N$ is an $\varepsilon$-embedding 
           with respect to $\mathbf{d}_N$, and 
           $\dim P_N 
           = \widim_\varepsilon(X, \mathbf{d}_N)$.
     \item $\dim Q_N < b N$.
     \item $\norm{\Phi_N\circ f_N(x)-F_N(x)}_\infty < \delta$
            for all $x\in X$.
     \item $\sup_{y\in Q_N} 
     \#\left(\Phi_N\left(g_N^{-1}(y)\right), 
     \norm{\cdot}_\infty,
     \delta \right) \leq 
     \exp\left(\left(S_\delta(\pi)+ \eta\right)N\right)$.
   \end{itemize} 
\end{proposition}  

The above map $f_N\colon X\to P_N$ represents 
orbit segments of length $N$ in $X$ 
by points of a finite-dimensional simplicial complex $P_N$.
The map $g_N\colon P_N\to Q_N$ 
serves as a finite-dimensional polyhedral model of 
the equivariant map $\pi\colon (X,T)\to(Y,S)$ 
over orbit segments of length $N$. 
The observable $\Phi_N$ is constructed so that 
$\Phi_N\circ f_N$ approximates the orbit-observation map
$F_N(x)=(F(x),F(Tx),\ldots,F(T^{N-1}x))$.

\begin{proof}
  We first notice that the following proof 
  of this proposition 
  does not require the assumption 
  $\mdim(X,T) < \infty$.

For a finite simplicial complex $P$, 
we denote the set of 
vertices of $P$ by $\mathrm{Ver}(P)$.
For a simplex $\sigma$ of $P$, 
its \textbf{interior}
means its relative interior in its affine hull,
and its \textbf{boundary} $\partial\sigma$
is the union of all its proper faces.
Thus $\sigma$ is the disjoint union of its
interior and boundary.
In particular, the interior of a vertex is
the vertex itself, and its boundary is empty.
For each $v\in \mathrm{Ver}(P)$, we define 
the \textbf{open star} $O_P(v)$ to be the union of 
interiors of all simplices of $P$ containing $v$.
This is an open neighborhood of $v$, and 
$\{O_P(v) \mid v\in \mathrm{Ver}(P)\}$
forms an open cover of $P$.

\smallskip
\noindent\emph{Step 1: Covers of neighborhoods of fibers.}
Choose $m\geq 1$ such that
\begin{equation}
\label{eq: modeling block entropy}
M:=\sup_{y\in Y}
\#\bigl(\pi^{-1}(y),\mathbf{d}_m,\delta\bigr)
\leq \exp\bigl((S_\delta(\pi)+\eta)m\bigr).
\end{equation}
For each $y\in Y$, choose at most $M$ open subsets of $X$
of $\mathbf{d}_m$-diameter less than $\delta$ covering
$\pi^{-1}(y)$.
(For an empty fiber, the empty family suffices.)
Since $X$ is compact, the same sets cover $\pi^{-1}(W_y)$
for some open neighborhood $W_y$ of $y$.

By compactness, finitely many of these $W_y$ cover $Y$.
Denote the resulting cover by $\mathcal{W}$.
Since only finitely many covering sets have been chosen,
there exists $0<\delta_0<\delta$ such that, for every
$W\in\mathcal{W}$, the set $\pi^{-1}(W)$ is covered by
at most $M$ open subsets of $X$, each of
$\mathbf{d}_m$-diameter less than $\delta_0$.

\smallskip
\noindent
\emph{Step 2: A polyhedral model of the base.}
Fix a compatible metric $\mathbf{d}^\prime$ on $Y$, 
and choose $\tau>0$
such that every subset of $Y$ 
of $\mathbf{d}^\prime$-diameter 
less than $\tau$
is contained in a member of $\mathcal{W}$.

Since $\mdim(Y,S)<b$, we can choose $L\in m\mathbb{N}$
such that
$\widim_\tau(Y,\mathbf{d}^\prime_L)<bL$.
Then we can take a finite polyhedron $Q$ 
of dimension less than $bL$ and 
a $\tau$-embedding 
$h \colon (Y, \mathbf{d}^\prime_L) \to Q$.
By subdividing $Q$ sufficiently fine, 
we can assume that 
$\diam \left(h^{-1}(O_Q(v)),\mathbf{d}^\prime_L\right)
< \tau$
for all vertices $v\in \mathrm{Ver}(Q)$.

Let $v\in \mathrm{Ver}(Q)$.
For each $0\leq i < L/m$, we have
$\diam \left(S^{im}h^{-1}(O_Q(v)),
\mathbf{d}^\prime_m\right)
< \tau$, and hence
we can take 
$W_i \in \mathcal{W}$ with 
\[ S^{i m} h^{-1}(O_Q(v)) \subset W_i. \]
Therefore 
\[ \pi^{-1}\left(h^{-1}(O_Q(v))\right)
   \subset \bigcap_{i=0}^{L/m-1} 
   T^{-im}\pi^{-1}(W_i), \]
and hence   
\begin{equation}
  \label{eq: fiber estimate in base polyhedron}
  \#\left(\pi^{-1}\left(h^{-1}(O_Q(v))\right), 
  \mathbf{d}_L, \delta_0\right) \leq
  \prod_{i=0}^{L/m-1} 
  \#\left(\pi^{-1}(W_i), \mathbf{d}_m, \delta_0\right)
  \leq M^{L/m}
\end{equation}

\smallskip
\noindent
\emph{Step 3: Construction of the polyhedral models.}
Choose $\varepsilon>0$ so small that
\begin{equation}
\label{eq: modeling epsilon}
\delta_0+2\varepsilon<\delta
\end{equation}
and every subset of $X$ of $\mathbf{d}$-diameter less than
$\varepsilon$ is contained in 
$\pi^{-1}\left(h^{-1}(O_Q(v))\right)$ for some
vertex $v$ of $Q$.
This is possible because 
$\{\pi^{-1}\left(h^{-1}(O_Q(v))\right)\mid 
v\in \mathrm{Ver}(Q)\}$ cover $X$.

Let $N=qL$, where $q\geq1$ is an integer.
Choose a finite simplicial complex $P_N$ of 
dimension $\widim_\varepsilon(X, \mathbf{d}_N)$
and an $\varepsilon$-embedding 
$f_N\colon (X, \mathbf{d}_N) \to P_N$. 
By subdividing $P_N$ sufficiently fine, 
we can assume that 
$\diam\left(f_N^{-1}(O_{P_N}(u)), \mathbf{d}_N\right) 
< \varepsilon$ for all vertices 
$u \in \mathrm{Ver}(P_N)$.
Moreover, by replacing $P_N$ with its suitable subcomplex, 
we can also assume that, 
if $u_1, \dots, u_k\in \mathrm{Ver}(P_N)$ 
form a simplex in $P_N$, then 
\[ f_{N}^{-1}
\left(O_{P_N}(u_1)\cap O_{P_N}(u_2)\cap \dots 
\cap O_{P_N}(u_k)\right) \neq \emptyset. \] 
In particular, $f_N^{-1}(O_{P_N}(u))$ 
is non-empty for every 
$u \in \mathrm{Ver}(P_N)$.

We set $Q_N := Q^q$.
This admits a structure of a finite 
simplicial complex with
\[
\dim Q_N=q\dim Q < b q L = bN.
\]

For every $u \in\mathrm{Ver}(P_N)$ and $0\leq j<q$, 
choose a vertex of $Q$ (denoted by 
$g_{N,j}(u)$) such that
\[
T^{jL}\left(f_N^{-1}(O_{P_N}(u))\right)\subset
\pi^{-1}\left(h^{-1}(O_Q(g_{N,j}(u)))\right).
\]
This is possible by the choice of $\varepsilon$.
The assignment 
\[ \mathrm{Ver}(P_N) \ni u \mapsto 
    g_{N,j}(u) \in \mathrm{Ver}(Q) \]
extends to
a simplicial map
$g_{N,j}\colon P_N\to Q$.
Indeed, if $u_1,\ldots,u_k$ form a simplex in $P_N$, 
then there is a point 
$x\in f_{N}^{-1}
\left(O_{P_N}(u_1)\cap \dots 
\cap O_{P_N}(u_k)\right)$.
We have 
\[ \pi(T^{jL}x) \in 
   h^{-1}\left(O_Q(g_{N,j}(u_1)) 
   \cap \dots 
   \cap O_Q(g_{N,j}(u_k))\right). \]
This shows that 
\[ O_Q(g_{N,j}(u_1))\cap \dots 
   \cap O_Q(g_{N,j}(u_k)) \neq \emptyset, \] 
and hence
the vertices 
$g_{N,j}(u_1), \dots, g_{N,j}(u_k)$ 
form a simplex in 
$Q$.    

Set
\[
g_N=(g_{N,0},\ldots,g_{N,q-1})\colon P_N\to Q_N.
\]

Choose $x(u)\in f_N^{-1}(O_{P_N}(u))$ 
for each $u \in \mathrm{Ver}(P_N)$.
Define $\Phi_N\colon P_N\to[0,1]^{\ell N}$ by
$\Phi_N(u)=F_N(x(u))$ on the vertices and extend it
affinely over each simplex.
If $p \in P_N$ lies in a simplex with vertices
$u_1,\ldots,u_k$, then, for every
$x\in f_N^{-1}
\left(O_{P_N}(u_1)\cap\cdots\cap O_{P_N}(u_k)\right)$, 
we have
$\mathbf{d}_N(x(u_i),x)<\varepsilon$ for every $i$
and hence
\begin{equation}
\label{eq: modeling simplex approximation}
\norm{\Phi_N(p)-F_N(x)}_\infty<\varepsilon.
\end{equation}
This follows because $F_N$ is $1$-Lipschitz with respect to
$\mathbf{d}_N$ and the max-norm.
Applying this to the smallest simplex 
containing $f_N(x)$ gives
\[
\norm{\Phi_N(f_N(x))-F_N(x)}_\infty
<\varepsilon<\delta
\qquad(x\in X).
\]

\smallskip
\noindent\emph{Step 4: Covering every model fiber.}
Fix $y=(y_0,\ldots,y_{q-1})\in Q_N$.
For each $j=0,\dots,q-1$, 
choose a vertex $v_j\in \mathrm{Ver}(Q)$ 
with $y_j \in O_Q(v_j)$.
Put
\[
X_y:=\bigcap_{j=0}^{q-1}
T^{-jL}\pi^{-1}\left(h^{-1}(O_Q(v_j))\right)
\subset X.
\]
By \eqref{eq: fiber estimate in base polyhedron}, 
we have (recall $N = qL = q(L/m)m$)
\begin{equation}
\label{eq: modeling orbit cover}
 \#\left(X_y, \mathbf{d}_N, \delta_0\right) 
 \leq M^{(L/m)q} = M^{N/m}.
\end{equation}

Let $p\in g_N^{-1}(y)$, 
and let $\Delta$ be the smallest
simplex containing $p$.
Choose
\[
x\in\bigcap_{u: \text{ a vertex of }\Delta} 
f_N^{-1}(O_{P_N}(u)).
\]
For each $j=0, \dots, q-1$, we have
$g_{N,j}(p)=y_j \in O_Q(v_j)$.
Consequently, some vertex $u$ of $\Delta$ satisfies
$g_{N,j}(u) = v_j$.
It follows that 
$\pi(T^{jL}x)\in h^{-1}\left(O_Q(v_j)\right)$ 
for every $j$,
so $x\in X_y$.
By \eqref{eq: modeling simplex approximation},
\[
\norm{\Phi_N(p)-F_N(x)}_\infty<\varepsilon.
\]
Thus $\Phi_N(g_N^{-1}(y))$ is contained in the open
$\varepsilon$-neighborhood of the set $F_N(X_y)$.

Recall that $\delta_0+2\varepsilon<\delta$ and
$F_N$ is $1$-Lipschitz with respect to 
$\mathbf{d}_N$ on $X$ and the max-norm on 
$[0,1]^{\ell N}$.
It follows from \eqref{eq: modeling orbit cover} and
\eqref{eq: modeling block entropy} that
\[
\#\bigl(\Phi_N(g_N^{-1}(y)),\norm{\cdot}_\infty,\delta\bigr)
\leq M^{N/m}
\leq\exp\bigl((S_\delta(\pi)+\eta)N\bigr).
\]
Taking the supremum over $y\in Q_N$ completes the proof.
\end{proof}

We now complete the proof of
Theorem \ref{theorem: main theorem}
when $\mdim(X,T)<\infty$.
The key is to choose the observable $F$ 
so that a sufficiently
accurate approximation of $F_N$ cannot factor through a
simplicial complex of too small a dimension.
This gives a uniform positive lower bound for the radius
appearing in Corollary 
\ref{corollary: ell^infty waist inequality}.

\begin{proof}[Completion of the 
  proof in the finite mean dimension case]
Put $D=\mdim(X,T)$.
Recall that we took arbitrary real numbers $a, b$ 
with $0<a<D$ and $\mdim(Y,S)<b$.
We would like to prove
\[  a \leq b + \lmdimM(\pi, T, \mathbf{d}). \]
If $a\leq b$, there is nothing to prove, 
so we assume $a>b$.
Choose real numbers $\alpha$ and $\beta$ such that
\begin{equation}
\label{eq: finite case dimension parameters}
\max\{a,D/2\}<\alpha<\beta<D.
\end{equation}
Notice that here we used our assumption $D<\infty$.

\smallskip
\noindent\emph{Step 1: Choosing an observable.}
Since $\beta<D$, there exists $0<\tau<1$ such that
\[
\lim_{N\to\infty}
\frac{\widim_\tau(X,\mathbf{d}_N)}{N}>\beta.
\]
In particular, for all sufficiently large $N$,
\begin{equation}
\label{eq: finite case width lower bound}
\widim_\tau(X,\mathbf{d}_N)>\beta N.
\end{equation}
Choose a $\tau/4$-spanning set\footnote{
For every $x\in X$, some $x_i$ satisfies 
$\mathbf{d}(x, x_i) < \tau/4$.}
 $\{x_1,\ldots,x_\ell\}$ of $X$
and set
\[
F(x)=\bigl(\min\{\mathbf{d}(x,x_i),1\}\bigr)_{i=1}^\ell.
\]
This is a $1$-Lipschitz map from $(X,\mathbf{d})$ to
$([0,1]^\ell,\norm{\cdot}_\infty)$.
By adding zero coordinates to $F$ if necessary, 
we may assume
that $\ell>D$.
Moreover,
\begin{equation}
\label{eq: finite case observable detects distance}
\norm{F(x)-F(x')}_\infty<\tau/2
\quad\Longrightarrow\quad
\mathbf{d}(x,x')<\tau.
\end{equation}
Indeed, choose $x_i$ with $\mathbf{d}(x,x_i)<\tau/4$.
The assumption implies
$\min\{\mathbf{d}(x',x_i),1\}<3\tau/4<1$,
so $\mathbf{d}(x',x_i)<3\tau/4$ and hence
$\mathbf{d}(x,x')<\tau$.
Applying this at each time gives
\begin{equation}
\label{eq: finite case orbit observable detects distance}
\norm{F_N(x)-F_N(x')}_\infty<\tau/2
\quad\Longrightarrow\quad
\mathbf{d}_N(x,x')<\tau.
\end{equation}

\smallskip
\noindent
\emph{Step 2: A uniform lower bound for the radius.}
Fix $0<\delta<\tau/8$.
Apply Proposition \ref{proposition: polyhedral modeling}
to $F$, $\delta$, and $\eta=1$, and let $\varepsilon>0$
and $L\in\mathbb{N}$ be given by the proposition.
For $N\in L\mathbb{N}$, take the resulting maps
\[
f_N\colon X\to P_N,\qquad
g_N\colon P_N\to Q_N,\qquad
\Phi_N\colon P_N\to[0,1]^{\ell N}.
\]
Set
\begin{equation}
\label{eq: finite case waist parameters}
d=\lfloor\alpha N\rfloor,\qquad
k=\lfloor(\beta-\alpha)N\rfloor,\qquad
s=(a-b)N.
\end{equation}

For all sufficiently large $N\in L\mathbb{N}$, 
we claim that
\begin{equation}
\label{eq: finite case radius lower bound}
\mathcal{R}^{\infty}_{d+k}(\Phi_N)\geq\tau/8.
\end{equation}
Suppose to the contrary 
that there are a finite simplicial
complex $K$ with $\dim K\leq d+k$ and continuous maps
\[
\varphi\colon P_N\to K,\qquad
\psi\colon K\to\mathbb{R}^{\ell N}
\]
such that
\[
\sup_{p\in P_N}
\norm{\Phi_N(p)-(\psi\circ\varphi)(p)}_\infty<\tau/8.
\]
Since $\norm{\Phi_N\circ f_N(x)-F_N(x)}_\infty<\delta$
for every $x\in X$, we have
\[
\norm{F_N(x)-(\psi\circ\varphi\circ f_N)(x)}_\infty
<\delta+\tau/8<\tau/4.
\]

Thus, if $(\varphi\circ f_N)(x)=(\varphi\circ f_N)(x')$,
then
\[
\norm{F_N(x)-F_N(x')}_\infty<\tau/2,
\]
and hence $\mathbf{d}_N(x,x')<\tau$ by
\eqref{eq: finite case orbit observable detects distance}.
This shows that $\varphi\circ f_N\colon X\to K$
is a $\tau$-embedding with respect to $\mathbf{d}_N$.
Consequently,
\[
\widim_\tau(X,\mathbf{d}_N)
\leq\dim K\leq d+k\leq\beta N,
\]
contradicting \eqref{eq: finite case width lower bound}.
This proves \eqref{eq: finite case radius lower bound}.

\smallskip
\noindent\emph{Step 3: Applying the waist inequality.}
Since $\varepsilon$ is fixed while $N$ tends to infinity,
\[
\lim_{\substack{N\to\infty\\N\in L\mathbb{N}}}
\frac{\dim P_N}{N}
=
\lim_{N\to\infty}
\frac{\widim_\varepsilon(X,\mathbf{d}_N)}{N}
\leq D<2\alpha.
\]
Therefore $\dim P_N\leq 2d-1$ for all sufficiently large
$N\in L\mathbb{N}$. (This is the crucial point where 
our assumption $D<\infty$ works.)
Also, $3\leq d<\ell N$ for all sufficiently large $N$,
since $\ell>D>\alpha$.

Since $d = \lfloor \alpha N\rfloor$, 
$s = (a-b)N$ and $\dim Q_N<bN$, we obtain
\[
d-\dim Q_N-2-s>(\alpha-a)N-3>0
\]
for all sufficiently large $N$.
In particular, $\dim Q_N \leq d-3$ and 
$0<s<d-\dim Q_N-2$.
Furthermore,
\[
\kappa :=1-\frac{s+2}{d-\dim Q_N}
=\frac{d-\dim Q_N-2-s}{d-\dim Q_N}
\geq \frac{(\alpha-a)N-3}{\alpha N} 
\geq\frac{\alpha-a}{2\alpha}>0
\]
for all sufficiently large $N$.
Since $k = \lfloor (\beta-\alpha)N\rfloor$,
\[
\sqrt{\frac{k+1}{\ell N}}
\geq\sqrt{\frac{\beta-\alpha}{\ell}}.
\]
All the hypotheses of
the $\ell^\infty$-version of the waist inequality
(Corollary 
\ref{corollary: ell^infty waist inequality})
are therefore satisfied with $n=\ell N$ and $r=\delta$.
Combining that corollary with
\eqref{eq: finite case radius lower bound} gives
\begin{equation}
\label{eq: finite case covering lower bound}
\sup_{y\in Q_N}
\#\bigl(\Phi_N(g_N^{-1}(y)),
\norm{\cdot}_\infty,\delta\bigr)
\geq
\frac{1}{\dim Q_N +1}
\left(\frac{C}{\delta}\right)^{(a-b)N},
\end{equation}
where
\[
C:=\frac{\tau}{8}\,
\frac{\alpha-a}{2\alpha}\,
\sqrt{\frac{\beta-\alpha}{\ell}}>0
\]
is independent of both $N$ and $\delta$.
On the other hand,
Proposition \ref{proposition: polyhedral modeling} 
(with $\eta = 1$) gives
\[
\sup_{y\in Q_N}
\#\bigl(\Phi_N(g_N^{-1}(y)),
\norm{\cdot}_\infty,\delta\bigr)
\leq\exp\bigl((S_\delta(\pi)+1)N\bigr).
\]
Since $\dim Q_N+1 < bN+1$, 
taking logarithms and dividing by $N$
yields
\[
S_\delta(\pi)+1
\geq(a-b)\log(C/\delta)-\frac{\log(bN+1)}{N}.
\]
For each fixed $\delta$, let $N\to\infty$ 
through multiples
of $L$ to obtain
\[
S_\delta(\pi)+1\geq(a-b)\log(C/\delta).
\]
Finally, divide by $\log(1/\delta)$ 
and take the lower limit
as $\delta\to0$.
Since $C>0$ is independent of $\delta$, 
we conclude that
\[
\lmdimM(\pi,T,\mathbf{d})
=\liminf_{\delta\to0}
\frac{S_\delta(\pi)}{\log(1/\delta)}
\geq a-b.
\]
As $a<D$ and $b>\mdim(Y,S)$ were arbitrary, 
this proves
Theorem \ref{theorem: main theorem} 
when $\mdim(X,T)<\infty$.
\end{proof}

\subsection{The case of $\mdim(X,T)=\infty$}
\label{subsection: the case of infinite mean dimension}

In this subsection, we assume that 
$\pi\colon (X, T)\to (Y, S)$ is an 
equivariant continuous map between dynamical systems
such that 
\[ \mdim(X,T) = \infty, \quad 
   \mdim(Y, S) < \infty. \]
We take a compatible metric $\mathbf{d}$ on $X$. 
We would like to prove that
$\lmdimM(\pi,T,\mathbf{d})=\infty$.
We reduce this to the finite mean dimension case 
by constructing 
an intermediate factor of finite but arbitrarily
large mean dimension,
equipped with a compatible metric for
which the lower conditional metric mean dimension 
does not
exceed that of $\pi$.

\begin{proof}[Completion of the proof 
in the infinite mean dimension case]
Let $a$ be an arbitrary positive number.
We will construct a dynamical system 
$(Z, R)$ with a compatible metric 
$\mathbf{D}$ on $Z$ 
and equivariant continuous maps
\[
(X,T)\xrightarrow{\Psi}(Z,R)
\xrightarrow{\widetilde{\pi}}(Y,S),
\qquad \widetilde{\pi}\circ\Psi=\pi,
\]
such that
\begin{equation}
\label{eq: infinite case factor bounds}
a<\mdim(Z,R)<\infty,
\qquad
\lmdimM(\widetilde{\pi},R,\mathbf{D})
\leq\lmdimM(\pi,T,\mathbf{d}).
\end{equation}

\smallskip
\noindent\emph{Step 1: Constructing an intermediate factor.}
Choose $0<\tau<1$ such that
\begin{equation}
\label{eq: infinite case width lower bound}
\lim_{N\to\infty}
\frac{\widim_\tau(X,\mathbf{d}_N)}{N}>a.
\end{equation}
As in Step~1 of the proof in
\S\ref{subsection: the case of finite mean dimension},
choose a finite $\tau/4$-spanning set
$\{x_1,\ldots,x_\ell\}$ in $X$ and define
\[
F(x)=\bigl(\min\{\mathbf{d}(x,x_i),1\}\bigr)_{i=1}^\ell.
\]
Then $F\colon (X,\mathbf{d})\to
([0,1]^\ell,\norm{\cdot}_\infty)$ is $1$-Lipschitz and
\begin{equation}
\label{eq: infinite case observable detects distance}
\norm{F(x)-F(x')}_\infty<\tau/2
\quad\Longrightarrow\quad
\mathbf{d}(x,x')<\tau.
\end{equation}

Let $\Sigma=([0,1]^\ell)^{\mathbb{Z}}$, and let
$\sigma\colon\Sigma\to\Sigma$ be the shift
$\sigma(u)_j=u_{j+1}$.
Define
\[
\Psi(x)=\left(\pi(x),\bigl(F(T^j x)\bigr)_{j\in\mathbb{Z}}\right)
\in Y\times\Sigma,
\qquad Z=\Psi(X).
\]
Then $Z$ is a closed $(S\times\sigma)$-invariant subset of
$Y\times\Sigma$.
Set $R=(S\times\sigma)|_Z$, and let
$\widetilde{\pi}\colon Z\to Y$ be the restriction of the
coordinate projection.
These maps satisfy $\widetilde{\pi}\circ\Psi=\pi$.

Fix a compatible metric $\mathbf{d}^\prime$ on $Y$.
For $(y,u),(y',u')\in Z$, set
\begin{equation}
\label{eq: infinite case factor metric}
\mathbf{D}\bigl((y,u),(y',u')\bigr)
=
\max\left\{
\mathbf{d}^\prime(y,y'),\,
\sup_{j\in\mathbb{Z}}2^{-|j|}\norm{u_j-u'_j}_\infty
\right\}.
\end{equation}
This is a compatible metric on $Z$.
We write $\mathbf{D}_N$ and $\mathbf{d}^\prime_{N}$ for the
Bowen metrics associated with $R$ and $S$, respectively.

\smallskip
\noindent
\emph{Step 2: Bounding the 
mean dimension of the factor.}
By construction, for all $N\geq 1$
\[
\norm{F_N(x)-F_N(x')}_\infty
\leq\mathbf{D}_N\bigl(\Psi(x),\Psi(x')\bigr).
\]
Consequently, composing any $(\tau/2)$-embedding of
$(Z,\mathbf{D}_N)$ into a finite simplicial 
complex with $\Psi$
gives a $\tau$-embedding of $(X,\mathbf{d}_N)$, by
\eqref{eq: infinite case observable detects distance}.
Thus
\[
\widim_\tau(X,\mathbf{d}_N)
\leq\widim_{\tau/2}(Z,\mathbf{D}_N).
\]
Together with \eqref{eq: infinite case width lower bound},
this implies
\begin{equation}
\label{eq: infinite case factor mean dimension lower bound}
\mdim(Z,R)>a.
\end{equation}

For the upper bound, we have
\begin{equation}
\label{eq: infinite case factor mean dimension upper bound}
\mdim(Z,R)\leq 
\mdim\left(Y\times \Sigma, S\times \sigma\right) \leq
\mdim(Y,S)+\ell<\infty.
\end{equation}
Here we have used 
$\mdim\left(Y\times \Sigma, S\times \sigma\right) \leq
\mdim(Y,S)+\mdim\left(\Sigma, \sigma\right)$; 
see \cite[\S 2]{Lindenstrauss}.

\smallskip
\noindent
\emph{Step 3: Comparing the conditional 
metric mean dimensions.}
Fix $r>0$, and choose $m\geq1$ with $2^{-m}<r$.
We claim that, for every $N\geq1$,
\begin{equation}
\label{eq: infinite case fiber covering comparison}
\sup_{y\in Y}
\#\bigl(\widetilde{\pi}^{-1}(y),\mathbf{D}_N,r\bigr)
\leq
\sup_{y\in Y}
\#\bigl(\pi^{-1}(y),\mathbf{d}_{N+2m},r\bigr).
\end{equation}
Indeed, if $x,x'\in\pi^{-1}(y)$, the $Y$-coordinates of
$R^i\Psi(x)$ and $R^i\Psi(x')$ coincide for every $i$.
Then
\begin{equation}
\label{eq: infinite case padded orbit estimate}
 \begin{split}
\mathbf{D}_N\bigl(\Psi(x),\Psi(x')\bigr)
 & = \max_{0\leq i <N} \sup_{j\in \mathbb{Z}}
   2^{-|j|}\norm{F(T^{i+j}x)-F(T^{i+j}x^\prime)}_\infty
   \\
 & \leq
\max\left\{
\mathbf{d}_{[-m, N+m)}(x,x'),\,2^{-m}
\right\},
 \end{split}
\end{equation}
where 
$\mathbf{d}_{[-m, N+m)}(x,x') 
:= \max_{-m \leq n < N+m}
\mathbf{d}(T^n x, T^n x^\prime)$.
Here we have used 
the $1$-Lipschitz property of $F$ for $|j|\leq m$
and the fact that $F$ takes values in $[0,1]^\ell$ for the
remaining coordinates.

Take an open cover $U_1,\ldots,U_M$ of
$\pi^{-1}(y)$ whose members have
$\mathbf{d}_{[-m, N+m)}$-diameter less than $r$, 
with
\[
M=\#\bigl(\pi^{-1}(y),\mathbf{d}_{[-m, N+m)},r\bigr).
\]
Since $\Psi$ is surjective, 
$\Psi(\pi^{-1}(y))=\widetilde{\pi}^{-1}(y)$ and 
hence
the sets $\Psi(U_i)$ cover
$\widetilde{\pi}^{-1}(y)$.
By \eqref{eq: infinite case padded orbit estimate},
each of them has $\mathbf{D}_N$-diameter strictly less than $r$.
Although these images need not be open, they can be slightly
enlarged to open subsets of $Z$ of
$\mathbf{D}_N$-diameter less than $r$.
It follows that
\begin{align*}
\#\bigl(\widetilde{\pi}^{-1}(y),\mathbf{D}_N,r\bigr)
&\leq
\#\bigl(\pi^{-1}(y),\mathbf{d}_{[-m, N+m)},r\bigr) \\
&=\#\bigl(\pi^{-1}(S^{-m}y),\mathbf{d}_{N+2m},r\bigr).
\end{align*}
In the last step, we used the fact that 
$T^{-m}\colon (\pi^{-1}(y), \mathbf{d}_{[-m, N+m)})
\to (\pi^{-1}(S^{-m}y),\mathbf{d}_{N+2m})$ is isometric.
Taking the supremum over $y\in Y$, we obtain
\eqref{eq: infinite case fiber covering comparison}.

The integer $m$ is independent of $N$.
Hence, taking logarithms, dividing by $N$, and letting
$N\to\infty$ in
\eqref{eq: infinite case fiber covering comparison} gives
\[
S_r(\widetilde{\pi})\leq S_r(\pi),
\]
where $S_r(\widetilde{\pi})$ is defined using $\mathbf{D}$.
Dividing by $\log(1/r)$ 
and taking the lower limit as $r\to0$
yields
\begin{equation}
\label{eq: infinite case conditional metric comparison}
\lmdimM(\widetilde{\pi},R,\mathbf{D})
\leq\lmdimM(\pi,T,\mathbf{d}).
\end{equation}
This completes the proof of
\eqref{eq: infinite case factor bounds}.

We can now apply the finite mean dimension case of
Theorem \ref{theorem: main theorem}, proved in
\S\ref{subsection: the case of finite mean dimension},
to $\widetilde{\pi}\colon(Z,R)\to(Y,S)$ and $\mathbf{D}$.
We obtain
\[
\begin{aligned}
a<\mdim(Z,R)
&\leq\mdim(Y,S)
  +\lmdimM(\widetilde{\pi},R,\mathbf{D})\\
&\leq \mdim(Y,S)+\lmdimM(\pi,T,\mathbf{d}).
\end{aligned}
\]
Since $\mdim(Y, S)<\infty$ and $a>0$ is arbitrary, 
this proves
$\lmdimM(\pi,T,\mathbf{d})=\infty$
and completes the proof of Theorem 
\ref{theorem: main theorem}.
\end{proof}

\section{Two main propositions and the proof of
Theorem \ref{theorem: new waist inequality}}
\label{section: two main proposition and
the proof of waist inequality}

In \S \ref{subsection: main propositions},
we state two main propositions
and use them to prove our new waist inequality
(Theorem \ref{theorem: new waist inequality})
in \S \ref{subsection: proof of new waist inequality}.
The two propositions are proved later,
in \S \ref{section: proof of proposition avoidance}
and \S
\ref{section: proof of proposition low dimensional approximation},
respectively.

In \S \ref{subsection: weaker waist inequality},
we derive a weaker waist inequality
and use it to give an alternative proof of
Theorem \ref{theorem: main theorem}.
This subsection is not needed for the rest of the paper
and may be skipped.

\subsection{Main propositions}
\label{subsection: main propositions}

It will be convenient below 
to work with \emph{locally finite}
simplicial complexes, because we will consider
the $d$-skeleton of a cubical grid in $\mathbb{R}^n$,
viewed as a simplicial complex after triangulation.
The use of locally finite complexes is only for
convenience of presentation:
For our application to the waist inequality,  
it would be possible to work with 
finite portions of the grids 
by restricting them to sufficiently large boxes. 
Thus, readers who prefer to work only 
with finite simplicial
complexes may do so throughout.

A simplicial complex is said to be \emph{locally finite}
if each vertex belongs to only finitely many simplices.
We identify a simplicial complex 
with its geometric realization,
equipped with the weak topology: 
a subset is closed if and only
if its intersection with every closed simplex 
is closed in that simplex.
There is no need to assume 
that the vertex set is countable,
although only the countable case 
will be needed in our arguments.

Recall that we have denoted the standard Euclidean norm on 
$\mathbb{R}^n$ by $\norm{x}$ and used 
\[  \rho_n(x) = \frac{\norm{x}}{\sqrt{n}}, \quad  
     (x\in \mathbb{R}^n). \]

\begin{proposition}[Quantitative avoidance]
  \label{proposition: quantitative avoidance}
  Let $m, n$ be positive integers. 
  Let $a$ be an integer and let $s$ be a real number 
  with $3\leq a \leq n-1$ and
  $0<s < a-2$.
  Set 
  \[ \kappa := 1 - \frac{s+2}{a} >0. \]   
  Let $K$ be a locally finite simplicial complex 
  satisfying
  \[ \dim K \leq n-a-1. \]
  Let
  $f, p_1,p_2, \dots, p_m\colon K\to \mathbb{R}^n$ be 
  continuous maps. If positive numbers $r$ and $R$
  satisfy 
  \[ m^{1/s}r < \kappa R, \]
  then there is a continuous map 
  $g\colon K\to \mathbb{R}^n$ such that, for every 
  $x\in K$ and every $i = 1,2,\dots, m$
  \[ \rho_n\left(g(x)-f(x)\right) \leq R, \qquad 
      \rho_n\left(g(x)-p_i(x)\right) > r. \]
\end{proposition}

\begin{figure}[htbp]
\centering
\includegraphics[width=0.75\textwidth]
{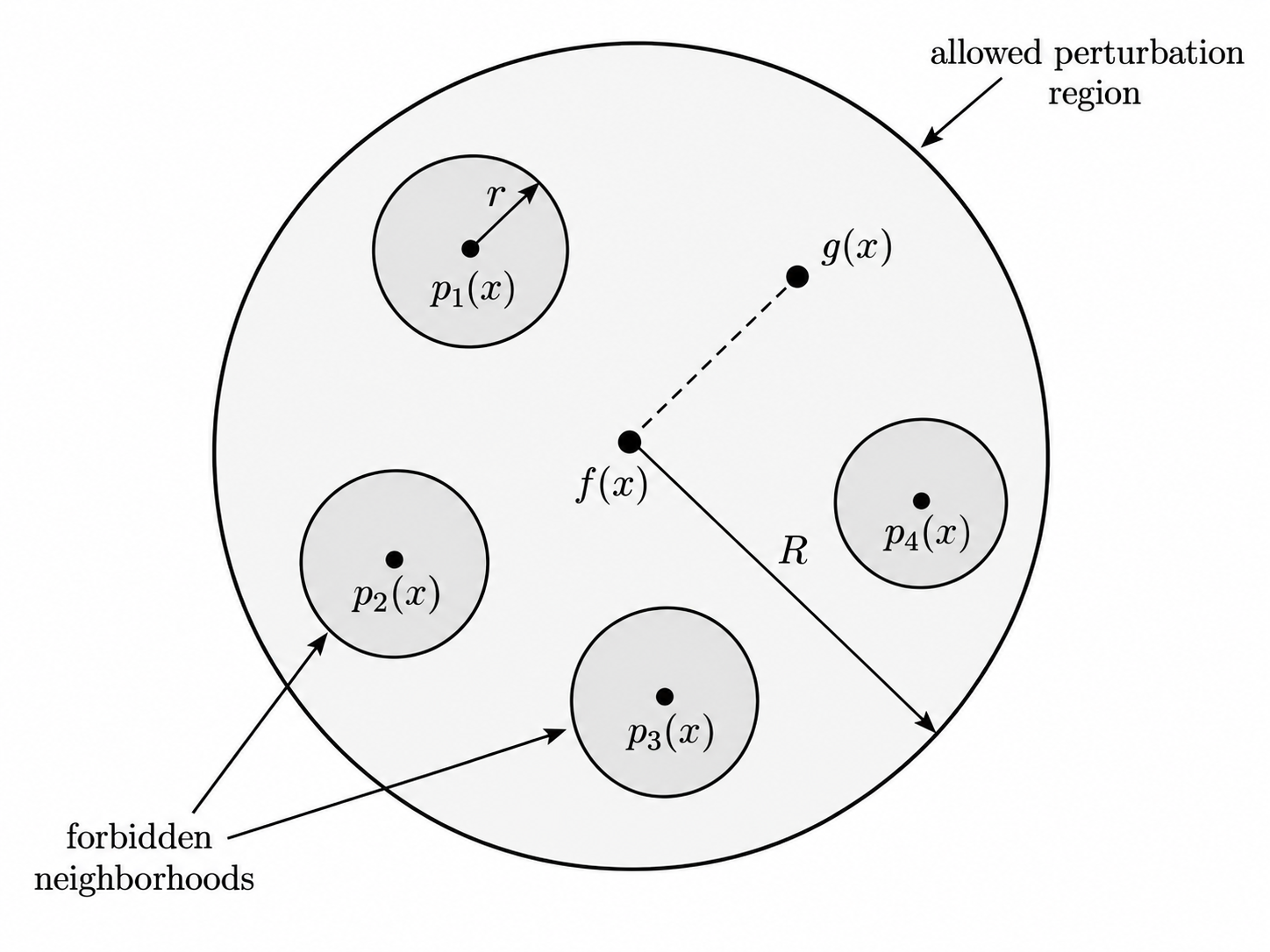}
\caption{
Schematic illustration of quantitative avoidance.
For each $x\in K$, we seek a point $g(x)$
within $\rho_n$-distance $R$ of $f(x)$
that stays at $\rho_n$-distance greater than $r$
from all the obstacles $p_1(x),\ldots,p_m(x)$.
The point $g(x)$ must be chosen continuously in $x$.
}
\label{fig: quantitative avoidance}
\end{figure}

Roughly speaking, this proposition says 
that a continuous map $f$
from a sufficiently low-dimensional simplicial complex 
can be perturbed to avoid finitely many 
\lq\lq{}moving obstacles\rq\rq{} $p_i$
while keeping the perturbation small.
The key feature is its quantitative nature: 
for fixed $a$ and $s$, 
the number $m$ of obstacles can grow 
on the order of $(R/r)^s$, 
where $R$ is the allowed perturbation size 
and $r$ is the required distance from the obstacles.
Figure~\ref{fig: quantitative avoidance} illustrates 
the basic idea of the proposition.

Fix $h>0$.
The \emph{cubical grid of mesh size $h$} in $\mathbb{R}^n$
consists of all the cubes
\[
\prod_{i=1}^n [hm_i,h(m_i+1)],
\qquad (m_1,\ldots,m_n)\in\mathbb{Z}^n,
\]
together with all their faces.
Here a face is obtained by replacing some of the interval
factors $[hm_i, h(m_i+1)]$
with one of their endpoints $hm_i, h(m_i+1)$.
These cubes and their faces form a \emph{cubical complex}:
the collection is closed under taking faces,
and the intersection of any two members is either empty
or a common face.
The vertex set of this cubical complex is 
$h\mathbb{Z}^n$.

For an integer $0\leq d\leq n$, the \emph{$d$-skeleton}
$K(n,d,h)$ of this grid is the union of all faces
of dimension at most $d$.
Equivalently,
\[
K(n,d,h)
=
\left\{
x\in\mathbb{R}^n \;\middle|\;
\#\{i:x_i\in h\mathbb{Z}\}\geq n-d
\right\}.
\]
We regard $K(n,d,h)$ as a locally finite simplicial complex
by choosing a triangulation compatible with
the cubical structure.\footnote{
Such a triangulation can be obtained by subdividing
the cubes into finitely many simplices
in a way that agrees on common faces.
The resulting simplicial complex is locally finite,
has dimension $d$, and carries the subspace topology
inherited from $\mathbb{R}^n$.}

\begin{figure}[htbp]
\centering
\includegraphics[width=0.78\textwidth]
{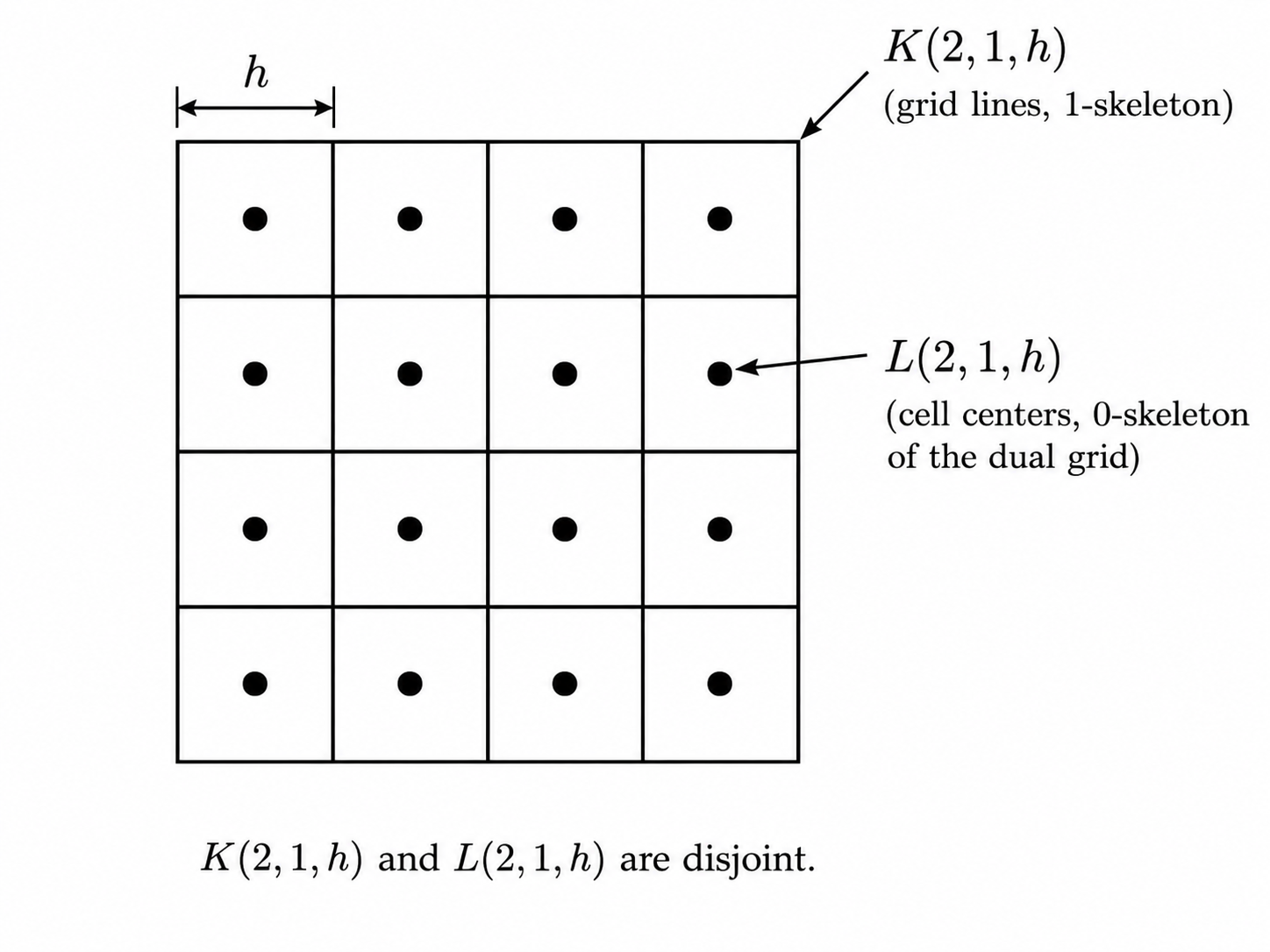}
\caption{
The $1$-skeleton $K(2,1,h)$ of the cubical 
grid in $\mathbb{R}^2$
and the dual skeleton $L(2,1,h)$.
Here $K(2,1,h)$ consists of the grid lines, while
$L(2,1,h)$ consists of the centers of the square cells
(i.e.\ the $0$-skeleton of the dual grid).
In particular, $K(2,1,h)$ and $L(2,1,h)$ are disjoint.
}
\label{fig:grid-dual-skeleton}
\end{figure}

The \emph{dual grid} is obtained by translating
the original cubical grid by $(h/2,\ldots,h/2)$.
Thus its vertices are the centers of the
$n$-dimensional cubes of the original grid.
For $0\leq d<n$, let $L(n,d,h)$ be the
$(n-d-1)$-skeleton of the dual grid:
\[
L(n,d,h)
=
\left\{
x\in\mathbb{R}^n \;\middle|\;
\#\left\{
i:x_i\in h\left(\mathbb{Z}+\frac12\right)
\right\}\geq d+1
\right\}.
\]
The two skeleta $K(n,d,h)$ and $L(n,d,h)$
are disjoint.
Indeed, a point in $K(n,d,h)$ has at least $n-d$
coordinates in $h\mathbb{Z}$,
whereas a point in $L(n,d,h)$ has at least $d+1$
coordinates in $h(\mathbb{Z}+1/2)$.
Since these two sets of coordinate values are disjoint
and $(n-d)+(d+1)>n$,
the two conditions cannot hold simultaneously.
Figure~\ref{fig:grid-dual-skeleton} illustrates
the configuration of $K(n,d,h)$ and $L(n,d,h)$ 
in the case
$n=2$ and $d=1$.

\begin{proposition}[From avoidance to low-dimensional 
  approximation]
  \label{proposition: low dimensional approximation}
  Let $1\leq d <n$ be integers, and
  let $h, R$ be positive real numbers.
  Let $P$ be a finite simplicial complex satisfying
  \[ \dim P \leq 2d -1. \]
  Let $\Phi\colon P\to \mathbb{R}^n$ and 
  $F\colon P\times L(n,d,h) \to \mathbb{R}^n$
  be continuous maps. Suppose that
  \[ \rho_n\left(F(x, u) -u\right) \leq R, \qquad 
      F(x, u) \neq \Phi(x) \]
  for all $x\in P$ and $u\in L(n,d,h)$.
  Then there exists a continuous map 
  $\varphi\colon P\to K(n,d,h)$ satisfying 
  \[ \sup_{x\in P} 
    \rho_n\left(\varphi(x)-\Phi(x)\right)
    < R + 4h. \]
\end{proposition}

This proposition converts 
an avoidance property into a low-dimensional approximation. 
For each $x\in P$, 
the map $u\mapsto F(x,u)$ is a perturbation of 
the inclusion of the dual skeleton $L(n,d,h)$ 
into $\mathbb{R}^n$, 
and its image avoids the point $\Phi(x)$.
The conclusion is that we can move $\Phi(x)$ 
into the original $d$-skeleton $K(n,d,h)$, 
with a controlled displacement, 
in a way that depends continuously on $x$. 
The essential point is the continuity of this choice, 
rather than the existence of a nearby point of 
$K(n,d,h)$ for each individual $x$.

Let us briefly explain how these two propositions 
fit together in the proof of 
Theorem \ref{theorem: new waist inequality}.
The quantitative avoidance proposition 
will be applied to a simplicial complex 
of the form $Q\times L(n,d,h)$. 
Setting $a=d-\dim Q$, we have
\[
\begin{aligned}
\dim\bigl(Q\times L(n,d,h)\bigr)
&\leq \dim Q+n-d-1\\
&=n-a-1.
\end{aligned}
\]
Moreover, the assumption $\dim Q\leq d-3$ 
is equivalent to $a\geq3$. 
Thus the dimensional conditions of 
quantitative avoidance are naturally satisfied.
The resulting avoidance map will be used to 
construct the map $F$ required 
in Proposition 
\ref{proposition: low dimensional approximation}. 
The latter proposition then produces a map from $P$ 
into the $d$-skeleton of the cubical grid, 
yielding the desired low-dimensional approximation.

We prove Propositions \ref{proposition: quantitative avoidance}
and \ref{proposition: low dimensional approximation}
in \S \ref{section: proof of proposition avoidance} and 
\S \ref{section: proof of proposition low dimensional approximation},
respectively.

\subsection{Proof of Theorem
\ref{theorem: new waist inequality}}
\label{subsection: proof of new waist inequality}

We now deduce Theorem \ref{theorem: new waist inequality}
from Propositions \ref{proposition: quantitative avoidance}
and \ref{proposition: low dimensional approximation}.
For the reader's convenience, we first restate the theorem.

\begin{restatedwaist}[New waist inequality]
Let $3\leq d<n$ be integers.
Let $P$ and $Q$ be finite simplicial complexes, and let
$f\colon P\to Q$ and $\Phi\colon P\to\mathbb{R}^n$
be continuous maps.
Suppose that
\[
\dim P\leq 2d-1,
\qquad
\dim Q\leq d-3.
\]
Let $s$ be a real number satisfying
$0<s<d-\dim Q-2$, and set
\[
\kappa=1-\frac{s+2}{d-\dim Q} >0.
\]
Then, for every $r>0$,
\[
\sup_{y\in Q}
\#\bigl(\Phi(f^{-1}(y)),\rho_n,r\bigr)
\geq
\frac{1}{\dim Q+1}
\left(\frac{\kappa\mathcal{R}_d(\Phi)}{r}\right)^s.
\]
\end{restatedwaist}

\begin{proof}
If $\mathcal{R}_d(\Phi)=0$, the assertion is immediate.
Hence we assume that $\mathcal{R}_d(\Phi)>0$.
Fix $r>0$, and put
\[
q=\dim Q,
\qquad
M=\sup_{y\in Q}
\#\bigl(\Phi(f^{-1}(y)),\rho_n,r\bigr).
\]
For $z\in \mathbb{R}^n$ and $r>0$,
we write
\[
B_{\rho_n}(z,r)=
\{w\in\mathbb{R}^n \mid \rho_n(w-z)<r\}.
\]

\smallskip
\noindent
\emph{Step 1: Covering nearby fibers.}
For each $y\in Q$, choose points
$z_{y,1},\ldots,z_{y,M}\in\mathbb{R}^n$ such that
\[
\Phi(f^{-1}(y))
\subset
\bigcup_{j=1}^M B_{\rho_n}(z_{y,j},r).
\]
(For an empty fiber,
choose the centers arbitrarily.)
Let $W_y$ denote the union of these balls.
Since $P$ is compact, we can take 
an open neighborhood $V_y$ of $y$ satisfying
\[
\Phi(f^{-1}(V_y))\subset W_y.
\]
Choose a finite subcover $\mathcal{V}$ of 
$\{V_y\mid y\in Q\}$.
Thus, each 
$\Phi\left(f^{-1}(V)\right)$, $(V\in \mathcal{V})$,
is covered by $M$ balls of radius $r$.

\smallskip
\noindent
\emph{Step 2: Constructing continuous moving centers.}
After sufficiently many barycentric subdivisions of $Q$,
we may assume that every open star
$O_Q(v)$ $(v\in \mathrm{Ver}(Q))$ 
is contained
in a member of $\mathcal{V}$.
Moreover, the vertices can be colored with $q+1$ colors
so that the vertices of each simplex have distinct colors:
in the last barycentric subdivision, color the barycenter
of each face of the preceding subdivision by its dimension.
See Figure~\ref{fig:barycentric-coloring-q2}.
Write 
\[
c\colon\mathrm{Ver}(Q)\to\{0,1,\ldots,q\}
\]
for this coloring.

\begin{figure}[htbp]
  \centering
  \includegraphics[width=0.8\textwidth]
  {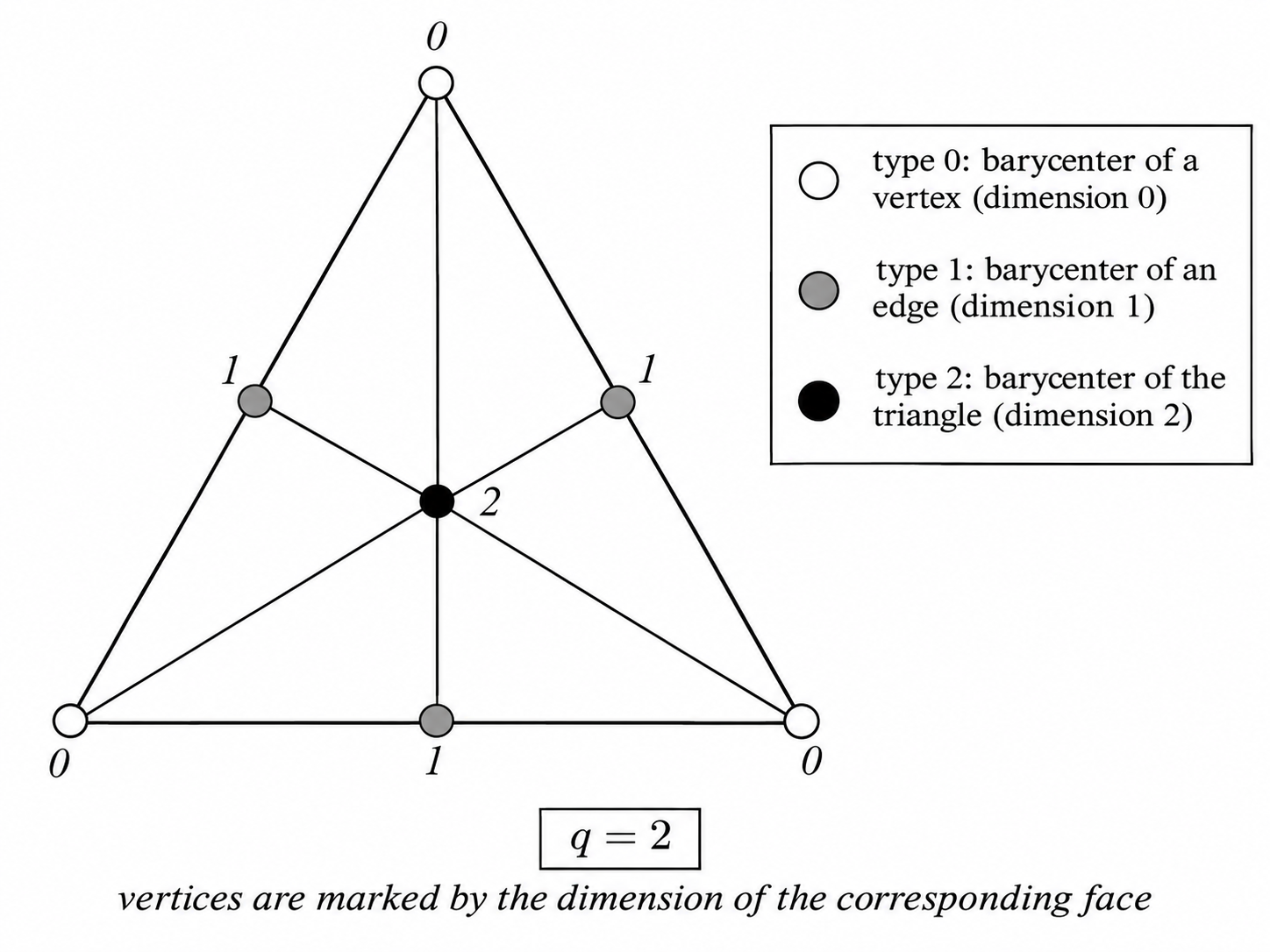}
  \caption{
    Coloring of the vertices in the barycentric subdivision for $q=2$.
    Each vertex is colored according to the dimension of the corresponding face
    in the preceding subdivision: type $0$ for vertices, type $1$ for edges,
    and type $2$ for the $2$-simplex.
  }
  \label{fig:barycentric-coloring-q2}
\end{figure}

For each vertex $v \in \mathrm{Ver}(Q)$,
choose points $z_{v,1},\ldots,z_{v,M}$ such that
\begin{equation}
\label{eq: waist local centers}
\Phi(f^{-1}(O_Q(v)))
\subset
\bigcup_{j=1}^M B_{\rho_n}(z_{v,j},r).
\end{equation}
This is possible by Step~1.

We take continuous functions 
$\lambda_v \colon Q\to [0,1]$, 
$(v\in \mathrm{Ver}(Q))$, such that 
 \begin{itemize}
   \item $\supp (\lambda_v) \subset O_Q(v)$,
   \item for any $y\in Q$, there exists a vertex $v$ with 
         $\lambda_v(y) =1$.
 \end{itemize}
For $0\leq c\leq q$ and $1\leq j\leq M$, define
continuous maps $p_{c,j}\colon Q\to\mathbb{R}^n$ by
\begin{equation}
\label{eq: waist moving centers definition}
p_{c,j}(y)=
\sum_{\substack{v\in\mathrm{Ver}(Q)\\c(v)=c}}
\lambda_v(y)z_{v,j}.
\end{equation}
At each point $y\in Q$, at most one summand in this sum
has a nonzero coefficient, since vertices of the same color
cannot belong to a common simplex.

We claim that, for every $x\in P$,
\begin{equation}
\label{eq: waist moving centers cover}
\min_{\substack{0\leq c\leq q\\1\leq j\leq M}}
\rho_n\bigl(\Phi(x)-p_{c,j}(f(x))\bigr)<r.
\end{equation}
To see this, put $y=f(x)$.
Take a vertex $v$ satisfying $\lambda_v(y)=1$.
Then $y\in O_Q(v)$.
For all other vertices $w$ of color $c(v)$, 
we have $O_Q(w)\cap O_Q(v) =\emptyset$ and hence
$\lambda_w(y)=0$. Therefore,
\[
p_{c(v),j}(y)=z_{v,j}
\qquad(1\leq j\leq M).
\]
Now \eqref{eq: waist moving centers cover} follows from
\eqref{eq: waist local centers}.
Thus $(q+1)M$ continuous moving centers suffice to cover
the images of all fibers by balls of radius $r$.

\smallskip
\noindent\emph{Step 3: Applying quantitative avoidance.}
Put
\[
a=d-q,
\qquad
m=(q+1)M.
\]
Then $3\leq a\leq n-1$, $0<s<a-2$, and
$\kappa=1-(s+2)/a$.
Choose any $R>0$ satisfying
\begin{equation}
\label{eq: waist avoidance radius}
m^{1/s}r<\kappa R,
\end{equation}
and let $h>0$ be arbitrary.
The product $Q\times L(n,d,h)$ admits a locally finite
triangulation with
\[
\dim\bigl(Q\times L(n,d,h)\bigr)
\leq q+n-d-1=n-a-1.
\]
We apply Proposition 
\ref{proposition: quantitative avoidance}
to this product and the map
\[ Q\times L(n,d,h) 
   \ni (y,u)\mapsto u \in \mathbb{R}^n, \]
and the $m$ maps
\[    Q\times L(n,d,h) 
   \ni (y,u)\mapsto p_{c,j}(y)\in \mathbb{R}^n. \]
Then we obtain a continuous map
\[
G\colon Q\times L(n,d,h)\to\mathbb{R}^n
\]
satisfying
\[
\rho_n(G(y,u)-u)\leq R,
\qquad
\rho_n(G(y,u)-p_{c,j}(y))>r
\]
for all $(y,u)\in Q\times L(n,d,h)$ and all $c,j$.

Define
\[
F\colon P\times L(n,d,h)\to\mathbb{R}^n,
\qquad
F(x,u)=G(f(x),u).
\]
Then $F$ is continuous and
$\rho_n(F(x,u)-u)\leq R$.
Also, \eqref{eq: waist moving centers cover} implies
$F(x,u)\neq\Phi(x)$ for every $x\in P$ and $u\in L(n,d,h)$:
otherwise $G(f(x),u)$ would be within distance less than $r$
of one of the points $p_{c,j}(f(x))$,
contrary to its construction.

\smallskip
\noindent
\emph{Step 4: Low-dimensional approximation.}
Since $\dim P\leq2d-1$, we can apply
Proposition \ref{proposition: low dimensional approximation}
to $\Phi$ and $F$.
It gives a continuous map
$\varphi\colon P\to K(n,d,h)$ such that
\[
\sup_{x\in P}\rho_n(\varphi(x)-\Phi(x))<R+4h.
\]
The image $\varphi(P)$ is compact and hence is contained
in a finite subcomplex of the locally finite complex
$K(n,d,h)$.
This subcomplex has dimension at most $d$.
Therefore, by the definition of the $d$-radius,
\[
\mathcal{R}_d(\Phi)<R+4h.
\]
Since $h>0$ is arbitrary, we obtain
$\mathcal{R}_d(\Phi)\leq R$.
Finally, letting $R$ decrease to $m^{1/s}r/\kappa$ in
\eqref{eq: waist avoidance radius} yields
\[
\mathcal{R}_d(\Phi)
\leq\frac{((q+1)M)^{1/s}r}{\kappa}.
\]
Rearranging this inequality gives
\[
M\geq\frac{1}{q+1}
\left(\frac{\kappa\mathcal{R}_d(\Phi)}{r}\right)^s,
\]
which is the desired conclusion.
\end{proof}

\subsection{A weaker waist inequality
and another proof of Theorem \ref{theorem: main theorem}}
\label{subsection: weaker waist inequality}

This subsection is optional and
may be skipped entirely.

The most technically demanding part of the proof of
our new waist inequality
(Theorem \ref{theorem: new waist inequality})
is the proof of Proposition
\ref{proposition: low dimensional approximation},
which will be given in \S
\ref{section: proof of proposition low dimensional approximation}.
It occupies more than one third of the paper.
Here we explain that
a weaker waist inequality can be derived
from Proposition \ref{proposition: quantitative avoidance}
alone,
without using Proposition
\ref{proposition: low dimensional approximation}.
This weaker estimate is aesthetically 
less satisfactory than
Theorem \ref{theorem: new waist inequality}.
Nevertheless, it is still sufficient to prove
Theorem \ref{theorem: main theorem}.
We hope that this alternative approach will be useful
to readers who are primarily interested in
applications to mean dimension.

The results of this subsection
will not be used elsewhere in the paper,
and hence we keep our presentation brief.

\begin{theorem}[A weaker waist inequality]
\label{theorem: weaker waist inequality}
There is a universal constant $C>1$ 
with the following property.
Let $P,Q$ be nonempty finite simplicial complexes, 
and let
$f\colon P\to Q$ and $\Phi\colon P\to\mathbb R^n$
be continuous maps.
Suppose that
\[
\dim Q +3\leq d<n,\qquad 0<s<d- \dim Q -2,
\]
where $d$ is an integer.
Suppose that $\mathcal{R}:=\mathcal R_d(\Phi)>0$, and set
\[
\Delta=\diam(\Phi(P), \rho_n),\qquad
\kappa=\frac{d-\dim Q -2-s}{n-\dim Q-1}.
\]
Then, for every $r>0$,
\begin{equation}
\label{eq: weaker waist inequality}
\sup_{y\in Q}
\#\bigl(\Phi(f^{-1}(y)),\rho_n,r\bigr)
\geq
\frac{1}{\dim Q +1}
\left(\frac{\kappa \mathcal{R}}{C(\mathcal{R}
+\Delta)}\right)^n
\left(\frac{\mathcal{R}}{r}\right)^s.
\end{equation}
\end{theorem}

Compared with \eqref{eq: new waist inequality}
in Theorem \ref{theorem: new waist inequality},
the estimate \eqref{eq: weaker waist inequality}
has less satisfactory constants in two respects:
the ambient dimension $n$ appears explicitly
in the exponent of the prefactor,
and the estimate depends on the diameter
$\Delta$ of the image $\Phi(P)$.
(Here $\Delta$ is the diameter of $\Phi(P)$
with respect to the norm $\rho_n$.)
On the other hand, it has the advantage that
no restriction is imposed on $\dim P$.

\begin{proof}
We first note that the following proof
does not use the simplicial structure of $P$.
In fact, 
it applies to any nonempty compact metric space $P$,
with $\mathcal R_d(\Phi)$ defined by the same formula
as in \eqref{eq: d-radius}.

Put $w=n-d-1$ and $\lambda=w+s$.
Then
\[
0<\lambda<n-\dim Q-3,\qquad
\kappa=1-\frac{\lambda+2}{n-\dim Q-1}\in(0,1).
\]
By translating $\Phi$, we may assume that
$\rho_n(\Phi(x))\leq\Delta$ for all $x\in P$.
Set
\[
K=K(n,d,\mathcal{R}),\qquad
L=L(n,d,\mathcal{R}).
\]
Note that $w$ is the dimension of $L$.

\smallskip
\noindent
\emph{Step 1: Reduction to a small translation.}
There is a continuous
retraction
\[
\mathfrak r\colon\mathbb R^n\setminus L\longrightarrow K,
\qquad
\rho_n\bigl(\mathfrak r(z)-z\bigr)\leq\mathcal{R}/2.
\]
See Lemma \ref{lemma: retraction}
for a proof.
Figure \ref{fig:radial-retraction} 
in \S \ref{subsection: motivation}
shows the case of 
$n=2$ and $d=1$.

Suppose that a continuous map $v\colon Q\to\mathbb R^n$
satisfies
\begin{equation}
\label{eq: weaker waist translation}
\rho_n(v(y))\leq\mathcal{R}/4\quad(y\in Q),\qquad
\Phi(x)-v(f(x))\notin L\quad(x\in P).
\end{equation}
Then
$\varphi(x):=\mathfrak r\bigl(\Phi(x)-v(f(x))\bigr)$
defines a continuous map $P\to K$ satisfying
\[
\sup_{x\in P}
\rho_n\bigl(\varphi(x)-\Phi(x)\bigr)
\leq 3\mathcal{R}/4.
\]
Since $P$ is compact and $K$ is locally finite,
$\varphi(P)$ lies in a finite subcomplex of $K$
of dimension at most $d$.
This contradicts
$\mathcal R_d(\Phi)=\mathcal{R}$.
We will show that a sufficiently small fiber covering number
would produce a map satisfying
\eqref{eq: weaker waist translation}.

\smallskip
\noindent
\emph{Step 2: A finite net of the dual skeleton.}
If $r\geq\mathcal{R}$, the right-hand side of
\eqref{eq: weaker waist inequality} is at most $1$,
whereas its left-hand side is at least $1$
because $P$ is nonempty.
So the inequality is trivial in this case.
Hence we fix $0<r<\mathcal{R}$.

The set
\begin{equation}
  \label{eq: portion of L in weaker waist}
L\cap \{z\in \mathbb{R}^n\mid 
\rho_n(z)\leq\Delta+\mathcal{R}/4\}
\end{equation}
has an $r$-net $u_1,\ldots,u_B\in L$ with
\begin{equation}
\label{eq: weaker waist skeleton net}
B\leq
\left[
C_0\left(1+\frac{\Delta}{\mathcal{R}}\right)
\right]^n
\left(\frac{\mathcal{R}}{r}\right)^w,
\qquad
C_0=32 \mathrm{e}.
\end{equation}
Here an $r$-net means that every point of
\eqref{eq: portion of L in weaker waist}
is within $\rho_n$-distance at most $r$ of some $u_j$.

To verify the bound, note that the dual $n$-cubes
have centers $\mathcal{R}m$, $m\in\mathbb Z^n$,
and $\rho_n$-radius $\mathcal{R}/2$.
Every such cube meeting the ball 
$\{z\in \mathbb{R}^n\mid 
\rho_n(z)\leq\Delta+\mathcal{R}/4\}$
satisfies
\[
\norm{m}\leq\sqrt n\,T,
\qquad
T:=\Delta/\mathcal{R}+3/4.
\]
We estimate the number of such cubes
by a simple counting argument.
By the Cauchy--Schwarz inequality,
\[
\sum_{i=1}^n |m_i|
\leq \sqrt n\,\norm{m}\leq nT.
\]
Set $H=\lfloor nT\rfloor$.
The number of nonnegative integer vectors
$(a_1,\ldots,a_n)$ satisfying
$a_1+\cdots+a_n\leq H$
is $\binom{n+H}{n}$.
Since each coordinate has at most two choices
of sign, we obtain
\[
\#\{m\in\mathbb Z^n \mid \norm{m}\leq\sqrt n\,T\}
\leq 2^n\binom{n+H}{n}
\leq 2^n
\left(\frac{\mathrm e(n+H)}{n}\right)^n
\leq \bigl[2\mathrm e(1+T)\bigr]^n.
\]
Here we used the standard estimate
$\binom{N}{n}\leq(\mathrm e N/n)^n$.

Each cube has $3^n$ faces.
Every face of dimension at most $w$
has an $r$-net obtained by subdividing
each free coordinate interval into
$\lceil\mathcal{R}/r\rceil$ equal subintervals
and taking the resulting grid vertices.
This is an $r$-net with respect to $\rho_n$,
since $\rho_n\leq\norm{\cdot}_\infty$,
and it has at most
$(\lceil\mathcal{R}/r\rceil+1)^w
\leq(3\mathcal{R}/r)^w$
points, where we used $r<\mathcal{R}$.

Since $w\leq n$ and
\[
1+T
=\frac74+\frac{\Delta}{\mathcal{R}}
\leq\frac74
\left(1+\frac{\Delta}{\mathcal{R}}\right),
\]
we obtain
\[
B\leq
\left[
C_0\left(1+\frac{\Delta}{\mathcal{R}}\right)
\right]^n
\left(\frac{\mathcal{R}}r\right)^w,
\]
where we may take $C_0=32 \mathrm{e}$.

\smallskip
\noindent
\emph{Step 3: Quantitative avoidance on $Q$.}
Put
\[ M=\sup_{y\in Q}
\#\bigl(\Phi(f^{-1}(y)),\rho_n,r\bigr).\]
The construction in Steps~1 and~2 of
\S\ref{subsection: proof of new waist inequality}
gives $(\dim Q+1)M$ continuous maps
$p_i\colon Q\to\mathbb R^n$ such that
\begin{equation}
\label{eq: weaker waist moving centers}
\min_i
\rho_n\bigl(\Phi(x)-p_i(f(x))\bigr)<r
\qquad(x\in P).
\end{equation}

Apply Proposition~\ref{proposition: quantitative avoidance}
on $Q$, with reference map $0$, parameter
$a=n-\dim Q-1$, exponent $\lambda$, and the
$(\dim Q+1)MB$ obstacles
\[
Q\ni y\longmapsto p_i(y)-u_j \in \mathbb{R}^n.
\]
The dimensional and exponent hypotheses hold because
\[
3\leq a\leq n-1,\qquad
\dim Q=n-a-1,\qquad
0<\lambda<a-2.
\]
If
\begin{equation}
\label{eq: weaker waist avoidance condition}
\left\{(\dim Q+1)MB\right\}^{1/\lambda}\,2r
<
\kappa\mathcal{R}/4,
\end{equation}
Proposition \ref{proposition: quantitative avoidance} 
gives a continuous map
$v\colon Q\to\mathbb R^n$ satisfying
\[
\rho_n(v(y))\leq\mathcal{R}/4,\qquad
\rho_n\bigl(v(y)-(p_i(y)-u_j)\bigr)>2r
\]
for all $y,i,j$.

This map satisfies \eqref{eq: weaker waist translation}.
Indeed, if
$u:=\Phi(x)-v(f(x))\in L$,
then $\rho_n(u)\leq\Delta+\mathcal{R}/4$.
Choose $j$ with $\rho_n(u-u_j)\leq r$ and $i$
as in \eqref{eq: weaker waist moving centers}.
Then
\[
\rho_n\bigl(v(f(x))-(p_i(f(x))-u_j)\bigr)
=\rho_n\bigl(\Phi(x)-p_i(f(x))+u_j-u\bigr)
<2r,
\]
a contradiction.

\smallskip
\noindent
\emph{Step 4: The covering estimate.}
By Step~1,
\eqref{eq: weaker waist avoidance condition}
cannot hold.
Together with \eqref{eq: weaker waist skeleton net},
this gives
\[
M
\geq
\frac{1}{(\dim Q+1)B}
\left(\frac{\kappa\mathcal{R}}{8r}\right)^\lambda
\geq
\frac{1}{(\dim Q+1)C_0^n}
\left(\frac{\kappa}{8}\right)^\lambda
\left(\frac{\mathcal{R}}{\mathcal{R}+\Delta}\right)^n
\left(\frac{\mathcal{R}}{r}\right)^s,
\]
where we used $\lambda-w=s$.
Since $0<\lambda<n$ and $0<\kappa/8<1$, we have
\[
(\kappa/8)^\lambda\geq(\kappa/8)^n.
\]
Thus \eqref{eq: weaker waist inequality} follows with
$C=8C_0 = 256\mathrm{e}$.

We emphasize that this proof does not use
Proposition \ref{proposition: low dimensional approximation}.
The key new ingredient is the estimate
\eqref{eq: weaker waist skeleton net}
for the size of an $r$-net covering
a bounded portion of the dual skeleton.
This allows us to apply
Proposition~\ref{proposition: quantitative avoidance}
directly on $Q$.
The resulting proof is more elementary,
at the cost of an additional factor
depending on the ambient dimension $n$
and the diameter of $\Phi(P)$
in the waist inequality.
\end{proof}

\medskip
\noindent
\emph{Application to dynamical systems.}
We now explain how
Theorem~\ref{theorem: weaker waist inequality}
gives another proof of
Theorem~\ref{theorem: main theorem}.
Let $\pi\colon(X,T)\to(Y,S)$ be an equivariant
continuous map, and let $\mathbf d$ be a compatible
metric on $X$.
It suffices to consider the case
\[
\mdim(Y,S)<\mdim(X,T),\qquad \mdim(Y,S)<\infty.
\]
Choose finite positive numbers satisfying
\[
\mdim(Y,S)<b<a<\alpha<\beta<\mdim(X,T).
\]

As in Step~1 of the proof of
Theorem~\ref{theorem: main theorem}
in \S\ref{subsection: the case of finite mean dimension},
choose $0<\tau<1$ and a $1$-Lipschitz continuous map
$F\colon X\to[0,1]^\ell$ such that
\[
\widim_\tau(X,\mathbf d_N)>\beta N
\]
for all sufficiently large $N$, and
\[
\norm{F(x)-F(x')}_\infty<\tau/2
\quad\Longrightarrow\quad
\mathbf d(x,x')<\tau.
\]
By adding zero coordinates, we may assume that $\ell>\beta$.

Fix $0<\delta<\tau/8$.
Apply Proposition~\ref{proposition: polyhedral modeling}
with $\eta=1$, and let $L\in\mathbb N$ be the resulting
integer.
For $N\in L\mathbb N$, take the maps
\[
f_N\colon X\to P_N,\qquad
g_N\colon P_N\to Q_N,\qquad
\Phi_N\colon P_N\to[0,1]^{\ell N}
\]
provided by Proposition~\ref{proposition: polyhedral modeling}, 
and put
\[
d=\lfloor\alpha N\rfloor,\qquad
k=\lfloor(\beta-\alpha)N\rfloor,\qquad
s=(a-b)N.
\]
The same argument as for
\eqref{eq: finite case radius lower bound} gives
\[
\mathcal R^\infty_{d+k}(\Phi_N)\geq\tau/8
\]
for all sufficiently large $N$.
Consequently, Lemma~\ref{lemma: rho_n and ell^infty}
implies
\[
\mathcal R_N:=\mathcal R_d(\Phi_N)
\geq
\sqrt{\frac{k+1}{\ell N}}\,
\mathcal R^\infty_{d+k}(\Phi_N)
\geq
\mathcal R_0:=
\frac{\tau}{8}\sqrt{\frac{\beta-\alpha}{\ell}}>0.
\]
Also, since $\Phi_N(P_N)\subset[0,1]^{\ell N}$,
\[
\Delta_N:=\diam (\Phi_N(P_N), \rho_{\ell N}) \leq 1.
\]

Since $\dim Q_N<bN$, we have
\[
d-\dim Q_N-2-s>(\alpha-a)N-3.
\]
Thus, for all sufficiently large $N$,
the dimensional and exponent hypotheses of
Theorem~\ref{theorem: weaker waist inequality}
hold with $n=\ell N$, and
\[
\kappa_N:=
\frac{d-\dim Q_N-2-s}{\ell N-\dim Q_N-1}
\geq
\frac{\alpha-a}{2\ell}
=:\kappa_0>0.
\]
Let $C>1$ be the universal constant in 
Theorem~\ref{theorem: weaker waist inequality},
and set
\[
B_0:=
\ell\log\left(
\frac{C(\mathcal R_0+1)}{\kappa_0\mathcal R_0}
\right)>0.
\]
Using
\[
\rho_{\ell N}\leq\norm{\cdot}_\infty,
\qquad
\frac{\mathcal R_N}{\mathcal R_N+\Delta_N}
\geq
\frac{\mathcal R_0}{\mathcal R_0+1},
\]
we obtain
\[
\sup_{y\in Q_N}
\#\bigl(\Phi_N(g_N^{-1}(y)),
       \norm{\cdot}_\infty,\delta\bigr)
\geq
\frac{e^{-B_0N}}{bN+1}
\left(\frac{\mathcal R_0}{\delta}\right)^{(a-b)N}.
\]
Here $\mathcal R_0$ and $B_0$ are independent of both
$N$ and $\delta$.

On the other hand,
Proposition~\ref{proposition: polyhedral modeling} gives
\[
\sup_{y\in Q_N}
\#\bigl(\Phi_N(g_N^{-1}(y)),
       \norm{\cdot}_\infty,\delta\bigr)
\leq
\exp\bigl((S_\delta(\pi)+1)N\bigr).
\]
Taking logarithms and dividing by $N$, we obtain
\[
S_\delta(\pi)+1
\geq
(a-b)\log(\mathcal R_0/\delta)
-B_0-\frac{\log(bN+1)}{N}.
\]
For each fixed $\delta$, let $N\to\infty$
through multiples of $L$.
Then
\[
S_\delta(\pi)+1
\geq
(a-b)\log(\mathcal R_0/\delta)-B_0.
\]
Dividing by $\log(1/\delta)$ and taking the lower limit
as $\delta\to0$, we conclude that
\[
\lmdimM(\pi,T,\mathbf d)\geq a-b.
\]
Since $a<\mdim(X,T)$ and $b>\mdim(Y,S)$
can be chosen arbitrarily, this proves
Theorem~\ref{theorem: main theorem}.

This argument also covers $\mdim(X,T)=\infty$.
Indeed, neither polyhedral modeling nor the radius
lower bound requires finite mean dimension,
and no upper bound on $\dim P_N$ is needed here.

\section{Proof of 
Proposition \ref{proposition: quantitative avoidance}}
\label{section: proof of proposition avoidance}

The purpose of this section is to prove 
Proposition \ref{proposition: quantitative avoidance} 
(Quantitative avoidance). 
We write points of $\mathbb{R}^n$ 
as $z=(z_1,\ldots,z_n)$, 
denote the standard Euclidean norm by $\norm{z}$, 
and recall that $\rho_n(z)=\norm{z}/\sqrt{n}$.

We recall the notion of \textit{connectivity} 
of topological spaces. 
A topological space $X$ is called 
\textbf{$(-1)$-connected} if it is nonempty, 
and \textbf{$0$-connected} if it is path-connected. 
For an integer $k\geq1$, 
we say that $X$ is \textbf{$k$-connected} 
if it is path-connected and every continuous map from the 
$i$-dimensional sphere $S^i$ to $X$ 
is null-homotopic for $1\leq i\leq k$.

We first explain a toy model illustrating 
the main task of the proof of 
Proposition \ref{proposition: quantitative avoidance}.
Let $K$ be a locally finite simplicial complex, 
and let $f,p\colon K\to\mathbb{R}^n$ be continuous maps.
Given $r,R>0$, we would like to find a continuous map 
$g\colon K\to\mathbb{R}^n$ satisfying
\begin{equation}
\label{eq: toy model g(x)}
\rho_n(f(x)-g(x)) \leq R,
\qquad
\rho_n(g(x)-p(x))>r
\end{equation}
for all $x\in K$.

For each $x\in K$, define the 
\lq\lq{}good region\rq\rq{}
 $G_x$ by
\[
G_x=\left\{z\in\mathbb{R}^n\;\middle|\;
\rho_n(f(x)-z)\leq R,\quad
\rho_n(z-p(x))>r
\right\}.
\]
Our problem is to choose a point $g(x)\in G_x$ 
continuously in $x$.

When $\dim K=0$, 
this is possible whenever $G_x$ is nonempty
for every $x\in K$, 
because $K$ is discrete. 
However, when $\dim K\geq1$, 
the nonemptiness of every $G_x$ is not 
sufficient in general: 
the choices $g(x)\in G_x$ 
must depend continuously on $x$.

For simplicity, suppose that 
$G:= \bigcup_{x\in K}\{x\}\times G_x$ forms a locally 
trivial fiber bundle\footnote{
Namely, every $x\in K$ has 
an open neighborhood $U_x$ such that
we can find a homeomorphism 
$\varphi_y\colon G_y\to G_x$ $(y\in U_x)$ 
that depends continuously 
on $y\in U_x$.
}  over $K$.
(This assumption does not hold in general, 
and the actual proof
requires a modification to handle this issue.)
Now, standard obstruction theory tells us 
that if every fiber $G_x$ is 
$(\dim K-1)$-connected, 
then there exists a continuous map 
$g\colon K\to\mathbb{R}^n$ 
such that $g(x)\in G_x$ for every $x\in K$.

Thus the main task is to estimate 
the connectivity of the regions $G_x$. 
In the actual proof of 
Proposition 
\ref{proposition: quantitative avoidance}, 
we will address the failure of local triviality 
and use \textit{Morse theory} 
to establish the required connectivity estimates.

\subsection{Morse theory on the sphere}
\label{subsection: Morse theory on the sphere}

In this subsection, 
we recall some basic differential-geometric formulas 
on spheres and establish a connectivity result.

We denote by $\langle\cdot,\cdot\rangle$ 
the standard inner product on $\mathbb{R}^n$. 
For a smooth function $H\colon\mathbb{R}^n\to\mathbb{R}$, 
we define its gradient $\nabla H(z)$ 
and Hessian matrix $\nabla^2 H(z)$ by
\[
\nabla H (z)=
\left(
\frac{\partial H}{\partial z_1}(z),\ldots,
\frac{\partial H}{\partial z_n}(z)
\right)^T,
\qquad
\nabla^2 H (z)=
\left(
\frac{\partial^2 H}{\partial z_i\partial z_j}(z)
\right)_{1\leq i,j\leq n}.
\]

Let $R>0$. 
Consider the sphere
\[
\Sigma :=\left\{z\in\mathbb{R}^n \mid
\norm{z}=R\right\}.
\]
For each $z\in\Sigma$, 
its tangent space is given by
\[
T_z\Sigma=
\left\{\xi\in\mathbb{R}^n\mid
\langle\xi,z\rangle=0\right\}.
\]
We equip $\Sigma$ with the Riemannian metric 
induced by the standard Euclidean inner product 
$\langle\cdot,\cdot\rangle$ (not by $\rho_n$).

Let $h\colon \Sigma \to \mathbb{R}$ be a
smooth function. 
We define the \textbf{Hessian of $h$} at $z\in\Sigma$ 
as the quadratic form
\[
\operatorname{Hess}_{z}h\colon 
T_z\Sigma\to\mathbb{R}
\]
given by
\[
\operatorname{Hess}_{z}h(\xi)
=
\left.\frac{d^2}{dt^2}h(\gamma(t))\right|_{t=0},
\]
where $\gamma$ is the (constant-speed)
geodesic on $\Sigma$ satisfying
$\gamma(0)=z$ and $\gamma'(0)=\xi$.

\begin{lemma}
 \label{lemma: hessian}
  Let $H\colon \mathbb{R}^n\to \mathbb{R}$ be a smooth 
  function, and let $h = H|_\Sigma \colon \Sigma\to \mathbb{R}$
  be the restriction of $H$ to $\Sigma$. Then 
  for $\xi \in T_z \Sigma$
  \[ \operatorname{Hess}_z h (\xi)
  = \langle \nabla^2 H(z)\xi, \xi\rangle -
    \frac{\langle \nabla H(z), z\rangle}{R^2} 
    \norm{\xi}^2.\]
\end{lemma}

\begin{proof}
Take a nonzero $\xi \in T_z \Sigma$.
The geodesic $\gamma$ satisfying $\gamma(0)= z$ and 
$\gamma^\prime(0)=\xi$ is given by 
\[ \gamma(t) =  \cos \left(\frac{\norm{\xi}t}{R}\right) z
   + \frac{R}{\norm{\xi}}
   \sin\left(\frac{\norm{\xi}t}{R}\right) \xi. \]
Differentiating this formula, we obtain
\[
\gamma'(0)=\xi,\qquad
\gamma''(0)=-\frac{\norm{\xi}^2}{R^2}z.
\]
By the chain rule,
\[
\frac{d}{dt}H(\gamma(t))
=
\left\langle
\nabla H(\gamma(t)),\gamma'(t)
\right\rangle.
\]
Differentiating once more, we have
\[
\begin{aligned}
\frac{d^2}{dt^2}H(\gamma(t))
={}&
\left\langle
\nabla^2H(\gamma(t))\gamma'(t),
\gamma'(t)
\right\rangle\\
&+
\left\langle
\nabla H(\gamma(t)),\gamma''(t)
\right\rangle.
\end{aligned}
\]
Evaluating at $t=0$, we obtain
\[
\begin{aligned}
\operatorname{Hess}_z h(\xi)
&=
\left.\frac{d^2}{dt^2}H(\gamma(t))
\right|_{t=0}\\
&=
\left\langle\nabla^2H(z)\xi,\xi\right\rangle
-
\frac{\norm{\xi}^2}{R^2}
\left\langle\nabla H(z),z\right\rangle.
\end{aligned}
\]
The formula is also valid for $\xi=0$.
This completes the proof.
\end{proof}

The following lemma is the key result of this subsection. 
Its proof is based on Morse theory.

\begin{lemma}
\label{lemma: Morse theory}
Let $2\leq\ell\leq n-1$ be integers, and let $c\in\mathbb{R}$.
Let $h\colon\Sigma\to\mathbb{R}$ be a smooth function.
Suppose that, for every $z\in\Sigma$ with $h(z)>c$,
there exists an $\ell$-dimensional linear subspace
$E_z\subset T_z\Sigma$ such that
\[
\operatorname{Hess}_z h(\xi)<0
\qquad(\forall \xi\in E_z\setminus\{0\}).
\]
Then, for every $t>c$, the sublevel set
\[
\{z\in\Sigma\mid h(z)<t\}
\]
is nonempty and $(\ell-2)$-connected.
\end{lemma}

\begin{proof}
First notice that $\Sigma$ is homeomorphic to $S^{n-1}$
and hence $(n-2)$-connected.

Fix $t>c$, and set
\[
U_t=\{z\in\Sigma \mid h(z)<t\}.
\]
Since $\Sigma$ is compact, $h$ attains its minimum
at some point $z_{\min}\in\Sigma$.
The Hessian of $h$ at $z_{\min}$ is positive semidefinite,
so the assumption implies $h(z_{\min})\leq c$.
Hence $U_t$ is nonempty.

We claim that every nonempty compact subset $C\subset U_t$
is contained in an $(\ell-2)$-connected subset $A\subset U_t$.
Choose a real number $a$ such that
\[
\max\left\{c,\max_{z\in C}h(z)\right\}<a<t.
\]
On the compact set
\[
B=\{z\in\Sigma \mid h(z)\geq a\},
\]
the Hessian of $h$ has at least $\ell$ negative eigenvalues
at every point.
This property is preserved by sufficiently small
$\mathcal{C}^2$-perturbations of $h$.
(Here the $\mathcal C^2$ topology is the topology
of uniform convergence of functions and their
derivatives up to order two.)
Indeed, the $\ell$-th smallest eigenvalue of the Hessian
is a continuous function of $z$ and is strictly negative
on $B$.

Choose $0<\varepsilon<(t-a)/2$.
By the Morse approximation theorem 
\cite[Corollary 6.8]{Milnor}, there exists a smooth
Morse function $\widetilde{h}\colon\Sigma\to\mathbb{R}$
sufficiently close to $h$ in the $\mathcal{C}^2$ 
topology that
\[
\sup_{z\in\Sigma}|\widetilde{h}(z)-h(z)|<\varepsilon
\]
and the Hessian of $\widetilde{h}$ has at least $\ell$
negative eigenvalues at every point of $B$.
Choose a regular value $b$ of $\widetilde{h}$ satisfying
\[
a+\varepsilon<b<t-\varepsilon,
\]
and set
\[
A=\{z\in\Sigma \mid \widetilde{h}(z)\leq b\}.
\]
Then
\[
C\subset A\subset U_t.
\]
Indeed, for $z\in C$ we have
$\widetilde{h}(z)<a+\varepsilon<b$,
whereas for $z\in A$ we have
$h(z)<b+\varepsilon<t$.

Every critical point $z$ of $\widetilde{h}$ with
$\widetilde{h}(z)>b$ satisfies
\[
h(z)>\widetilde{h}(z)-\varepsilon
>b-\varepsilon>a.
\]
Consequently, every such critical point 
has Morse index\footnote{
The Morse index of a nondegenerate critical point 
is the number of negative eigenvalues 
of the Hessian at that point, counted with multiplicity.  
}
at least $\ell$.
By the cell-attachment theorem of Morse theory
\cite[Theorems 3.1 and 3.2, Remark 3.3]{Milnor},
$\Sigma$ is obtained from $A$, up to homotopy relative
to $A$, by attaching cells of dimensions at least $\ell$.
Attaching cells of dimensions at least $\ell$ 
does not change $(\ell-2)$-connectivity. 
Since $\Sigma$ is $(n-2)$-connected 
and $\ell\leq n-1$, 
it follows that $A$ is $(\ell-2)$-connected. 
This proves the claim.

Taking $C$ to consist of any two points of $U_t$,
we see that $U_t$ is path-connected.
Moreover, for $1\leq j\leq\ell-2$, the image of any
continuous map $\alpha\colon S^j\to U_t$ is compact.
By the claim, this image is contained in an
$(\ell-2)$-connected subset $A\subset U_t$.
Thus $\alpha$ is null-homotopic in $A$, and hence in $U_t$.
Therefore $U_t$ is $(\ell-2)$-connected.
\end{proof}

\subsection{Proof of Proposition 
\ref{proposition: quantitative avoidance}}
\label{subsection: proof of quantitative avoidance}

We now prove Proposition \ref{proposition: quantitative avoidance}.
For the reader's convenience, we first restate it.

\begin{restatedavoidance}[Quantitative avoidance]
Let $m,n$ be positive integers, and let $a$ be an integer
and $s$ a real number satisfying
\[
3\leq a\leq n-1,\qquad 0<s<a-2.
\]
Set
\[
\kappa=1-\frac{s+2}{a} >0.
\]
Let $K$ be a locally finite simplicial complex satisfying
\[
\dim K\leq n-a-1.
\]
Let $f,p_1,\ldots,p_m\colon K\to\mathbb{R}^n$
be continuous maps.
Suppose that $r,R>0$ satisfy
\[
m^{1/s}r<\kappa R.
\]
Then there exists a continuous map
$g\colon K\to\mathbb{R}^n$ such that
\[
\rho_n(g(x)-f(x))\leq R,\qquad
\rho_n(g(x)-p_i(x))>r
\]
for every $x\in K$ and every $1\leq i\leq m$.
\end{restatedavoidance}

Although the proposition only requires
\[
\rho_n(g(x)-f(x))\leq R,
\]
we will actually construct a map $g$ satisfying
\[
\rho_n(g(x)-f(x))=R
\]
for every $x\in K$.
The reason is technical: 
this allows us to work on a sphere rather 
than a closed ball. 
Since the sphere is a smooth manifold without boundary, 
this simplifies our application of Morse theory.

The proof has two parts.
First, we construct smooth auxiliary functions 
on a fixed sphere and use Lemma \ref{lemma: Morse theory} 
to estimate the connectivity of their sublevel sets.
We then construct a continuous selection 
by induction over the skeleta of a sufficiently 
fine subdivision of $K$.
By allowing a slight increase in the sublevel thresholds 
at each step, we can 
carry out this construction 
without assuming local triviality.

\begin{proof}
We may assume that $K$ is nonempty, and put $k=\dim K$.

\smallskip
\noindent
\emph{Step 1: Normalization and auxiliary functions.}
For $1\leq i\leq m$, define
\[
q_i(x)=\frac{p_i(x)-f(x)}{\sqrt n},
\qquad
\Sigma =\{z\in\mathbb{R}^n \mid \norm{z}=R\}.
\]
It suffices to construct a continuous map
$v\colon K\to\Sigma$ such that
\begin{equation}
\label{eq: avoidance normalized separation}
\norm{v(x)-q_i(x)}>r
\qquad(x\in K,\ 1\leq i\leq m).
\end{equation}
Indeed, we can then set $g(x)=f(x)+\sqrt n\,v(x)$.
Notice that $R$ is the Euclidean radius of $\Sigma$;
the normalization above accounts for the factor $\sqrt n$
in the definition of $\rho_n$.

Set
\[
c_0=\frac{m}{(\kappa R)^s}.
\]
Our assumption gives $c_0<r^{-s}$.
Choose $\varepsilon>0$ so small that
\begin{equation}
\label{eq: avoidance threshold gap}
c_0<t_*:=(r^2+\varepsilon^2)^{-s/2}.
\end{equation}
For $x\in K$, define the smooth function
\[
H_x(z)=\sum_{i=1}^m
\bigl(\norm{z-q_i(x)}^2+\varepsilon^2\bigr)^{-s/2}
\qquad(z\in\mathbb{R}^n),
\]
and let $h_x=H_x|_{\Sigma}$.
The function $H_x(z)$ has the following crucial 
property:
\begin{equation}
\label{eq: avoidance sublevel implies separation}
H_x(z)<t_*
\quad\Longrightarrow\quad
\norm{z-q_i(x)}>r\quad(1\leq i\leq m).
\end{equation}
Indeed, if $\norm{z-q_i(x)}\leq r$ for some $i$,
then the corresponding summand is at least $t_*$.

\smallskip
\noindent
\emph{Step 2: A Hessian estimate on the sphere.}
Fix $x\in K$ and $z\in\Sigma$.
For brevity, write
\[
u_i=z-q_i(x),\qquad
D_i=\bigl(\norm{u_i}^2+\varepsilon^2\bigr)^{1/2},
\qquad
A=\sum_{i=1}^m D_i^{-s-2},\qquad
B=\sum_{i=1}^m D_i^{-s-1}.
\]
We have $H_x(z) = \sum_{i=1}^m D_i^{-s}$.
The Cauchy--Schwarz inequality 
and the H\"older inequality with
conjugate exponents
$\frac{s+2}{2}$ and $\frac{s+2}{s}$ 
give
\[
B^2\leq H_x(z)A,
\qquad
H_x(z)\leq m^{2/(s+2)}A^{s/(s+2)}, \quad 
\text{respectively}.
\]
Consequently,
\begin{equation}
\label{eq: avoidance moment estimate}
\frac{B}{A}
\leq\sqrt{\frac{H_x(z)}{A}}
\leq\left(\frac{m}{H_x(z)}\right)^{1/s}.
\end{equation}
In particular, if $H_x(z)>c_0 = m/(\kappa R)^s$, 
then
\begin{equation}
\label{eq: avoidance moment strict bound}
\frac{B}{RA}<\kappa.
\end{equation}

Direct differentiation gives
\[
\nabla H_x(z)=-s\sum_{i=1}^m D_i^{-s-2}u_i,
\]
and
\[
\nabla^2H_x(z)
=s(s+2)\sum_{i=1}^m D_i^{-s-4}u_i u_i^T-sA I_n,
\]
where $I_n$ is the identity matrix.
Let $\xi_1,\ldots,\xi_a$ be any orthonormal vectors
in $T_z\Sigma$.
(Notice that $a$ is the positive integer given in 
the statement of the proposition.)
By Lemma \ref{lemma: hessian}, 
\begin{align*}
\sum_{j=1}^a\operatorname{Hess}_z h_x(\xi_j)
&=s(s+2)\sum_{i=1}^m D_i^{-s-4}
    \sum_{j=1}^a\langle u_i,\xi_j\rangle^2
    -saA
    +\frac{sa}{R^2}\sum_{i=1}^m D_i^{-s-2}
      \langle u_i,z\rangle\\
&\leq s(s+2)A-saA+\frac{sa}{R}B, 
\qquad  \text{by \eqref{eq: u_i xi_i D_i} below}\\
&=saA\left(\frac{B}{RA}-\kappa\right), \qquad 
\text{by $\kappa = 1-\frac{s+2}{a}$}.
\end{align*}
Here we used
\begin{equation}
  \label{eq: u_i xi_i D_i}
\sum_{j=1}^a\langle u_i,\xi_j\rangle^2
\leq\norm{u_i}^2\leq D_i^2,
\qquad
\langle u_i,z\rangle\leq R\norm{u_i}\leq RD_i.
\end{equation}
It follows from 
\eqref{eq: avoidance moment strict bound} that
\begin{equation}
\label{eq: avoidance negative trace}
\sum_{j=1}^a\operatorname{Hess}_z h_x(\xi_j)<0
\qquad\text{whenever }H_x(z)>c_0.
\end{equation}

Since $\dim T_z\Sigma =n-1$, this implies that
$\operatorname{Hess}_z h_x$ has at least $n-a$ strictly
negative eigenvalues whenever $H_x(z)>c_0$.
Otherwise, it would have at least $a$ 
nonnegative eigenvalues,
and corresponding orthonormal eigenvectors 
would contradict
\eqref{eq: avoidance negative trace}.

\smallskip
\noindent
\emph{Step 3: Connectivity of the sublevel sets.}
For $x\in K$ and $t>c_0$, set
\[
U_x(t)=\{z\in\Sigma \mid H_x(z)<t\}.
\]
We claim that $U_x(t)$ is nonempty and
$(n-a-2)$-connected.
If $n-a\geq2$, this follows from
Lemma \ref{lemma: Morse theory}, with $\ell=n-a$,
by the Hessian estimate in Step~2.
If $n-a=1$, only nonemptiness is required.
In this case, $h_x$ attains a minimum 
on the compact sphere
$\Sigma$, and the Hessian at a minimum is positive
semidefinite.
Step~2 therefore implies $\min h_x\leq c_0<t$,
which proves nonemptiness.
Thus the claim holds in all cases.
Since $k\leq n-a-1$, each $U_x(t)$ is, in particular,
$(k-1)$-connected.

\smallskip
\noindent
\emph{Step 4: Constructing a continuous selection.}
Choose $\eta>0$ such that
\begin{equation}
\label{eq: avoidance selection margin}
c_0+(2k+1)\eta<t_*.
\end{equation}
By continuity of $H_x$ and compactness of $\Sigma$,
every point of $K$ has an open neighborhood $V$ such that
\[
\sup_{z\in\Sigma}|H_x(z)-H_{x'}(z)|<\eta
\qquad(x,x'\in V).
\]
Choose a locally finite subdivision of $K$ subordinate
to this open cover, and keep the notation $K$ for it.
Then every simplex $\sigma$ satisfies
\begin{equation}
\label{eq: avoidance variation on simplices}
|H_x(z)-H_{x'}(z)|<\eta
\qquad(x,x'\in\sigma,\ z\in\Sigma).
\end{equation}
For $0\leq j \leq k$, we denote the 
$j$-skeleton of $K$ by $K^{(j)}$.

We construct $v$ inductively over the skeleta $K^{(j)}$
so that
\begin{equation}
\label{eq: avoidance skeletal bound}
H_x(v(x))<c_0+(2j+1)\eta
\qquad(x\in K^{(j)}).
\end{equation}
For each vertex $x \in \mathrm{Ver}(K)$, 
choose $v(x)\in U_x(c_0+\eta)$,
which is possible by Step~3.
This proves \eqref{eq: avoidance skeletal bound} for $j=0$.

Suppose that $1\leq j\leq k$ and that $v$ has been
constructed on $K^{(j-1)}$.
Let $\sigma$ be a $j$-simplex and choose $x_\sigma\in\sigma$.
For every $x\in\partial\sigma$,
\eqref{eq: avoidance skeletal bound} and
\eqref{eq: avoidance variation on simplices} imply
\[
H_{x_\sigma}(v(x))
<H_x(v(x))+\eta
<c_0+2j\eta.
\]
Thus $v|_{\partial\sigma}$ takes values in the fixed set
$U_{x_\sigma}(c_0+2j\eta)$.
This set is $(k-1)$-connected and $j\leq k$,
so the map $v|_{\partial\sigma}$ extends continuously to
\[
v_\sigma\colon\sigma\to U_{x_\sigma}(c_0+2j\eta).
\]
For $j=1$, this uses path-connectedness; for $j\geq2$,
it uses the vanishing of the $(j-1)$-st homotopy group.
For every $x\in\sigma$, we have
\[
H_x(v_\sigma(x))
<H_{x_\sigma}(v_\sigma(x))+\eta
<c_0+(2j+1)\eta.
\]
Performing this construction on every $j$-simplex gives
an extension to $K^{(j)}$ satisfying
\eqref{eq: avoidance skeletal bound}.
The extensions agree on common faces and are continuous
in the weak topology of the simplicial complex.

After $k$ steps we obtain a continuous map
$v\colon K\to\Sigma$ satisfying
\[
H_x(v(x))<c_0+(2k+1)\eta<t_*
\qquad(x\in K).
\]
By \eqref{eq: avoidance sublevel implies separation},
this map satisfies 
\eqref{eq: avoidance normalized separation}.
Finally, set
\[
g(x)=f(x)+\sqrt n\,v(x).
\]
Then $g$ is continuous and, 
for every $x\in K$ and $1\leq i\leq m$,
\[
\rho_n(g(x)-f(x))=\norm{v(x)}=R,
\qquad
\rho_n(g(x)-p_i(x))=\norm{v(x)-q_i(x)}>r.
\]
This completes the proof.
\end{proof}

\begin{remark}
The Hessian estimates for Riesz-type potentials 
used above
are closely related to the theory of
$p$-plurisubharmonic functions studied by
Harvey and Lawson~\cite{Harvey--Lawson13}.
For Morse-theoretic aspects of inverse-power potentials
generated by point charges, see also
Gabrielov, 
Novikov and Shapiro~\cite{Gabrielov--Novikov--Shapiro07}.
\end{remark}

\section{Proof of Proposition
\ref{proposition: low dimensional approximation}}
\label{section: proof of proposition low 
dimensional approximation}

In this section, we prove Proposition
\ref{proposition: low dimensional approximation}.
The proof relies on ideas from algebraic topology
that may be unfamiliar to some readers
working in dynamical systems.
We therefore begin by explaining the motivation
behind the proof in
\S\ref{subsection: motivation}
and then review the necessary background
from algebraic topology in
\S\ref{subsection: preliminaries on algebraic topology}.
We complete the proof in
\S\ref{subsection: proof of low-dimensional approximation}.

Recall that, for $x\in\mathbb{R}^n$, 
$\norm{x}$ denotes the Euclidean norm 
and $\rho_n(x)=\norm{x}/\sqrt{n}$.
For $h>0$ and integers $1\leq d<n$, 
let $K(n,d,h)$ denote the $d$-skeleton of the cubical grid 
of mesh size $h$ in $\mathbb{R}^n$ 
with vertex set $h\mathbb{Z}^n$. 
Similarly, let $L(n,d,h)$ denote the $(n-d-1)$-skeleton of 
the dual grid with vertex set $h(\mathbb{Z}+1/2)^n$.

\subsection{Motivation}
\label{subsection: motivation}

Let us recall the setting of Proposition 
\ref{proposition: low dimensional approximation}.
Set $K=K(n,d,h)$ and $L=L(n,d,h)$.
Let $P$ be a finite simplicial complex satisfying
\[
\dim P\leq 2d-1.
\]
Suppose that we are given continuous maps
\[
\Phi\colon P\to\mathbb{R}^n,\qquad
F\colon P\times L\to\mathbb{R}^n
\]
and a constant $R>0$ such that
\[
\rho_n(F(p,x)-x)\leq R,
\qquad
F(p,x)\neq\Phi(p)
\]
for every $p\in P$ and $x \in L$.
Our goal is to construct a continuous map
$\varphi\colon P\to K$ satisfying
\[
\sup_{p\in P}
\rho_n(\varphi(p)-\Phi(p))<R+4h.
\]

\begin{figure}[htbp]
\centering
\includegraphics[width=\textwidth]
{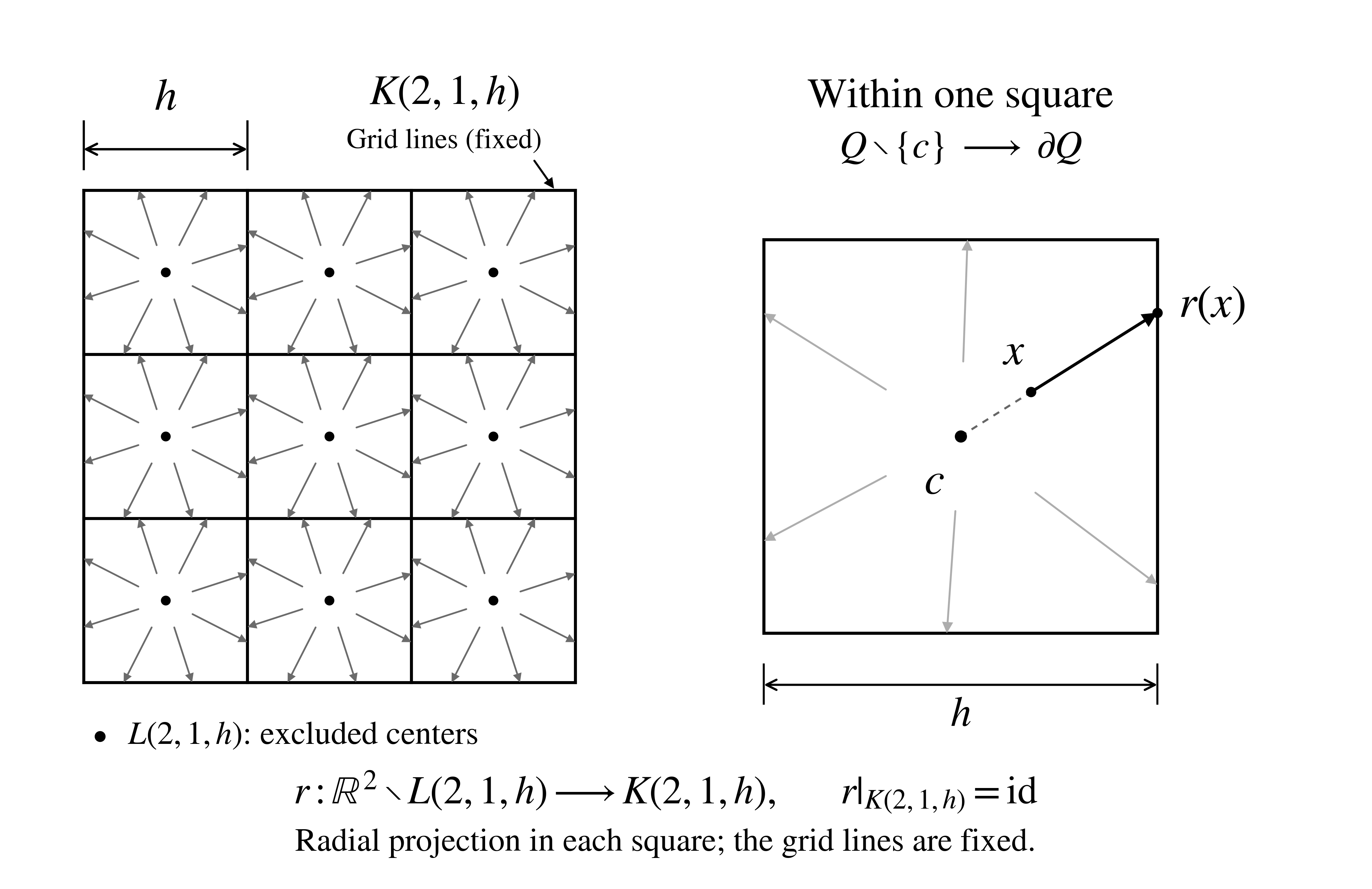}
\caption{
The radial retraction
$r\colon \mathbb{R}^2\setminus L(2,1,h)\to K(2,1,h)$.
Within each square, a point $x$ distinct from the center $c$
is mapped to the intersection $r(x)$ of the ray from $c$
through $x$ with the boundary of the square.
The grid lines are fixed, while the cell centers
are excluded from the domain.
}
\label{fig:radial-retraction}
\end{figure}

We first notice that 
there is a retraction 
$r\colon \mathbb{R}^n \setminus L \to K$
satisfying $\rho_n(r(x)-x) \leq h/2$.
Figure \ref{fig:radial-retraction} shows the case of 
$n=2$ and $d=1$.
(See Lemma \ref{lemma: retraction}
in \S
\ref{subsection: proof of low-dimensional approximation}
for the precise construction of $r$.)

By applying the Tietze extension theorem 
to the map 
\[ P\times L \ni (p, x)\mapsto 
   F(p, x)-x \in \mathbb{R}^n, \]
we can extend $F$ to a continuous map 
$F\colon P\times \mathbb{R}^n\to \mathbb{R}^n$ 
so that 
\[ \rho_n\left(F(p,x) - x\right) \leq R \]
for every $(p, x)\in P\times \mathbb{R}^n$.
(See Lemma \ref{lemma: Tietze extension}
in \S 
\ref{subsection: proof of low-dimensional approximation}
for the details.)
For each $p\in P$, we define 
$F_p\colon \mathbb{R}^n \to \mathbb{R}^n$
by $F_p(x) = F(p, x)$.

\smallskip
\noindent
\emph{Approach 1.}
Suppose that $F_p$ is a homeomorphism of $\mathbb{R}^n$
for every $p\in P$.
Then we can define $f\colon P\to \mathbb{R}^n$ by 
\[
f(p)=F_p^{-1}(\Phi(p)).
\]
This is a continuous map\footnote{
Indeed, consider the map
\[
\widehat F\colon P\times\mathbb{R}^n
\longrightarrow P\times\mathbb{R}^n,
\qquad
\widehat F(p,x)=(p,F_p(x)).
\]
This is a continuous bijection.
Moreover, it is proper because $P$ is compact and
$\rho_n(F_p(x)-x)\leq R$.
Hence $\widehat F$ is a homeomorphism.
This proves the continuity of $f$.}.

Since $F_p(x)\neq\Phi(p)$ for all $x\in L$,
we have
\[
f(p)\notin L
\qquad(p\in P).
\]
Moreover,
\[
\rho_n(f(p)-\Phi(p))
=\rho_n(f(p)-F_p(f(p)))
\leq R.
\]

Define
$\varphi:=r\circ f\colon P\longrightarrow K$.
Then $\varphi$ is continuous, and
\[
\rho_n(\varphi(p)-\Phi(p))
 \leq
\rho_n(r(f(p))-f(p))
+\rho_n(f(p)-\Phi(p))
\leq \frac{h}{2}+R.
\]
Thus we obtain the desired low-dimensional
approximation under the additional assumption
that every $F_p$ is a homeomorphism.
This is a prototype for the proof.

This approach is used in
\S \ref{subsection: weaker waist inequality},
where the maps $F_p$ can be chosen to be translations.
In general, however, $F_p$ need not be a homeomorphism,
so a different approach is needed.

\smallskip
\noindent
\emph{Approach 2.}
The next idea is to use an 
\lq\lq{}approximate solution\rq\rq{}
of the equation $F_p(x) = \Phi(p)$.

We can choose $0<\delta<h$ such that
\[
\rho_n(F_p(x)-\Phi(p))>\delta
\]
for every $(p,x)\in P\times L$.
Indeed, the function
$(p,x)\mapsto\rho_n(F_p(x)-\Phi(p))$
is strictly positive on $P\times L$.
Moreover, it tends to infinity as $\norm{x}\to\infty$,
uniformly in $p\in P$, because $P$ is compact and
$\rho_n(F_p(x)-x)\leq R$.
Hence it is bounded away from zero on $P\times L$.

For $p\in P$, set 
\[ U_p = \{x\in \mathbb{R}^n\mid 
   \rho_n\left(F_p(x) -\Phi(p)\right) < \delta \}. \]
Roughly, this is the set of 
approximate solutions
of the equation $F_p(x) = \Phi(p)$.
Each $F_p$ is a bounded perturbation of the identity,
hence is proper of degree one and therefore surjective.
In particular, $U_p$ is nonempty.

Suppose that there is a continuous map 
$f\colon P\to \mathbb{R}^n$ satisfying 
$f(p) \in U_p$ for all $p\in P$.
Then $f(p)\notin L$ by the choice of $\delta$.
Hence we can define
$\varphi:=r\circ f\colon P\longrightarrow K$.
This satisfies
 \begin{equation*}
  \begin{split}
   \rho_n\left(\varphi(p)-\Phi(p)\right) & \leq
   \rho_n\left(\varphi(p)-f(p)\right) + 
   \rho_n\left(f(p)-F_p(f(p))\right) +
   \rho_n\left(F_p(f(p))-\Phi(p)\right) \\
   &< \frac{h}{2}+R+\delta < R+ 2h.
  \end{split}
 \end{equation*}
Thus we obtain the desired map also in this case.
The crucial point is 
how to choose $f(p)\in U_p$ continuously in 
$p\in P$.
If every $U_p$ is $(\dim P-1)$-connected, 
this may be possible 
as in \S \ref{section: proof of proposition avoidance}.
However, the sets $U_p$ need not even be
path-connected in general.
So this approach
does not provide a general construction.

\smallskip
\noindent
\emph{Approach 3.}
In the third approach, we try to construct
$\varphi\colon P\to K$
inductively over the skeleta $P^{(j)}$ $(j\geq0)$.

For each vertex $p\in P^{(0)}$,
choose a point $x_p\in U_p$ and set
$\varphi(p):=r(x_p)\in K$.
We would like to extend $\varphi$ successively
to the higher-dimensional skeleta of $P$.
The key observation is that
$K=K(n,d,h)$ is $(d-1)$-connected.
Hence we can extend $\varphi$ to $P^{(j)}$
for $j\leq d$.
(Here we ignore the problem of controlling
the approximation error
$\rho_n(\varphi(p)-\Phi(p))$.
This issue will be addressed by using suitable
finite subcomplexes of $K$ in the actual proof.)

If $\dim P\leq d$, this would solve
the extension problem.
However, we only assume that $\dim P\leq2d-1$,
and the $(d-1)$-connectivity of $K$
does not guarantee extensions beyond $P^{(d)}$.
Thus this approach alone is not sufficient.

Notice that this argument makes little use
of the avoidance map $F$.
The main question is how to exploit $F$
to overcome the higher-dimensional extension obstructions.
The obstruction to extending $\varphi$
from $P^{(j)}$ to $P^{(j+1)}$
is represented by elements of the homotopy group $\pi_j(K)$.
These groups vanish for $j\leq d-1$,
but need not vanish for $j\geq d$.
We therefore need additional information to overcome
the obstructions for $d\leq j\leq2d-2$.

Our idea is to introduce an auxiliary space
$\mathcal{E}_n(K)$ and a natural map
$\iota_K\colon K\to\mathcal{E}_n(K)$
inducing isomorphisms on homotopy groups
up to degree $2d-2$.
The avoidance map $F$ will be used to construct
a continuous map $\chi \colon P\to\mathcal{E}_n(K)$.
(The construction of $\chi$ is inspired by 
the viewpoint of Approach 2.)
The crucial point is to establish 
a \lq\lq{}correspondence\rq\rq{}
between the extension obstructions for $\varphi$
and those associated with $\chi$.
Since $\chi$ is defined on all of $P$,
the latter obstructions vanish.
Through this correspondence,
we eliminate
the obstructions for $\varphi$ as well,
and thereby extend $\varphi$ over
$P=P^{(2d-1)}$.

More precisely, 
at each stage of the inductive construction of 
$\varphi$, 
we construct $\varphi$ on $P^{(j)}$ together
with a homotopy
\[
  \iota_K\circ\varphi \simeq \chi|_{P^{(j)}}.
\]
For each $(j+1)$-simplex $\sigma$,
this homotopy identifies the image of the boundary
class under $(\iota_K)_*$ with the class of
$\chi|_{\partial\sigma}$.
The latter is zero because $\chi$ is already defined
on all of $\sigma$.
The injectivity of $(\iota_K)_*$ therefore implies
that the boundary map extends over $\sigma$.
(See Lemma \ref{lemma: relative extension} below
for the details.)

This describes the topological mechanism of the proof.
To control the approximation error, we carry out
a local version of this construction, using suitable
$(d-1)$-connected subcomplexes of $K$
near the images of the simplices under $\Phi$.

\subsection{Preliminaries on algebraic topology}
\label{subsection: preliminaries on algebraic topology}

In this subsection, we first recall some basic notions
and results from algebraic topology.
We then introduce an auxiliary space $\mathcal E_n(C)$
and establish a connectivity estimate for its canonical map.
These results will be used in the proof of Proposition
\ref{proposition: low dimensional approximation}.

Let $I=[0,1]$ denote the unit interval, and let
$D^n=\{x\in\mathbb{R}^n\mid\norm{x}\leq1\}$
denote the closed unit ball in $\mathbb{R}^n$.
We often identify the $n$-dimensional sphere $S^n$ 
with the one-point compactification of $\mathbb{R}^n$:
\[
S^n=\mathbb{R}^n\cup\{\infty\}.
\]
We regard $\infty$ as the basepoint of $S^n$.

For topological spaces $X$ and $Y$, let
$\mathrm{Map}(X,Y)$ denote the space of continuous maps
from $X$ to $Y$, equipped with 
the compact-open topology\footnote{
In general, it is subtle how to 
define a topology on $\mathrm{Map}(X,Y)$.
However, for the proof of Proposition 
\ref{proposition: low dimensional approximation},
we only need to consider the case that 
$X$ is a compact metrizable space and 
$Y$ is a metrizable space. In this case, 
$\mathrm{Map}(X,Y)$ 
is also metrizable and 
the situation becomes simple.   
Indeed, when $X$ is compact and $Y$ is metrizable,
the compact-open topology on 
$\mathrm{Map}(X,Y)$ is induced by the metric
\[
D(f,g)=\sup_{x\in X}
\mathbf{d}_Y\bigl(f(x),g(x)\bigr),
\]
where $\mathbf{d}_Y$ is any compatible metric on $Y$.
}.
When $X$ and $Y$ have basepoints $x_0$ and
$y_0$ respectively, let $\mathrm{Map}_*(X,Y)$ denote
the subspace consisting of continuous maps 
$f\colon X\to Y$ satisfying
$f(x_0)=y_0$.

We will use the following standard property 
of the compact-open topology: 
if $X$ is a locally compact Hausdorff space and 
$Y,Z$ are topological spaces, 
then a map $f\colon Z\to\mathrm{Map}(X,Y)$ is continuous 
if and only if the associated map
\[
\widehat f\colon Z\times X\to Y,\qquad
\widehat f(z,x)=f(z)(x),
\]
is continuous; see \cite[Proposition~A.14]{Hatcher}.
In particular, this applies when $X=I$ or $X=S^n$.

Let $f\colon X\to Y$ be a continuous map
between topological spaces.
For $y\in Y$, we define the \textbf{homotopy fiber}
of $f$ over $y$ by
\[
\mathrm{hofib}_y(f)
=
\{(x,\gamma)\in X\times\mathrm{Map}(I,Y)
\mid \gamma(0)=f(x),\ \gamma(1)=y\}.
\]
We equip this space with the subspace topology
inherited from $X\times\mathrm{Map}(I,Y)$.
Here $x$ ranges over those points of $X$ for which
$f(x)$ lies in the same path component as $y$.
The homotopy fiber is a homotopy-invariant replacement
for the ordinary fiber $f^{-1}(y)$:
the equality $f(x)=y$ is replaced by the data of a path
from $f(x)$ to $y$.

For an integer $k\geq 0$, 
the map $f\colon X\to Y$ is said to be
\textbf{$k$-connected} if
$\mathrm{hofib}_y(f)$ is $(k-1)$-connected
for every $y\in Y$.
(Here, $(-1)$-connected means nonempty.)
For $k=0$, this is equivalent to the surjectivity of
the induced map
\[
f_*\colon \pi_0(X)\longrightarrow\pi_0(Y),
\]
where $\pi_0$ denotes the set of path components.
For $k\geq1$, it is equivalent to the condition that
$f_*\colon\pi_0(X)\to\pi_0(Y)$ is bijective and,
for every $x\in X$, the induced homomorphism
\[
f_*\colon \pi_j(X,x)\longrightarrow\pi_j(Y,f(x))
\]
is an isomorphism for $1\leq j<k$ and
surjective for $j=k$.
(Here $\pi_j(X,x)$ denotes the $j$-th homotopy group
of $X$ with basepoint $x$.)
The case $k=0$ follows directly from the definition
of the homotopy fiber.
For $k\geq1$, the equivalence follows from the long
exact homotopy sequence 
of the mapping path fibration, applied over
each path component of $Y$;
see \cite[Theorem~4.41 and Proposition~4.64]{Hatcher}.

This characterization implies that, 
if $f\colon X\to Y$ and $g\colon Y\to Z$ are 
$k$-connected continuous 
maps, then 
$g\circ f\colon X\to Z$ is also $k$-connected.
Note also that, if $k^\prime \geq k$ 
and $f$ is $k^\prime$-connected, 
then it is $k$-connected.

The next lemma will be crucial 
in the proof of Proposition 
\ref{proposition: low dimensional approximation}.
Here $\partial D^j = \{x\in \mathbb{R}^j \mid 
\norm{x}=1\}$ is the boundary of $D^j$.

\begin{lemma}[Relative extension lemma]
\label{lemma: relative extension}
Let $1\leq j\leq k$ be integers, and let
$f\colon X\to Y$ be a $k$-connected continuous map.
Suppose we are given continuous maps
\[
\varphi\colon \partial D^j\to X,\qquad
\chi\colon D^j\to Y,\qquad
H\colon \partial D^j\times I\to Y
\]
satisfying
\[
H(x,0)=f\bigl(\varphi(x)\bigr),\qquad
H(x,1)=\chi(x)
\]
for every $x\in\partial D^j$.
Then $\varphi$ and $H$ admit continuous extensions
\[
\varphi\colon D^j\to X,\qquad
H\colon D^j\times I\to Y
\]
satisfying
\[
H(x,0)=f\bigl(\varphi(x)\bigr),\qquad
H(x,1)=\chi(x)
\]
for every $x\in D^j$.
\end{lemma}

\begin{proof}
This is a special case of the homotopy extension
and lifting property (HELP) for the pair
$(D^j,\partial D^j)$;
see \cite[Chapter~10, \S3]{May}.
We include a proof for 
completeness.

Set $y_0=\chi(0)$ and
$F=\mathrm{hofib}_{y_0}(f)$.
For each $x\in\partial D^j$, define a path
$\alpha_x\colon I\to Y$ by
\[
\alpha_x(t)=
\begin{cases}
H(x,2t), & 0\leq t\leq \frac12,\\
\chi((2-2t)x), & \frac12\leq t\leq1.
\end{cases}
\]
This path starts at $f(\varphi(x))$ and ends at $y_0$.
Hence
\[
\partial D^j\longrightarrow F,\qquad
x\longmapsto \bigl(\varphi(x),\alpha_x\bigr)
\]
is a continuous map, where the path space is equipped
with the compact-open topology.

Since $f$ is $k$-connected, the space $F$ is
$(k-1)$-connected.
As $j\leq k$, the above map extends to a continuous map
\[
D^j\longrightarrow F,\qquad
x\longmapsto
\bigl(\widetilde\varphi(x),\beta_x\bigr).
\]

In particular,
\[
\widetilde\varphi|_{\partial D^j}=\varphi,
\qquad
\beta_x=\alpha_x
\quad\text{for }x\in\partial D^j.
\]

Define $G\colon D^j\times I\to Y$ by
\[
G(x,t)=
\begin{cases}
\beta_x(2t), & 0\leq t\leq\frac12,\\
\chi((2t-1)x), & \frac12\leq t\leq1.
\end{cases}
\]
(The two formulas agree at $t=1/2$, since
$\beta_x(1)=y_0=\chi(0)$.)
This map is continuous and satisfies
\[
G(x,0)=f\bigl(\widetilde\varphi(x)\bigr),
\qquad
G(x,1)=\chi(x).
\]
For $x\in\partial D^j$, the path $G(x,\cdot)$
first follows $H(x,\cdot)$, then follows the path
$t\mapsto \chi((1-t)x)$ from $\chi(x)$ to $y_0$,
and finally follows its reverse.
(Here we ignore the parametrizations 
of these paths.\footnote{
More precisely, 
for $x\in\partial D^j$, we have
\[
G(x,t)=
\begin{cases}
H(x,4t),
& 0\leq t\leq\frac14,\\
\chi((2-4t)x),
& \frac14\leq t\leq\frac12,\\
\chi((2t-1)x),
& \frac12\leq t\leq1.
\end{cases}
\]
})
Contracting this backtracking part and reparametrizing
gives a homotopy from
$G|_{\partial D^j\times I}$ to $H$,
relative to $\partial D^j\times\{0,1\}$.
This homotopy depends continuously on $x$.

Set
\[
A=(\partial D^j\times I)\cup(D^j\times\{0,1\}).
\]
Together with the constant homotopies on
$D^j\times\{0,1\}$, the preceding homotopy defines
a homotopy on $A$ starting at $G|_A$.
The pair $(D^j\times I,A)$ is a CW pair and therefore
has the homotopy extension property
\cite[Proposition~0.16]{Hatcher}.\footnote{
Set $B = D^j \times I$. Then 
$B\times \{0\}\cup A\times I$ is a deformation 
retract of $B\times I$.
The homotopy extension property follows from this.
}
Consequently, this homotopy extends to a homotopy
on $D^j\times I$ starting at $G$.
Its final map
$\widetilde H\colon D^j\times I\to Y$ satisfies
\[
\widetilde H|_{\partial D^j\times I}=H,
\qquad
\widetilde H(x,0)=f\bigl(\widetilde\varphi(x)\bigr),
\qquad
\widetilde H(x,1)=\chi(x).
\]
Thus $\widetilde\varphi$ and $\widetilde H$
are the required extensions.
\end{proof}

Let $X$ and $Y$ be topological spaces with
basepoints $x_0$ and $y_0$, respectively.
We define their \textbf{wedge sum} $X\vee Y$
as the disjoint union with the two basepoints
identified:
\[
X\vee Y=(X\sqcup Y)/(x_0\sim y_0).
\]
We identify $X\vee Y$ with the subspace
\[
(X\times\{y_0\})\cup(\{x_0\}\times Y)
\]
of $X\times Y$.
We define their \textbf{smash product} by
collapsing this subspace to a single point:
\[
X\wedge Y=(X\times Y)/(X\vee Y).
\]
Both spaces are equipped with the quotient topology.
The identified point in $X\vee Y$ and the collapsed
point in $X\wedge Y$ serve as their basepoints.

For an integer $n\geq1$, we define the
\textbf{$n$-fold reduced suspension} of $X$ by
\[
\Sigma^n X=S^n\wedge X,
\]
where $S^n=\mathbb{R}^n\cup\{\infty\}$ is the
one-point compactification of $\mathbb{R}^n$,
with basepoint $\infty$.
We write $\Sigma X=\Sigma^1 X$.
There is a natural homeomorphism
\begin{equation}
  \label{eq: iterated reduced suspension}
 \underbrace{\Sigma(\Sigma 
  \dots (\Sigma}_{\text{$n$-fold}} X)\dots) 
 \cong
\Sigma^nX
\end{equation}
which we use to identify these spaces.
(For example, $\Sigma(\Sigma(\Sigma X))
\cong\Sigma^3 X$.)
Explicitly, this homeomorphism sends
\[
[t_1,[t_2,\ldots,[t_n,x]\ldots]]
\longmapsto [(t_1,\ldots,t_n),x]
\qquad
(t_1,\ldots,t_n\in\mathbb{R},\ x\in X),
\]
and sends the basepoint to the basepoint.

If $X$ is a $k$-connected CW complex with 
a chosen $0$-cell $x_0$ as its basepoint,
then $\Sigma^n X$ is $(n+k)$-connected; 
see \cite[(6.7.10) Corollary]{tomDieck}.

The \textbf{$n$-fold loop space} of $X$ 
is defined to be 
the space of basepoint-preserving 
continuous maps 
from $S^n$ to $X$:
\[ \Omega^n X = \mathrm{Map}_*(S^n, X). \]
When $n=1$, we simply write 
$\Omega X := \Omega^1 X$.
For $n\geq2$, we naturally identify
$\Omega(\Omega^{n-1}X)$ with $\Omega^nX$
by sending $\alpha\colon S^1\to\Omega^{n-1}X$
to the map
\[
S^n=S^1\wedge S^{n-1}\longrightarrow X,
\qquad [t,v]\longmapsto\alpha(t)(v).
\]

Let $c\in \Omega^n X$ be the constant map 
that sends $S^n$ to $\{x_0\}$.
For every integer $k\geq 1$, there is a 
natural isomorphism
\[ \pi_k\left(\Omega^n X, c\right)
   \to \pi_{k+n}(X, x_0). \]
It sends the homotopy class of 
$\alpha\colon S^k \to \Omega^n X$ to 
the homotopy class of
\[ S^{n+k} \cong S^n \wedge S^k \to X, \quad 
  [u, v] \mapsto \alpha(v)(u). \]
For $k=0$, there is a natural bijection 
between the set of 
path components of $\Omega^n X$ 
and the $n$-th homotopy group of $(X,x_0)$:
\[ \pi_0\left(\Omega^n X\right) \cong 
   \pi_n(X, x_0). \]

If $k\geq n \geq 1$ and $X$ is $k$-connected, 
then $\Omega^n X$ is $(k-n)$-connected.
This follows from the above isomorphisms 
of homotopy groups 
and the fact that $\pi_j(X,x_0)=0$ for $j\leq k$.

If $f\colon (X,x_0)\to (Y,y_0)$ is a $k$-connected
basepoint-preserving continuous map for an integer
$k\geq1$, then
\[
\Omega f\colon \Omega X\longrightarrow\Omega Y,
\qquad
\alpha\longmapsto f\circ\alpha,
\]
is $(k-1)$-connected.
Indeed, let $c_X$ and $c_Y$ denote the constant loops
at $x_0$ and $y_0$, respectively.
For every $j\geq0$, there is a natural commutative
diagram
\[
\begin{tikzcd}[column sep=large]
\pi_j(\Omega X,c_X)
  \arrow[r,"(\Omega f)_*"]
  \arrow[d,"\cong"']
&
\pi_j(\Omega Y,c_Y)
  \arrow[d,"\cong"]
\\
\pi_{j+1}(X,x_0)
  \arrow[r,"f_*"']
&
\pi_{j+1}(Y,y_0).
\end{tikzcd}
\]
For $j=0$, the vertical arrows are bijections
identifying the path components of the loop spaces
with the elements of the fundamental groups.
If $k=1$, then
$f_*\colon\pi_1(X,x_0)\to\pi_1(Y,y_0)$ is surjective.
The diagram with $j=0$ therefore shows that
$\Omega f$ induces a surjection on path components.
Thus $\Omega f$ is $0$-connected.
Now suppose that $k\geq2$.
Since $f_*\colon\pi_1(X,x_0)\to\pi_1(Y,y_0)$ is
an isomorphism, the diagram with $j=0$ shows that
$\Omega f$ induces a bijection on path components.
For $j\geq1$, the same diagram shows that
\[
(\Omega f)_*\colon
\pi_j(\Omega X,c_X)\longrightarrow
\pi_j(\Omega Y,c_Y)
\]
is an isomorphism for $1\leq j<k-1$ and
surjective for $j=k-1$.
These assertions also hold at every other basepoint
of $\Omega X$, since concatenation with a fixed loop
is a homotopy equivalence and $\Omega f$ preserves
concatenation.
Hence $\Omega f$ is $(k-1)$-connected.

The following lemma follows by iterating
the Freudenthal suspension theorem.

\begin{lemma}[Suspension theorem]
\label{lemma: suspension theorem}
Let $X$ be a finite CW complex
with a chosen $0$-cell $x_0$ as its basepoint.
Suppose that $X$ is $k$-connected
for some integer $k\geq0$.
Then, for every integer $n\geq1$, the map
\[
\theta_n\colon X\longrightarrow\Omega^n\Sigma^nX,
\qquad
\theta_n(x)(z)=[z,x]
\quad (z\in S^n),
\]
is $(2k+1)$-connected.
\end{lemma}

\begin{proof}
This statement appears in 
\cite[Remark 5.10]{BlomquistHarper}.

First notice that
$\Sigma^n X$ is $(n+k)$-connected.
The standard Freudenthal suspension theorem 
states that 
\[ \theta_1 \colon X\to \Omega \Sigma X \]
is $(2k+1)$-connected
(\cite[(6.10.4) Theorem]{tomDieck}, 
\cite[Corollary 4.24]{Hatcher}).

By applying the Freudenthal theorem to
$\Sigma X$, we obtain that the natural map 
\[ \Sigma X \to \Omega \Sigma^2 X \]
is $(2k+3)$-connected, and hence 
\begin{equation}
  \label{eq: Omega Sigma to Omega^2 Sigma^2}
   \Omega \Sigma X \to \Omega^2 \Sigma^2 X
\end{equation}   
is $(2k+2)$-connected.

The composite of $\theta_1$ with the map
\eqref{eq: Omega Sigma to Omega^2 Sigma^2}
sends $x\in X$ to the map given by
\[
(t_1,t_2)\longmapsto [t_2,[t_1,x]]
=[(t_2,t_1),x]
\qquad ((t_1,t_2)\in\mathbb{R}^2\subset S^2).
\]
Thus this composite agrees with $\theta_2$ up to
precomposition with the coordinate-switching homeomorphism
\[
\tau\colon S^2\to S^2,\qquad
\tau(t_1,t_2)=(t_2,t_1),\qquad \tau(\infty)=\infty.
\]
Since precomposition with $\tau$ induces a homeomorphism
of $\Omega^2\Sigma^2X$, the composite and $\theta_2$
have the same connectivity.
It follows that $\theta_2$ is $(2k+1)$-connected.

We can proceed similarly and prove that all 
$\theta_n$ are $(2k+1)$-connected.
\end{proof}

Let $C$ be a non-empty finite CW complex 
(\emph{initially with no 
chosen basepoint}). 
Add a disjoint basepoint $*$ to it and form 
\[ C_+ = C \sqcup\{*\}. \]
We consider the $n$-fold reduced suspension
$\Sigma^n C_+ = S^n\wedge C_+$.
Then
\begin{align*}
\Sigma^n C_+
&= (S^n\times C_+)/(S^n\vee C_+) \\
&\cong (S^n\times C)/(\{\infty\}\times C).
\end{align*}
This can be seen as the one-point compactification of 
$\mathbb{R}^n\times C$.
We may think of a point of $\Sigma^n C_+$ as a
point of $S^n$ carrying a \lq\lq{}label\rq\rq{} 
in $C$.
More precisely, every non-basepoint is uniquely
represented by a pair $(x,c)$ with
$x\in \mathbb{R}^n$ and $c\in C$,
where $c$ is regarded as the label attached to $x$.
At $\infty$, all labels are identified:
the pairs $(\infty,c)$, for $c\in C$, collapse
to the single basepoint of $\Sigma^n C_+$.

Let
\[
q_C\colon \Sigma^n C_+ 
= S^n\wedge C_+\longrightarrow S^n,
\qquad [x,c]\longmapsto x,
\]
be the map that forgets the label.
Here $x\in S^n$ and $c\in C$; in particular,
$q_C$ sends the basepoint to $\infty$.

Fix a basepoint-preserving continuous map
$b\colon S^n\to S^n$ that is homotopic to
the identity map $\mathrm{id}_{S^n}$ 
through basepoint-preserving maps.
Equivalently, choose a point $b \in\Omega^n S^n$
in the path component containing $\mathrm{id}_{S^n}$.
Consider the map 
\[ \Omega^n(q_C)\colon 
   \Omega^n(\Sigma^n C_+) \longrightarrow
   \Omega^n(S^n), \qquad
   u \longmapsto  q_C\circ u. \]
Define $\mathcal{E}_n(C)$ to be the homotopy fiber of 
this map over the point $b$:
\[ \mathcal{E}_n(C) 
= \mathrm{hofib}_b \left(\Omega^n (q_C)\right). \]
Explicitly,
\[
\mathcal{E}_n(C)
=
\left\{
(u,\gamma)\in
\Omega^n \Sigma^n C_+
\times \mathrm{Map}(I,\Omega^n S^n)
\;\middle|\;
\gamma(0)=q_C\circ u,\quad \gamma(1)=b
\right\}.
\]
Thus an element of $\mathcal{E}_n(C)$ consists of
a basepoint-preserving map
$u\colon S^n\to \Sigma^n C_+$
together with a specified homotopy from
$q_C\circ u$ to $b$ through basepoint-preserving maps.
More concretely, the path $\gamma$ corresponds to
the continuous map
\[
H\colon S^n\times I\longrightarrow S^n,
\qquad H(x,t)=\gamma(t)(x),
\]
which satisfies
\[
H(x,0)=q_C(u(x)),\qquad
H(x,1)=b(x),\qquad
H(\infty,t)=\infty
\]
for all $x\in S^n$ and $t\in I$.
We sometimes write $\mathcal{E}_n(C,b)$ to emphasize
the dependence on the choice of $b$.
Since all such choices of $b$ lie in the same path
component of $\Omega^n S^n$, the homotopy type of
$\mathcal{E}_n(C,b)$ is independent of this choice.

There is a canonical map
\begin{equation}
  \label{eq: canonical map in low dimensional approximation}
\iota_C\colon C\longrightarrow\mathcal{E}_n(C),
\qquad
c\longmapsto(u_c,\gamma_b),
\end{equation}
where
\[
u_c\colon S^n\longrightarrow\Sigma^n C_+,
\qquad u_c(x)=[b(x),c],
\]
and $\gamma_b\colon I\to \Omega^n S^n$ 
is the constant path at
$b\in\Omega^n S^n$.
Since $q_C\circ u_c=b$, this map is well-defined.
Informally, $\iota_C(c)$ is obtained by attaching
the same label $c$ to every non-basepoint value of $b$.

A similar idea of introducing an auxiliary space
and a canonical map appears in
\cite[\S 3]{KleinWilliams07}.
See also \cite{Larmore72}.

The following lemma is a key ingredient 
of the proof of Proposition 
\ref{proposition: low dimensional approximation}.

\begin{lemma}
  \label{lemma: connectivity of canonical map}
Let $k\geq 0$ and $n\geq k+2$ be integers.
If $C$ is a nonempty, $k$-connected  
finite CW complex, then the canonical map 
$\iota_C \colon C\to \mathcal{E}_n(C)$ is 
$(2k+1)$-connected. 
\end{lemma}

\begin{proof}
Choose a $0$-cell $c_0$ of $C$
as its basepoint.
The basepoint of $C_+=C\sqcup\{*\}$ remains
the adjoined point $*$, which is distinct from $c_0$.

Let
\[
\rho\colon\Sigma^n C_+\longrightarrow\Sigma^n C
\]
be the quotient map that collapses the subspace
\[
\{[z,c_0]\mid z\in S^n\}\cong S^n
\]
to the basepoint.
We consider the map
\[
F=(q_C,\rho)\colon
\Sigma^n C_+\longrightarrow S^n\times\Sigma^n C.
\]
The purpose of Steps 1 and 2 below is to show that 
$F$ is $(2n+k)$-connected.

\smallskip
\noindent
\emph{Step 1: A wedge decomposition compatible with $F$.}

We construct a basepoint-preserving homotopy equivalence
\[
\chi\colon
\Sigma^n C_+\longrightarrow S^n\vee\Sigma^n C
\]
such that
\[
j\circ\chi\simeq F,
\]
where
\[
j\colon S^n\vee\Sigma^n C
\hookrightarrow S^n\times\Sigma^n C
\]
is the natural inclusion.
Here $\simeq$ means that the two maps are homotopic.
All homotopies in this step preserve the basepoint.

We first carry out the construction for one suspension.
Let
\[
q_1\colon\Sigma C_+\longrightarrow S^1,
\qquad
\rho_1\colon\Sigma C_+\longrightarrow\Sigma C
\]
denote the corresponding maps.

Let $SC$ be the \textbf{unreduced suspension} 
of $C$,
namely the quotient of $I\times C$ obtained by
collapsing $\{0\}\times C$ and $\{1\}\times C$
to two distinct points 
$\mathbf{N}$ and $\mathbf{S}$, 
respectively.

Inside $SC$, consider the arc
\[
L=\{[t, c_0] \mid 0\leq t\leq1\}.
\]
Attach another interval $e$ to $SC$ by identifying
$e(0)$ with $\textbf{N}$ and 
$e(1)$ with $\textbf{S}$, and write
\[
Z=SC\cup e.
\]
We take $\mathbf{N}$ 
as the basepoint of $Z$.
See Figure \ref{fig:suspension-two-arcs}.

\begin{figure}[!hbp]
  \centering
  \includegraphics[width=0.55\textwidth]
    {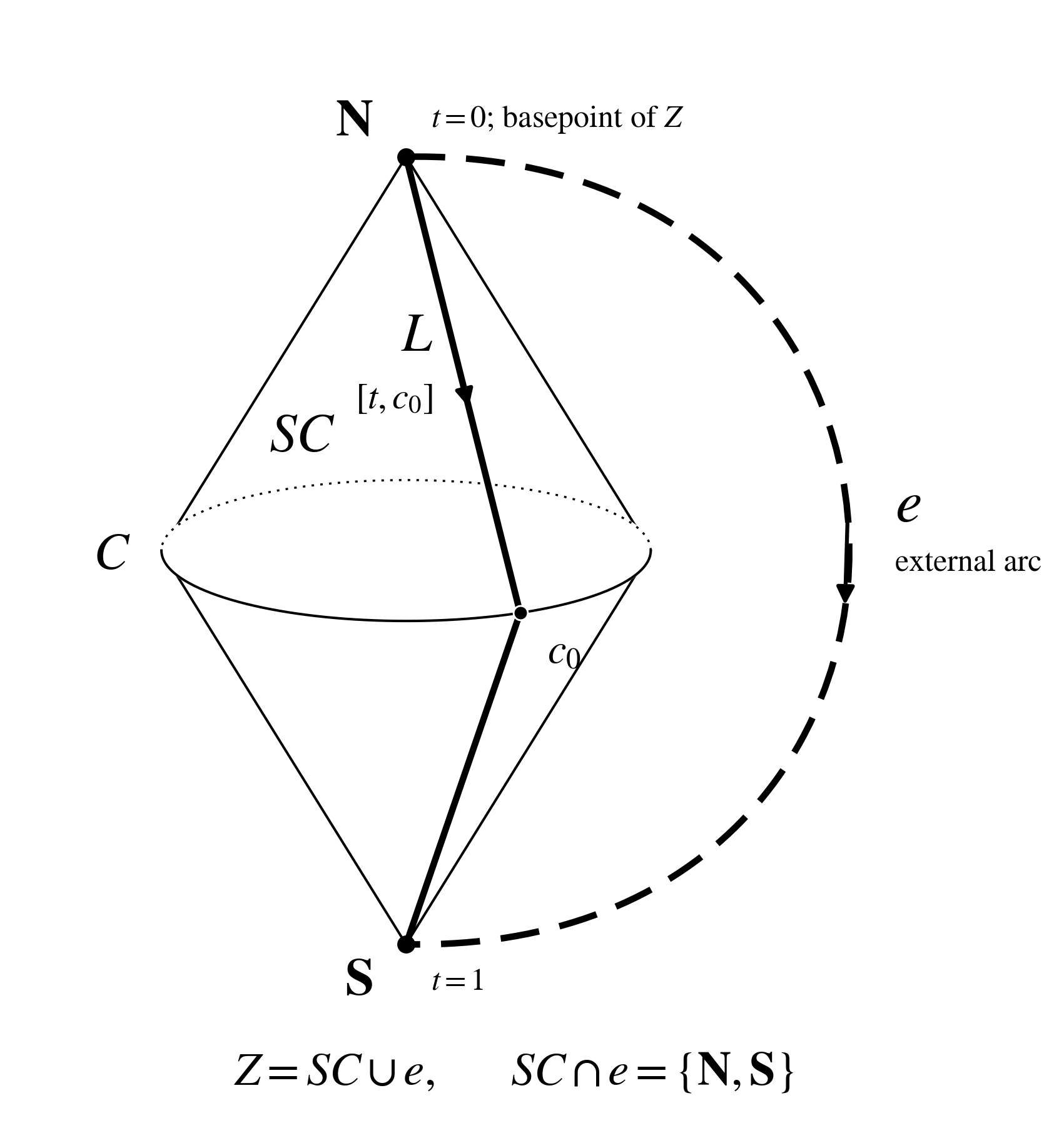}
  \caption{
    A schematic illustration of $Z=SC\cup e$.
    The internal arc
    $L=\{[t,c_0]\mid 0\leq t\leq1\}$
    is drawn as a thick solid line,
    and the external arc $e$ as a thick dashed line.
    Both arcs join $\mathbf{N}$ to $\mathbf{S}$,
    and the arrows indicate increasing $t$.
    We take $\mathbf{N}$ as the basepoint of $Z$.
  }
  \label{fig:suspension-two-arcs}
\end{figure}

Collapsing the external arc $e$ identifies the two
poles, giving a based homeomorphism
\[
Z/e\cong\Sigma C_+.
\]
Collapsing the internal arc $L$ turns $SC$ into
$\Sigma C$ and the external arc $e$ into a circle,
giving a based homeomorphism
\[
Z/L\cong S^1\vee\Sigma C.
\]
Let
\[
p\colon Z\longrightarrow\Sigma C_+,
\qquad
v\colon Z\longrightarrow S^1\vee\Sigma C
\]
be the corresponding quotient maps,
with the identifications specified below.
Both $p$ and $v$ are based homotopy equivalences,
since each collapsed arc is a subcomplex of $Z$
that contracts to $\mathbf{N}$ 
while fixing $\mathbf{N}$;
see \cite[Propositions~0.16--0.17]{Hatcher}.

For the following comparison, identify $S^1$ with
the unit circle in $\mathbb C$, based at $1$.
Choose the quotient identifications so that
\[
p([t,c])=[e^{2\pi i t},c],
\qquad
v(e(t))=e^{-2\pi i t}
\in S^1\subset S^1\vee\Sigma C.
\]
Here $e(t)$ denotes the point on the external arc $e$
corresponding to the parameter $t\in[0,1]$.
On $SC$, the map $v$ is the usual quotient map
to $\Sigma C$, given by
\[
v([t,c])=[e^{2\pi i t},c]
\in\Sigma C\subset S^1\vee\Sigma C.
\]
Note that we choose to parametrize the circle
arising from the external arc $e$
in the reverse direction.

Let
\[
r_S\colon S^1\vee\Sigma C\longrightarrow S^1,
\qquad
r_\Sigma\colon S^1\vee\Sigma C\longrightarrow\Sigma C
\]
be the maps collapsing the other wedge summand.

The $\Sigma C$-components agree:
both $\rho_1\circ p$ and $r_\Sigma\circ v$
collapse $e$ to the basepoint and restrict to the
usual quotient map $SC\to\Sigma C$.
Hence
\[
\rho_1\circ p=r_\Sigma\circ v.
\]

For the $S^1$-components, the two maps are given by
\[
\begin{array}{c|cc}
 & [t,c]\in SC & e(t)\in e \\ \hline
q_1\circ p & e^{2\pi i t} & 1 \\[2pt]
r_S\circ v & 1 & e^{-2\pi i t}
\end{array}
\]
Interpolating between the exponents, define
$H_s\colon Z \to S^1$ by
\[
H_s([t,c])=e^{2\pi i(1-s)t},
\qquad
H_s(e(t))=e^{-2\pi i st}
\qquad (0\leq s\leq1).
\]
Throughout this homotopy, 
both formulas send $\mathbf{N}$
to $1$ and $\mathbf{S}$ 
to $e^{-2\pi i s}$.
They therefore fit together to define a continuous
homotopy on $Z$ fixing the basepoint $\mathbf{N}$.
Since
\[
H_0=q_1\circ p,
\qquad
H_1=r_S\circ v,
\]
we obtain
\[
q_1\circ p\simeq r_S\circ v.
\]

Choose a based homotopy inverse
$h\colon\Sigma C_+\to Z$ of $p$, and put
\[
\chi_1=v\circ h\colon
\Sigma C_+\longrightarrow S^1\vee\Sigma C.
\]
Then $\chi_1$ is a based homotopy equivalence, and
\[
r_S\circ\chi_1\simeq q_1,
\qquad
r_\Sigma\circ\chi_1\simeq\rho_1.
\]

Suspend these two relations separately $n-1$ times.
Using the natural identifications
\[
\Sigma^{n-1}(\Sigma C_+)\cong\Sigma^n C_+,
\qquad
\Sigma^{n-1}(S^1\vee\Sigma C)
\cong S^n\vee\Sigma^n C,
\]
we obtain a based homotopy equivalence
\[
\chi=\Sigma^{n-1}\chi_1\colon
\Sigma^n C_+\longrightarrow S^n\vee\Sigma^n C
\]
whose two projections are homotopic to $q_C$
and $\rho$, respectively.
Combining these homotopies gives
\[
j\circ\chi\simeq(q_C,\rho)=F.
\]

\smallskip
\noindent
\emph{Step 2: The connectivity of $F$.}

Since $C$ is $k$-connected, $\Sigma^n C$ is
$(n+k)$-connected.
By the standard CW approximation construction,
there is a based CW complex $W$ and a based
homotopy equivalence
\[
g\colon W\longrightarrow\Sigma^n C
\]
such that the basepoint of $W$ is its single $0$-cell
and all its other cells have dimension at least $n+k+1$;
see \cite[Proposition~4.15]{Hatcher}.
The complex $W$ need not be finite.

Give $S^n$ its CW structure with one $0$-cell
$\{\infty\}$
and one $n$-cell.
The relative cells of
\[
(S^n\times W,\,S^n\vee W)
\]
are precisely the products of the $n$-cell of
$S^n$ with the positive-dimensional cells of $W$.
Their dimensions are therefore at least
\[
n+(n+k+1)=2n+k+1.
\]
The usual product topology has this CW structure
because $S^n$ is a finite CW complex.
By cellular approximation, the inclusion
\[
S^n\vee W\hookrightarrow S^n\times W
\]
is $(2n+k)$-connected;
see \cite[Corollary~4.12]{Hatcher}.
Consequently,
\[
j\colon S^n\vee\Sigma^n C
\hookrightarrow S^n\times\Sigma^n C
\]
is also $(2n+k)$-connected.
Since $\chi$ is a homotopy equivalence and
$j\circ\chi\simeq F$, the map 
$F = (q_C, \rho)\colon \Sigma^n C_+
\to S^n \times \Sigma^n C$ is
$(2n+k)$-connected.

\smallskip
\noindent
\emph{Step 3: Comparison with the suspension map.}

Applying the loop functor $n$ times, we obtain an
$(n+k)$-connected map
\[
\Omega^n F\colon
\Omega^n(\Sigma^n C_+)\longrightarrow
\Omega^n S^n\times\Omega^n\Sigma^n C,
\]
where we use the natural identification of the
loop space of a product with the product of
the loop spaces.

Define
\[
T\colon\mathcal E_n(C,b)
\longrightarrow\Omega^n\Sigma^n C,
\qquad
T(u,\gamma)=\rho\circ u.
\]
For each $w\in\Omega^n\Sigma^n C$, 
there is a natural
homeomorphism
\[
\mathrm{hofib}_w(T)
\cong
\mathrm{hofib}_{(b,w)}(\Omega^n F).
\]
Indeed, a point of the left-hand side is a triple
$(u,\gamma,\delta) \in \Omega^n \Sigma^n C_+
\times
\mathrm{Map}(I, \Omega^n S^n)
\times
\mathrm{Map}(I, \Omega^n \Sigma^n C)$, 
where
\[
\gamma(0)=q_C\circ u,\quad \gamma(1)=b,
\qquad
\delta(0)=\rho\circ u,\quad \delta(1)=w.
\]
The corresponding point of the right-hand side is
\[
\bigl(u,t\longmapsto(\gamma(t),\delta(t))\bigr).
\]
It follows that $T$ is $(n+k)$-connected.

By the definition 
\eqref{eq: canonical map in low dimensional approximation}
of $\iota_C\colon 
C\to\mathcal E_n(C,b)$,
\[
T(\iota_C(c))(z)=[b(z),c]\in\Sigma^n C.
\]
Recall that $b\in \Omega^n S^n$ is homotopic 
to $\mathrm{id}_{S^n}$ through 
basepoint-preserving maps.
Then, 
a based homotopy from $b$ to $\mathrm{id}_{S^n}$
gives a based homotopy
\[
T\circ\iota_C\simeq\theta_n,
\]
where
\[
\theta_n\colon C\longrightarrow\Omega^n\Sigma^n C,
\qquad
\theta_n(c)(z)=[z,c].
\]
Here $c_0$ is mapped to the constant map throughout
the homotopy, since $[z,c_0]$ is the basepoint
of $\Sigma^n C$ for every $z\in S^n$.
By Lemma~\ref{lemma: suspension theorem},
$T\circ\iota_C$ is $(2k+1)$-connected.

Finally, the assumption $n\geq k+2$ gives
\[
n+k\geq2k+2.
\]
Hence $T$ induces isomorphisms on homotopy groups
in every degree $1\leq j\leq2k+1$.
Also, $\Omega^n\Sigma^n C$ is path-connected,
since $\Sigma^n C$ is $(n+k)$-connected.
The connectivity of $T$ therefore implies that
$\mathcal E_n(C,b)$ is path-connected.

Taking basepoints $c_0\in C$, 
$\iota_C(c_0) \in \mathcal{E}_n(C,b)$, 
and the
constant map in $\Omega^n\Sigma^n C$, we have
\[
(T\circ\iota_C)_*=T_*\circ(\iota_C)_*\colon 
\pi_j(C) \to \pi_j(\Omega^n\Sigma^n C).
\]
It follows that $(\iota_C)_* \colon 
\pi_j(C) \to \pi_j(\mathcal{E}_n(C,b))$ 
is an isomorphism
for $1\leq j<2k+1$ and surjective for $j=2k+1$.
Thus $\iota_C$ is $(2k+1)$-connected.
\end{proof}

Let $C$ and $C'$ be nonempty finite CW complexes,
and let $f\colon C\to C'$ be a continuous map.
Extend $f$ to a basepoint-preserving map
\[
f_+\colon C_+\longrightarrow C'_+
\]
by setting $f_+(*)=*$.
We also set
\[
\Sigma^n(f_+)
:=\mathrm{id}_{S^n}\wedge f_+
\colon \Sigma^n C_+\longrightarrow\Sigma^n C'_+.
\]
Explicitly, this map sends $[z,c]$ to $[z,f(c)]$
for $z\in S^n$ and $c\in C$, and sends the
basepoint to the basepoint.
This map leaves the $S^n$-coordinate
unchanged and replaces each label $c\in C$
with the label $f(c)\in C'$.

Since
\[
q_{C'}\circ\Sigma^n(f_+)=q_C,
\]
we obtain a continuous map
\[
\mathcal E_n(f)\colon
\mathcal E_n(C,b)\longrightarrow\mathcal E_n(C',b),
\qquad
(u,\gamma)\longmapsto
\bigl(\Sigma^n(f_+)\circ u,\gamma\bigr).
\]
Moreover, the canonical maps satisfy
\begin{equation}
  \label{eq: functoriality of canonical map}
  \mathcal E_n(f)\circ\iota_C
=
\iota_{C'}\circ f.
\end{equation}

For the application below, it is important that
this construction does not require a choice of
basepoint in either $C$ or $C'$.
In the proof of Proposition
\ref{proposition: low dimensional approximation},
we choose a CW complex $K_\sigma$ for each simplex
$\sigma$ of $P$ to control the approximation error.
The naturality relation
\eqref{eq: functoriality of canonical map}
is used to make
the maps constructed on individual simplices
compatible on their common faces, so that they
glue to a continuous map on all of $P$.

\subsection{Proof of Proposition
\ref{proposition: low dimensional approximation}}
\label{subsection: proof of low-dimensional approximation}

Here we complete the proof of Proposition 
\ref{proposition: low dimensional approximation}.
Before going into the proof, we present two 
simple lemmas.
Throughout this subsection, we fix a real number $h>0$ and 
integers $1\leq d <n$.
We set $K = K(n, d, h)$ and $L = L(n, d, h)$.

\begin{lemma}[Retraction]
 \label{lemma: retraction}
  There exists a continuous map 
  $r\colon \mathbb{R}^n\setminus L\to K$ such that 
  \begin{itemize}
    \item $r(x) = x$ for all $x\in K$,
    \item $\rho_n(x-r(x)) \leq \frac{h}{2}$ for all 
          $x\in \mathbb{R}^n \setminus L$.  
  \end{itemize}
  See Figure \ref{fig:radial-retraction} 
  for the case of 
  $n=2$ and $d=1$.
\end{lemma}

\begin{proof}
Let $Q$ be a closed $n$-dimensional 
cube of the original cubical
grid\footnote{
Namely, $Q = p + [0, h]^n$ for some 
$p\in h\mathbb{Z}^n$.  
}, 
and let $c=(c_1,\ldots,c_n)$ be its center.
For $x\in Q$, write
\[
x=c+u,\qquad |u_i|\leq h/2 \quad (1\leq i\leq n).
\]
Let $t_Q(x)$ be the $(d+1)$-st smallest number among
$|u_1|,\ldots,|u_n|$.
This is a continuous function of $x$.
Within $Q$, a point belongs to $L$ if and only if
at least $d+1$ of the coordinates $u_i$ vanish.
Hence $t_Q(x)>0$ for $x\in Q\setminus L$.

For $x\in Q\setminus L$, define
\[
(r_Q(x))_i
=
c_i+
\max\left\{
-\frac h2,\,
\min\left\{
\frac{h u_i}{2t_Q(x)},\,\frac h2
\right\}
\right\},
\qquad 1\leq i\leq n.
\]
This defines a continuous map on $Q\setminus L$.
At least $n-d$ of the numbers $|u_i|$ are
greater than or equal to $t_Q(x)$.
For each such coordinate,
$(r_Q(x))_i$ is an endpoint of the corresponding
coordinate interval of $Q$.
Thus $r_Q(x)$ lies in a face of $Q$ of dimension
at most $d$, and hence belongs to $K$.

We check that these maps agree on common faces.
If $x$ belongs to two grid cubes, the numbers
$|x_i-c_i|$ are independent of which cube is used:
whenever the centers differ in the $i$-th coordinate,
$x_i$ is a common endpoint and
$|x_i-c_i|=h/2$ for both cubes.
Consequently, the value of $t_Q(x)$ is also independent
of this choice.
Since $t_Q(x)\leq h/2$, every coordinate corresponding
to a common endpoint is fixed by the above formula.
The remaining coordinates have the same formula
for both cubes.
Therefore the maps $r_Q$ agree on their overlaps.
The grid cubes form a locally finite closed cover,
so these maps define a continuous map
\[
r\colon\mathbb{R}^n\setminus L\longrightarrow K.
\]

If $x\in K\cap Q$, at least $n-d$ of the numbers
$|u_i|$ equal $h/2$.
It follows that $t_Q(x)=h/2$, and the formula gives
$r(x)=x$.

Finally, since $0<t_Q(x)\leq h/2$, each coordinate
moves away from the center of $Q$ without passing
the corresponding endpoint.
In particular,
\[
|(r_Q(x))_i-x_i|\leq \frac h2
\qquad (1\leq i\leq n).
\]
Therefore
\[
\rho_n(x-r(x))
=
\left(\frac1n\sum_{i=1}^n|x_i-r_i(x)|^2\right)^{1/2}
\leq \frac h2.
\]
\end{proof}

\begin{lemma}
  \label{lemma: Tietze extension}
  Let $R>0$, 
  and let $P$ be a finite simplicial complex.
  Let $F\colon P\times L\to\mathbb{R}^n$ be
  a continuous map satisfying
  \[
    \rho_n(F(p,x)-x)\leq R
  \]
  for all $p\in P$ and $x\in L$.
  Then $F$ admits a continuous extension
  (still denoted by $F$)
  \[
    F\colon P\times\mathbb{R}^n\to\mathbb{R}^n,
  \]
  satisfying
  \[
    \rho_n(F(p,x)-x)\leq R
  \]
  for all $p\in P$ and $x\in\mathbb{R}^n$.
\end{lemma}

\begin{proof}
Define
\[
v\colon P\times L\to\mathbb{R}^n,
\qquad v(p,x)=F(p,x)-x.
\]
Since $P\times\mathbb{R}^n$ is metrizable and
$P\times L$ is closed, the Tietze extension theorem
applies to each coordinate of $v$.
Combining these coordinate extensions, we obtain
a continuous map
\[
\widetilde v\colon P\times\mathbb{R}^n\to\mathbb{R}^n
\]
extending $v$.

Let $\Pi\colon\mathbb{R}^n\to\mathbb{R}^n$
be the radial projection onto the closed Euclidean
ball of radius $R\sqrt n$ centered at the origin:
\[
\Pi(y)=
\begin{cases}
y, & \norm{y}\leq R\sqrt n,\\[2mm]
\dfrac{R\sqrt n}{\norm{y}}\,y,
& \norm{y}>R\sqrt n.
\end{cases}
\]
This map is continuous and satisfies
$\rho_n(\Pi(y))\leq R$ for every $y\in\mathbb{R}^n$.
Moreover, the assumption implies
$\norm{v(p,x)}\leq R\sqrt n$ on $P\times L$,
so $\Pi(v(p,x))=v(p,x)$ there.

We now define
\[
\widetilde F(p,x)=x+\Pi(\widetilde v(p,x)),
\qquad (p,x)\in P\times\mathbb{R}^n.
\]
Then $\widetilde F$ is continuous, agrees with $F$
on $P\times L$, and satisfies
\[
\rho_n(\widetilde F(p,x)-x)
=\rho_n(\Pi(\widetilde v(p,x)))
\leq R
\]
for every $(p,x)\in P\times\mathbb{R}^n$.
This proves the lemma.
\end{proof}

We now prove Proposition
\ref{proposition: low dimensional approximation},
which we restate for convenience.

\begin{restatedapproximation}
[From avoidance to low-dimensional approximation]
Let $1\leq d<n$ be integers, and let $h,R>0$.
Let $P$ be a finite simplicial complex satisfying
\[
\dim P\leq 2d-1.
\]
Let
\[
\Phi\colon P\longrightarrow\mathbb{R}^n,
\qquad
F\colon P\times L(n,d,h)\longrightarrow\mathbb{R}^n
\]
be continuous maps satisfying
\[
\rho_n(F(p,x)-x)\leq R,
\qquad
F(p,x)\neq\Phi(p)
\]
for all $p\in P$ and $x\in L(n,d,h)$.
Then there exists a continuous map
\[
\varphi\colon P\longrightarrow K(n,d,h)
\]
such that
\[
\sup_{p\in P}\rho_n(\varphi(p)-\Phi(p))<R+4h.
\]
\end{restatedapproximation}

\begin{proof}
Set $K=K(n,d,h)$ and $L=L(n,d,h)$.
The construction below implements 
the idea outlined in
Approach~3 of \S\ref{subsection: motivation}.

\smallskip
\noindent
\emph{Step 1: Extending $F$ and choosing a collapse map.}

By Lemma~\ref{lemma: Tietze extension},
we can extend $F$ to a continuous map
$F\colon P\times\mathbb{R}^n
\to\mathbb{R}^n$
such that
\[
\rho_n(F(p,x)-x)\leq R
\qquad
(p\in P,\ x\in\mathbb{R}^n).
\]
We continue to denote this extension by $F$.
For $p\in P$, put
\[
F_p(x)=F(p,x),
\qquad
A_p(x)=F_p(x)-\Phi(p).
\]
Let
$M=\max_{p\in P}\rho_n(\Phi(p))$.
Then
\[ \rho_n(A_p(x)-x)\leq R+M, \]
and consequently
\begin{equation}
\label{eq: approximation properness estimate}
\rho_n(A_p(x))\geq\rho_n(x)-R-M.
\end{equation}

The function $(p,x)\mapsto\rho_n(A_p(x))$
is strictly positive on $P\times L$.
By \eqref{eq: approximation properness estimate},
it tends to infinity as $\rho_n(x)\to\infty$,
uniformly in $p$.
Since $P$ is compact and $L$ is closed,
we can therefore choose
\[
0<\delta<\frac h2
\]
such that
\begin{equation}
\label{eq: approximation uniform avoidance}
\rho_n(A_p(x))>\delta
\qquad
(p\in P,\ x\in L).
\end{equation}
In particular, the set
\[
V=
\{(p,x)\in P\times\mathbb{R}^n:
  \rho_n(A_p(x))\leq\delta\}
\]
is compact and disjoint from $P\times L$.

Throughout this proof, 
let $S^n=\mathbb{R}^n\cup\{\infty\}$ be 
the one-point compactification of $\mathbb{R}^n$, 
with $\infty$ as its basepoint.

Define a basepoint-preserving map
$b\colon S^n\to S^n$ by
\[
b(x)=
\begin{cases}
\displaystyle
\frac{x}{1-\rho_n(x)/\delta},
& \rho_n(x)<\delta,\\[3mm]
\infty,
& \rho_n(x)\geq\delta,
\end{cases}
\qquad
b(\infty)=\infty.
\]
Thus $b$ maps the open $\rho_n$-ball of radius
$\delta$ homeomorphically onto $\mathbb{R}^n$
and collapses its complement to $\infty$.

This map is homotopic to the identity through
basepoint-preserving maps.
Indeed, for $0\leq t\leq1$, define
\[
b_t(x)=
\begin{cases}
\displaystyle
\frac{x}{1-t\rho_n(x)/\delta},
& t\rho_n(x)<\delta,\\[3mm]
\infty,
& t\rho_n(x)\geq\delta,
\end{cases}
\qquad
b_t(\infty)=\infty.
\]
These formulas define a continuous based homotopy
with $b_0=\mathrm{id}_{S^n}$ and $b_1=b$.
For continuity at infinity, note that
$\rho_n(b_t(x)) \geq \rho_n(x)$
whenever
$b_t(x)$ is finite.

In the rest of the proof, all auxiliary spaces
$\mathcal E_n(C,b)$ are defined using this same map $b$.

\smallskip
\noindent
\emph{Step 2: Choosing the local complexes $K_\sigma$.}

We first describe finite cubical neighborhoods
whose connectivity is easy to control.
For $a\in\mathbb{R}^n$ and $s>0$, write
\[
\overline B_\rho(a,s)
=
\{x\in\mathbb{R}^n\mid \rho_n(x-a)\leq s\}.
\]
Let $D(a,s)$ be the union of all closed
$n$-dimensional grid cubes meeting
$\overline B_\rho(a,s)$, together with their faces.
See Figure \ref{fig:cubical-neighborhood} (a).
This is a finite cubical complex, and
\begin{equation}
\label{eq: approximation cubical neighborhood}
\overline B_\rho(a,s)
\subset D(a,s)
\subset \overline B_\rho(a,s+h),
\end{equation}
because every grid cube has $\rho_n$-diameter $h$.

\begin{figure}[tb]
  \centering
  \includegraphics[width=\textwidth]
    {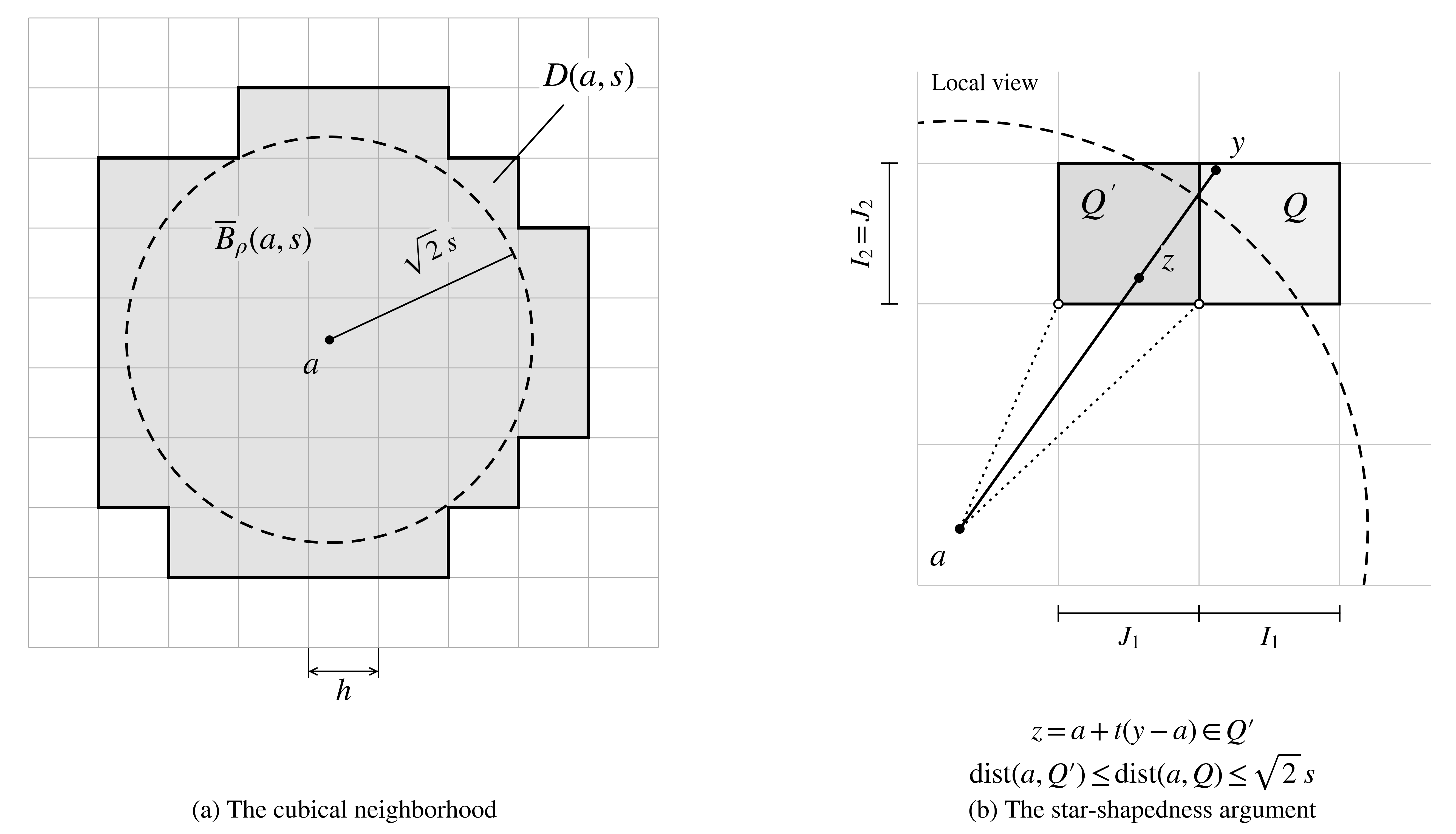}
  \caption{
    The cubical neighborhood $D(a,s)$ and its star-shapedness
    in dimension $n=2$.
    (a) The shaded region is the union of all closed grid squares
    of side length $h$ meeting $\overline B_\rho(a,s)$;
    the dashed circle is the boundary of this ball.
    (b) A local view of the construction: $y\in Q$ and
    $z=a+t(y-a)\in Q'$, where $Q=I_1\times I_2$
    and $Q'=J_1\times J_2$.
    In the illustrated case, $J_2=I_2$, while $J_1$ is closer
    to $a_1$ than $I_1$.
    The dotted segments realize the Euclidean distances
    from $a$ to the two squares.
    Since $\operatorname{dist}(a,Q')\leq\operatorname{dist}(a,Q)
    \leq\sqrt{2}\,s$, the square $Q'$ also meets the ball,
    and hence $z\in D(a,s)$.
  }
  \label{fig:cubical-neighborhood}
\end{figure}

We claim that $D(a,s)$ is star-shaped with respect
to $a$, and hence contractible.\footnote{
A subset $D\subset\mathbb{R}^n$ is said to be
\emph{star-shaped with respect to $a\in D$} if
\[
a+t(y-a)\in D
\qquad (y\in D,\ 0\leq t\leq1).
\]
In other words, the line segment joining $a$
to any point of $D$ is contained in $D$.
}
Let $y\in D(a,s)$, and choose a grid cube
\[
Q=\prod_{i=1}^n I_i
\]
containing $y$ and meeting $\overline B_\rho(a,s)$.
Fix $t\in[0,1]$ and set $z=a+t(y-a)$.

For each $i$, choose a grid interval $J_i$
containing $z_i$ as follows.
If $z_i\in I_i$, take $J_i=I_i$.
Otherwise, $z_i$ lies between $a_i$ and the nearest
endpoint of $I_i$, so
\[
|a_i-z_i|\leq\operatorname{dist}(a_i,I_i).
\]
In this case, choose any grid interval $J_i$
containing $z_i$.
See Figure \ref{fig:cubical-neighborhood} (b).
In either case,
\[
\operatorname{dist}(a_i,J_i)
\leq\operatorname{dist}(a_i,I_i).
\]
Here $\operatorname{dist}$ denotes Euclidean distance.
The cube $Q'=\prod_{i=1}^nJ_i$ therefore contains $z$
and satisfies
\[
\operatorname{dist}(a,Q')^2
=\sum_{i=1}^n\operatorname{dist}(a_i,J_i)^2
\leq\sum_{i=1}^n\operatorname{dist}(a_i,I_i)^2
=\operatorname{dist}(a,Q)^2
\leq ns^2.
\]
Thus $Q'$ also meets $\overline B_\rho(a,s)$,
which proves that $z\in D(a,s)$.
This establishes the claim.

Since $D(a,s)$ is contractible, its $d$-skeleton
with respect to the cubical CW structure is
$(d-1)$-connected;
see \cite[Corollary~4.12]{Hatcher}.
This skeleton is precisely $D(a,s)\cap K$.

Put
\[
\varepsilon=\frac{h}{2(\dim P +1)}.
\]
After subdividing $P$, we may assume that
\begin{equation}
\label{eq: approximation small simplices}
\rho_n(\Phi(p)-\Phi(p'))<\varepsilon
\qquad
(p,p'\in\sigma)
\end{equation}
for every simplex $\sigma$ of $P$.

For each nonempty simplex $\sigma$
of $P$, choose
$p_\sigma\in\sigma$ and set
\[
a_\sigma=\Phi(p_\sigma),
\qquad
s_\sigma=R+2h+(\dim\sigma)\varepsilon.
\]
Define
\[
D_\sigma=D(a_\sigma,s_\sigma),
\qquad
K_\sigma=D_\sigma\cap K.
\]
Each $K_\sigma$ is a nonempty finite
$(d-1)$-connected CW complex.

These complexes are nested along faces.
Indeed, if $\tau$ is a proper face of $\sigma$, then
\[
\rho_n(a_\tau-a_\sigma)<\varepsilon
\leq s_\sigma-s_\tau.
\]
Hence
$\overline B_\rho(a_\tau,s_\tau)
\subset \overline B_\rho(a_\sigma,s_\sigma)$,
which implies
\begin{equation}
\label{eq: approximation nested complexes}
D_\tau\subset D_\sigma,
\qquad
K_\tau\subset K_\sigma.
\end{equation}

Let
\[ r\colon\mathbb{R}^n\setminus L
   \longrightarrow  K \]
be the retraction provided by
Lemma~\ref{lemma: retraction}.
If $p\in\sigma$ and $\rho_n(A_p(x))\leq\delta$,
then $x\notin L$ by
\eqref{eq: approximation uniform avoidance}, and
(recall that $A_p(x) = F_p(x)-\Phi(p)$ and 
$0<\delta<h/2$)
\[
\begin{aligned}
\rho_n(r(x)-a_\sigma)
&\leq
\rho_n(r(x)-x)
+\rho_n(x-F_p(x))\\
&\qquad
+\rho_n(F_p(x)-\Phi(p))
+\rho_n(\Phi(p)-a_\sigma)\\
&\leq \frac h2+R+\delta+\varepsilon\\
&<R+h+\varepsilon
<R+2h
\leq s_\sigma.
\end{aligned}
\]
Since $r(x)\in K$, it follows that
\begin{equation}
\label{eq: approximation labels in local complex}
r(x)\in K_\sigma
\quad\text{whenever}\quad
p\in\sigma,\ \rho_n(A_p(x))\leq\delta.
\end{equation}

Finally, for every $p\in\sigma$ and $y\in K_\sigma$,
\eqref{eq: approximation cubical neighborhood}
and \eqref{eq: approximation small simplices} give
\begin{equation}
\label{eq: approximation local error bound}
\begin{aligned}
\rho_n(y-\Phi(p))
&\leq s_\sigma+h+\varepsilon
\leq R+3h+(\dim P +1)\varepsilon\\
&=R+\frac{7h}{2}.
\end{aligned}
\end{equation}
Thus any map taking $\sigma$ into $K_\sigma$
will satisfy the required approximation estimate.

\smallskip
\noindent
\emph{Step 3: Constructing maps into the auxiliary spaces.}

For each simplex $\sigma$ of $P$
and each $p\in\sigma$,
define a based map
\[
u_{\sigma,p}\colon
S^n\longrightarrow\Sigma^n(K_\sigma)_+
\]
by
\[
u_{\sigma,p}(x)=
\begin{cases}
[b(A_p(x)),r(x)],
& \rho_n(A_p(x))<\delta,\\
*,
& \rho_n(A_p(x))\geq\delta,
\end{cases}
\qquad
u_{\sigma,p}(\infty)=*.
\]
Here $*$ denotes the basepoint of
$\Sigma^n(K_\sigma)_+$.
The first case is well-defined by
\eqref{eq: approximation labels in local complex}.
These maps depend continuously on $p$.\footnote{
For completeness, we provide a proof of continuity. 
The set
\[
V_\sigma=
\{(p,x)\in\sigma\times\mathbb{R}^n:
  \rho_n(A_p(x))\leq\delta\}
\]
is compact and disjoint from $\sigma\times L$.
On $V_\sigma$, the formula
\[
(p,x)\longmapsto[b(A_p(x)),r(x)]
\]
is continuous with values in $\Sigma^n(K_\sigma)_+$.
When $\rho_n(A_p(x))=\delta$, its value is the
basepoint, because $b(A_p(x))=\infty$.
On the closed set
\[
(\sigma\times S^n)\setminus
\{(p,x)\in\sigma\times\mathbb{R}^n:
  \rho_n(A_p(x))<\delta\},
\]
the map is constant at the basepoint.
The two formulas agree on the overlap.
Therefore we obtain the continuity of
\[
\sigma\times S^n\longrightarrow\Sigma^n(K_\sigma)_+,
\qquad
(p,x)\longmapsto u_{\sigma,p}(x).
\]
By the compact-open topology, this gives a continuous
map $p\mapsto u_{\sigma,p}$ into
$\Omega^n\Sigma^n(K_\sigma)_+$.}

For $p\in P$, define a path
\[
\gamma_p\colon I\longrightarrow\Omega^n S^n
\]
by
\[
\gamma_p(t)(x)
=
b\bigl((1-t)A_p(x)+tx\bigr)
\qquad(x\in\mathbb{R}^n),
\]
and set $\gamma_p(t)(\infty)=\infty$.
The expression inside $b$ satisfies
\[
\rho_n\bigl((1-t)A_p(x)+tx-x\bigr)
\leq R+M.
\]
Consequently, for all sufficiently large $\rho_n(x)$,
the value of $\gamma_p(t)(x)$ is $\infty$,
uniformly in $(p,t)\in P\times I$.
This proves continuity at infinity.
Thus $(p,t,x)\mapsto\gamma_p(t)(x)$ is continuous,
and so is $p\mapsto\gamma_p$.

The endpoints satisfy
\[
\gamma_p(0)=q_{K_\sigma}\circ u_{\sigma,p},
\qquad
\gamma_p(1)=b.
\]
We therefore obtain a continuous map
\[
\chi_\sigma\colon
\sigma\longrightarrow\mathcal E_n(K_\sigma,b),
\qquad
\chi_\sigma(p)=(u_{\sigma,p},\gamma_p).
\]
The map $u_{\sigma,p}$ records the approximate
solutions of $F_p(x)=\Phi(p)$ together with their
labels $r(x)$
(cf. Approach 2 of 
\S \ref{subsection: motivation}).
The path $\gamma_p$ compares the map obtained by
forgetting these labels with the fixed map $b$.

For a face $\tau\subset\sigma$, let
\[
i_{\tau\sigma}\colon K_\tau\hookrightarrow K_\sigma
\]
denote the inclusion.
Since the formulas above use the same maps $b$, $r$,
and $\gamma_p$ for every simplex, we have the 
compatibility relation
\begin{equation}
\label{eq: approximation compatible auxiliary maps}
\chi_\sigma|_\tau
=
\mathcal E_n(i_{\tau\sigma})\circ \chi_\tau.
\end{equation}

\smallskip
\noindent
\emph{Step 4: Constructing the approximation 
over the skeleta.}

Recall $1\leq d <n$.
Since $K_\sigma$ is $(d-1)$-connected,
Lemma~\ref{lemma: connectivity of canonical map},
applied with $k=d-1\geq 0$, shows that
\[
\iota_{K_\sigma}\colon
K_\sigma\longrightarrow\mathcal E_n(K_\sigma,b)
\]
is $(2d-1)$-connected.

We construct, 
for every simplex $\sigma$ of $P$,
continuous maps
\[
\varphi_\sigma\colon\sigma\longrightarrow K_\sigma,
\qquad
H_\sigma\colon
\sigma\times I\longrightarrow\mathcal E_n(K_\sigma,b)
\]
such that
\begin{equation}
\label{eq: approximation local homotopy}
H_\sigma(p,0)
=\iota_{K_\sigma}(\varphi_\sigma(p)),
\qquad
H_\sigma(p,1)=\chi_\sigma(p),
\end{equation}
and such that, for every face $\tau\subset\sigma$,
\begin{equation}
\label{eq: approximation face compatibility}
\varphi_\sigma|_\tau
=i_{\tau\sigma}\circ\varphi_\tau,
\qquad
H_\sigma|_{\tau\times I}
=\mathcal E_n(i_{\tau\sigma})\circ H_\tau.
\end{equation}
We note that only the maps $\varphi_\sigma$ are needed
for the final approximation.
However, to carry out the inductive construction,
we must construct the homotopies $H_\sigma$
simultaneously.
These homotopies provide the boundary data needed
to apply Lemma~\ref{lemma: relative extension}
when extending the construction to
higher-dimensional simplices.

We begin with the vertices.
For a vertex $v$ of $P$, 
the map $\iota_{K_v}$ is
surjective on path components.
Hence we can choose $\varphi_v(v)\in K_v$
and a path
\[
H_v\colon\{v\}\times I\longrightarrow
\mathcal E_n(K_v,b)
\]
from $\iota_{K_v}(\varphi_v(v))$ 
to $\chi_v(v)$.
This gives the required construction on the
$0$-skeleton.

We next explain explicitly the extension over an edge
$\sigma=[v_0,v_1]$ of $P$.
Using the inclusions $K_{v_i}\subset K_\sigma$,
the two chosen vertex values define a map
\[
\varphi_{\partial\sigma}\colon
\partial\sigma\longrightarrow K_\sigma.
\]
Likewise, the two paths $H_{v_i}$ give a map
\[
H_{\partial\sigma}\colon
\partial\sigma\times I
\longrightarrow\mathcal E_n(K_\sigma,b)
\]
by
\[
H_{\partial\sigma}(v_i,t)
=
\mathcal E_n(i_{v_i\sigma})(H_{v_i}(v_i,t)).
\]
The naturality relation
\eqref{eq: functoriality of canonical map}
and \eqref{eq: approximation compatible auxiliary maps}
imply
\[
H_{\partial\sigma}(v_i,0)
=
\iota_{K_\sigma}(\varphi_{\partial\sigma}(v_i)),
\qquad
H_{\partial\sigma}(v_i,1)=\chi_\sigma(v_i).
\]
We can therefore apply
Lemma~\ref{lemma: relative extension}
with $j=1$, 
$f=\iota_{K_\sigma} \colon K_\sigma 
\to \mathcal{E}_n(K_\sigma)$
and $\chi=\chi_\sigma$.
It extends both $\varphi_{\partial\sigma}$
and $H_{\partial\sigma}$ over $\sigma$,
giving $\varphi_\sigma$ and $H_\sigma$
with the required properties.

The higher-dimensional steps are analogous.
Suppose that the construction has been completed
on all simplices of dimension less than $j$,
and let $\sigma$ be a $j$-simplex.
For $p\in\tau\subset\partial\sigma$, define
\[
\varphi_{\partial\sigma}(p)
=i_{\tau\sigma}(\varphi_\tau(p)),
\qquad
H_{\partial\sigma}(p,t)
=
\mathcal E_n(i_{\tau\sigma})(H_\tau(p,t)).
\]
The inductive compatibility conditions ensure that
these formulas agree on intersections of faces.
They therefore define continuous maps
\[
\varphi_{\partial\sigma}\colon
\partial\sigma\longrightarrow K_\sigma,
\qquad
H_{\partial\sigma}\colon
\partial\sigma\times I
\longrightarrow\mathcal E_n(K_\sigma,b).
\]
Again, naturality and
\eqref{eq: approximation compatible auxiliary maps}
give
\[
H_{\partial\sigma}(p,0)
=\iota_{K_\sigma}(\varphi_{\partial\sigma}(p)),
\qquad
H_{\partial\sigma}(p,1)=\chi_\sigma(p).
\]

Since
\[
j\leq\dim P\leq2d-1
\]
and $\iota_{K_\sigma}$ is $(2d-1)$-connected,
Lemma~\ref{lemma: relative extension}, applied to
$(\sigma,\partial\sigma)\cong(D^j,\partial D^j)$,
extends these maps to $\varphi_\sigma$ and $H_\sigma$
satisfying
\eqref{eq: approximation local homotopy}
and \eqref{eq: approximation face compatibility}.
This completes the induction.

\smallskip
\noindent
\emph{Step 5: Gluing and estimating the error.}

The maps $\varphi_\sigma$ agree on common faces
when regarded as maps into $K$.
Since $P$ is a finite simplicial complex,
they glue to a continuous map
\[
\varphi\colon P\longrightarrow K.
\]
For every simplex $\sigma$,
\[
\varphi(\sigma)\subset K_\sigma.
\]
Applying \eqref{eq: approximation local error bound}
with $y=\varphi(p)$ $(p\in \sigma)$, we obtain
\[
\sup_{p\in P}\rho_n(\varphi(p)-\Phi(p))
\leq R+\frac{7h}{2}
<R+4h.
\]
This proves the proposition.
\end{proof}


%
%

\end{document}